\documentclass[preprint,11pt,numbers,sort&compress]{elsarticle}

\usepackage{soul}
\usepackage{psfrag}
\usepackage{a4wide}
\usepackage{rotating}
\usepackage{color}
\usepackage{amssymb}
\usepackage{amsmath}
\usepackage{amsthm}
\usepackage{verbatim}
\usepackage{cancel}
\usepackage{nicematrix}
\usepackage{multirow}
\usepackage{tikz}

\newcommand{\f}[1]{\mathbf{#1}}
\newcommand{\ab}[1]{\boldsymbol{#1}}
\def\bfm#1{\boldsymbol{#1}}
\newcommand{\bb}[1]{\bfm{#1}}
\newcommand{\N}{\mathbb N}

\newcommand{\R}{\mathbb R}

\newcommand{\V}{\mathcal{V}}

\newcommand{\W}{\mathcal{W}}

\newcommand{\LL}{i_0}
\newcommand{\RR}{i_1}
\newcommand{\g}{\gamma}
\newcommand{\gC}{\gamma}

\newcommand{\Side}{\tau}

\newcommand{\JH}{\mathcal{J}_{\mathcal{H}}^{(i)}}
\newcommand{\JPH}{\mathcal{J}_{\partial,\mathcal{H}}^{(i)}}

\newcommand{\MM}{\ab{M}_B}

\DeclareMathOperator{\Rank}{rank}
\DeclareMathOperator{\Diag}{diag}

\DeclareMathOperator{\Span}{span}

\theoremstyle{definition}
\newtheorem{ex}{Example}

\newtheorem{rem}{Remark}

\newproof{pf}{proof}
\advance\textheight by 0.4cm
\advance\topmargin by -0.2cm

\definecolor{gold}{rgb}{1,0.7,0}
\definecolor{dred}{rgb}{0.92,0,0}
\definecolor{dgreen}{rgb}{0,0.6,0}
\def\MK{\color{dgreen}}

\def\new{\color{black}}

\begin{document}

\begin{frontmatter}

\title{
{\new An adaptive Dual-Primal Isogeometric Tearing and Interconnecting (IETI-DP) method for multi-patch structured Kirchhoff--Love plates}
}

\cortext[cor]{Corresponding author}

\author[vil]{Mario Kapl}
\ead{m.kapl@fh-kaernten.at}

\author[slo1]{Alja\v z Kosma\v c}
\ead{aljaz.kosmac@iam.upr.si}

\author[slo1]{Vito Vitrih\corref{cor}}
\ead{vito.vitrih@upr.si}

\address[vil]{Department of Engineering $\&$ IT, Carinthia University of Applied Sciences, Villach, Austria}

\address[slo1]{IAM and FAMNIT, University of Primorska, Koper, Slovenia}

\begin{abstract}
We present a novel adaptive isogeometric method for solving the 
{\new Kirchhoff--Love plate problem for plates whose mid-surface is defined by a} 
planar multi-patch domain with possibly extraordinary vertices, parametrized by an analysis-suitable $G^1$ multi-patch geometry. The proposed technique relies on the concept of Dual-Primal Isogeometric Tearing and Interconnecting (IETI-DP), which enforces the required $C^1$-smoothness of the solution across a common edge of two neighboring patches by imposing appropriate continuity conditions by means of Lagrange multipliers. The resulting saddle point problem is solved using a dual-primal formulation, first by a small linear problem for the Lagrange multipliers and then by local, parallelizable linear problems on the individual patches for the coefficients of the numerical solution. While for the local problems on the single patches standard diagonally scaling is used as preconditioner, a suitable preconditioner for the problem of finding the Lagrange multipliers is introduced. To perform adaptive refinement, the solution of the 
{\new Kirchhoff--Love plate problem on} the single patches of the multi-patch domain is discretized by employing (truncated) hierarchical B-splines, and an appropriate refinement strategy of the underlying mesh is presented. Finally, the potential of the developed adaptive IETI-DP method for solving the 
{\new Kirchhoff--Love plate problem} on planar multi-patch geometries is numerically tested on the basis of several numerical examples. Thereby, the numerical results show on the one hand optimal convergence behavior with respect to adaptive refinement, and on the other hand a good performance of the proposed preconditioner for the linear problem of determining the Lagrange multipliers.
\end{abstract}

\end{frontmatter}

\section{Introduction}

The concept of Isogeometric Tearing and Interconnecting (IETI) was introduced in~\cite{KlPeSaJu12} in the field of isogeometric analysis~\cite{HuCoBa04} as a method for solving second-order linear partial differential equations (PDEs) {\new on} planar multi-patch geometries with possibly extraordinary vertices, i.e.~{\new on} planar unstructured quadrilateral meshes. The basic idea is to solve the given PDE by means of smaller local and parallelizable problems on the individual patches. To ensure the required smoothness of the solution of the PDE across a common edge of two neighboring patches, appropriate continuity conditions are imposed using Lagrange multipliers. This is achieved by transforming the original problem into a dual-primal formulation in order to first determine the Lagrange multipliers, followed by solving the local problems to get the coefficients of the numerical solution of the PDE. Due to the use of a dual-primal (DP) formulation, the developed IETI methods, e.g~\cite{KaKoVi26,KlPeSaJu12,SoglTakacs_IETI_Elasticity,IETI_LowRank2024,HoLa17,SoTa23,MoSaSchTaTa23,Ta25,WiZaScPa21}, are also referred to as IETI-DP methods.

The described strategy of the IETI-DP techniques for solving PDEs originates from the concept of Finite Element Tearing and Interconnecting (FETI), e.g.~\cite{Farhat1991,FarhatFETI-DP,Klawonn2000,Klawonn2006}, and was introduced as a domain decomposition solver for PDEs. While most of the developed FETI techniques, e.g.~\cite{,FarhatFETI-DP,Klawonn2000,Klawonn2006}, are based on a dual-primal formulation and are also referred to as FETI-DP methods, the original work~\cite{Farhat1991} only uses a dual formulation instead of a dual-primal formulation and therefore requires the use of a pseudoinverse to solve the local problems on the individual patches. A closely related approach to FETI-DP is Balancing Domain Decomposition by constraints (BDDC), e.g. \cite{Do03,MaDo03} in the context of finite element analysis, and e.g.~\cite{VePaScWiZa14,VePaScWiZa17,WiZaScPa21} in the context of isogeometric analysis. From a theoretical point of view, the FETI-DP and BDDC methods are equivalent, cf.~\cite{MaDoTe05}, but they apply two different approaches to domain decomposition. While the FETI-DP (and IETI-DP) techniques enforce in general just the minimal required smoothness across the inner edges of the multi-patch domain, such as $C^0$ continuity in case of second-order PDEs and $C^1$ continuity in case of fourth-order problems, the BDDC methods were also employed to enforce the maximal possible smoothness, i.e. $C^{p-1}$ continuity, across the inner edges (cf.~\cite{VePaScWiZa14,VePaScWiZa17}), where $p$ is the degree of the used discretization space. Another Lagrange multiplier-based method is the mortar method, e.g. \cite{BeLoSaTa23,Bouclier2017,DiSchWoHe19,SchDiWoKlHe19,MiZoScBoTh21}, where the smoothness across the inner edges of the multi-patch domain is not strongly enforced as in the FETI-DP, IETI-DP and BDDC methods, but only in a weak sense. However, this leads to a numerical solution that is only approximately smooth.

FETI-DP, IETI-DP and BDDC have in common that most of the developed methods are mainly limited to structured meshes or, in the case of unstructured meshes, to second-order PDEs, where the coupling between the patches needs to be only $C^0$-smooth. Two examples of these methods for solving a fourth-order PDE, namely the biharmonic equation, {\new on} planar structured meshes are the FETI-DP technique~\cite{FaMa98} and the IETI-DP framework~\cite{Ta25}. To the best of our knowledge, the only multiplier-based method, which can also solve high-order PDEs {\new on} (planar) unstructured meshes by strongly enforcing the coupling across neighboring patches, i.e.~by computing a numerical solution which is exactly smooth with the desired order of smoothness across the inner edges, is the IETI-DP method~\cite{KaKoVi26}. There, a IETI-DP technique is presented, which allows to solve the polyharmonic equation of any order~${\new m} \geq 1$ {\new on} multi-patch domains with possibly extraordinary vertices, parametrized by bilinear-like $G^s$ multi-patch geometries~\cite{KaVi20b} with $s \geq {\new m-1}$.

In general, solving higher-order PDEs, i.e.~PDEs of order~$2m$, $m \geq 2$, {\new on} planar multi-patch domains with possibly extraordinary vertices by strongly enforcing the required $C^s$ smoothness, $s \geq m-1$, of the computed solution across adjacent patches is a challenging task. Until now, the most common approach to solving this issue has been to develop methods that employ a globally $C^s$-smooth basis to discretize the numerical solution of the given PDE. Thereby, the constructed techniques can be classified according to the type of parametrization used to represent the multi-patch domain. While the case $s=1$ comprises the use of many different multi-patch parameterizations such as $C^1$-smooth multi-patch parameterizations with singularities at the extraordinary vertices, e.g.~\cite{NgPe16,ToSpHu17,WeLiQiHuZhCa22}, $C^1$-smooth multi-patch parameterizations with $G^1$-smooth caps in the vicinity of the extraordinary vertices, e.g.~\cite{KaPe17,KaPe18,WeFaLiWeCa23}, particular class of $G^1$-smooth multi-patch parameterizations called analysis-suitable $G^1$ multi-patch parameterizations, e.g.~\cite{CoSaTa16,FaJuKaTa22,KaSaTa17b}, general $G^1$-smooth multi-patch parameterizations, e.g.~\cite{ChAnRa18,ChAnRa19,MaMaMo2024,mourrain2015geometrically} and $G^1$-smooth multi-patch parameterizations with singularities at so-called scaling centers, e.g.~\cite{ArReKlSi23,ReArSiKl23}, the case $s >1$ is so far mainly limited to the class of bilinear-like $G^s$ multi-patch geometries, e.g.~\cite{KaVi17b,KaVi17c,KaVi19a,KaVi20} for $s=2$ and \cite{KaVi20b,KaKoVi24b,KaKoVi24c,KaKoVi24} for an arbitrary~$s \geq 1$. The recently developed IETI-DP method~\cite{KaKoVi26} represents a promising alternative which, unlike the techniques described above, does not require a globally $C^s$-smooth basis, whose construction can be technically complex and time-consuming, but works only with the ``standard'' isogeometric spline basis functions on the individual patches. Moreover, the IETI-DP framework~\cite{KaKoVi26} directly enables the parallel solving of the PDE by decomposing the large problem into small and local problems. 

Besides solving PDEs by uniform refinement, solving them using adaptive refinement, i.e. by refining locally via an error estimator, has become popular in isogeometric analysis. The idea is to refine only those parts of the underlying mesh where the error is high. This allows to significantly reduce the needed degrees of freedom to represent a solution with the same accuracy compared to uniform refinement. For more details about the general concept of adaptive isogeometric analysis, we refer to the review article~\cite{BuGaGiPrVa22}. Anyway, only a few methods exist that solve high-order PDEs {\new on} multi-patch geometries with extraordinary vertices using adaptive refinement with exactly $C^s$-smooth isogeometric spline spaces, in particular with $C^1$-smooth spline spaces. Examples are the T-spline based frameworks~\cite{ScSiEv13,WeZhLiHu17}, the techniques~\cite{CaWeToLiHuKiZh20, We22} using the concept of D-patches~\cite{Re97}, {\new see also~\cite{ToSpHu17} for their use in isogeometric analysis}, and the methods~\cite{BrGiKaVa20, BrGiKaVa23, BrFaGiKaVa24,BrGiReToVa23} which generalize the construction of the $C^1$-smooth isogeometric spline spaces~\cite{KaSaTa19a,FaJuKaTa22} defined {\new on} analysis-suitable $G^1$ multi-patch parameterizations to a hierarchical and adaptive setting. {\new While, the techniques~\cite{ScSiEv13,WeZhLiHu17} employ $C^1$-smooth functions to solve still second order PDEs, the methods~\cite{BrGiKaVa20, BrGiKaVa23, BrFaGiKaVa24,BrGiReToVa23,CaWeToLiHuKiZh20,We22} solve fourth order PDEs, namely the biharmonic equation~\cite{BrGiKaVa20,BrGiKaVa23,We22}, the Kirchhoff-Love shell problem~\cite{BrFaGiKaVa24,CaWeToLiHuKiZh20}, and the Cahn-Hilliard equation~\cite{BrGiReToVa23}.} 

The goal of this paper is to demonstrate that IETI-DP is a suitable framework for solving high-order PDEs {\new on} multi-patch domains with extraordinary vertices via adaptive refinement. More specifically, we will not only show that the proposed IETI-DP technique is an innovative and practical approach to tackling this challenging task, but will also make it possible to directly solve the problem efficiently by decomposing it into smaller local and parallelizable problems. For this purpose, we will develop an adaptive IETI-DP method to determine the solution of a particular fourth-order PDE, given by the 
{\new Kirchhoff--Love plate problem}, {\new on} planar analysis-suitable $G^1$ multi-patch geometries. Thereby, the presented approach will rely on the recently developed IETI-DP method~\cite{KaKoVi26} {\new for the biharmonic equation}, and will enhance it in 
{\new four} directions. {\new First, we will generalize the technique~\cite{KaKoVi26} so that instead of the biharmonic equation with clamped boundary conditions, also the more general Kirchhoff--Love plate problem with clamped, simply supported, and free boundary conditions, including arbitrary combinations of these boundary types, can be handled. Second}, we will {\new significantly simplify the underlying IETI-DP algorithm of the technique~\cite{KaKoVi26} by using} another formulation of the corresponding weak form of the problem. 
{\new Third}, we will provide a suitable Dirichlet  preconditioner for the linear problem of determining the Lagrange multipliers, which has not been previously available. 
{\new Fourth,} as the main contribution of this work, we will generalize the method~\cite{KaKoVi26}, which has been limited so far to uniform refinement, also to the case of adaptive refinement.  

An important feature of the proposed adaptive IETI-DP method will be that the underlying $C^1$-smooth hierarchical discretization space will coincide with the global $C^1$-smooth hierarchical spline space {\new on planar analysis-suitable $G^1$ multi-patch geometries} constructed in~\cite{BrGiKaVa23} {\new (and extended to the case of multi-patch structured Kirchoff-Love shells in~\cite{BrFaGiKaVa24})} in most parts of the multi-patch domain except in the vicinity of an inner vertex. The reason for this small deviation will be that the adaptive method~\cite{BrGiKaVa23} may require the refinement of numerous additional elements of the underlying mesh near an inner vertex, whereas in the adaptive method presented in this work the number of these elements will be smaller in the case of hierarchical B-splines, or zero in the case of truncated hierarchical B-splines. Thus, the numerical solution computed by the proposed adaptive IETI-DP method will convergence with the same optimal convergence behavior as in~\cite{BrGiKaVa23}. However, in contrast to~\cite{BrGiKaVa23}, the explicit construction of the global $C^1$-smooth hierarchical spline space, which is based on the design and use of a quite technical refinement mask for the $C^1$-smooth one-level
basis functions, will not be needed. Instead, just the ``standard'' (truncated) hierarchical B-splines (e.g.~\cite{Kr97,GiJuSp2012}) with their relatively simple refinement mask will be employed to represent the numerical solution on the single patches of the multi-patch domain, and the required $C^1$-smoothness across the inner edges will be enforced by means of Lagrange multipliers. 
Several numerical examples will demonstrate the practical applicability of the new adaptive method for the solving of the 
{\new Kirchhoff--Love plate problem} {\new on} analysis-suitable $G^1$ multi-patch geometries, and will numerically verify its optimal convergence behavior with respect to adaptive refinement as well as the good performance of the proposed preconditioner for the linear problem of computing the Lagrange multipliers.

The remainder of the paper is organized as follows. Section~\ref{sec:prel} will provide the necessary tools and concepts, in particular the setting of the considered multi-patch geometries called analysis-suitable $G^1$ multi-patch parameterizations, the basic idea of (truncated) hierarchical B-splines as well as the concept of $C^1$-smooth isogeometric hierarchical spline functions {\new on} the given multi-patch domains. In Section~\ref{sec:IETI}, we will present the IETI-DP framework that allows to solve the 
{\new Kirchhoff--Love plate problem} {\new on} analysis-suitable $G^1$ multi-patch geometries via adaptive refinement. The implementation details of the developed adaptive IETI-DP method, such as the refinement procedure of the underlying mesh and the computation of the constraint matrices to couple the numerical solution with $C^1$-smoothness across the inner edges, will be described in Section~\ref{sec:implementation_details}. {\new The details regarding the computational realization of the adaptive IETI-DP method will be given in~\ref{app:computational-realization}}. Section~\ref{sec:numerical_examples} will show several numerical examples to demonstrate the potential of the proposed adaptive isogeometric method for solving the 
{\new Kirchhoff--Love plate problem}. Finally, we conclude the paper in Section~\ref{sec:Conclusion}. 

\section{Preliminaries} \label{sec:prel}

We introduce some required preliminaries including the planar multi-patch domains used throughout the paper {\new to define the mid-surface of the plate}, the specific class of multi-patch geometries, called analysis-suitable $G^1$ multi-patch parameterizations~\cite{CoSaTa16}, to represent these multi-patch domains, the concept of (truncated) hierarchical B-splines as well as the idea of $C^1$-smooth isogeometric hierarchical spline functions {\new on} the considered multi-patch domains. 
\subsection{The planar multi-patch domain} \label{subsec:planar_multi-patch_domain}

We consider 
an open domain $\Omega \subset \R^2$ whose closure~$\overline{\Omega}$ 
is a planar multi-patch domain given as the disjoint union
\[
\overline{\Omega}
=
\bigcup_{i \in \mathcal{I}_{\Omega}} \Omega^{(i)}
\;\dot{\cup}\;
\bigcup_{i \in \mathcal{I}_{\Gamma}} \Gamma^{(i)}
\;\dot{\cup}\;
\bigcup_{i \in \mathcal{I}_{\Xi}} \bfm{\Xi}^{(i)}
\]
with open quadrilateral patches~$\Omega^{(i)}$, $i \in \mathcal{I}_{\Omega}$, open edges~$\Gamma^{(i)}$, $i \in \mathcal{I}_{\Gamma}$, and vertices~$\bfm{\Xi}^{(i)}$, $i \in \mathcal{I}_{\Xi}$. The index sets $\mathcal{I}_{\Omega}$, $\mathcal{I}_{\Gamma}$ and $\mathcal{I}_{\Xi}$ are further divided into $\mathcal{I}_{\Omega} = \mathcal{I}_{\Omega}^I \dot{\cup} \mathcal{I}_{\Omega}^B$, $\mathcal{I}_{\Gamma} = \mathcal{I}_{\Gamma}^I \dot{\cup} \mathcal{I}_{\Gamma}^B$ and $\mathcal{I}_{\Xi} = \mathcal{I}_{\Xi}^I \dot{\cup} \mathcal{I}_{\Xi}^B$, respectively, where $\mathcal{I}_{\Omega}^I$, $\mathcal{I}_{\Gamma}^I$ and $\mathcal{I}_{\Xi}^{I}$ collect the indices of the inner patches, edges and vertices, and $\mathcal{I}_{\Omega}^B$, $\mathcal{I}_{\Gamma}^B$ and $\mathcal{I}_{\Xi}^{B}$ collect the indices of the boundary patches, edges and vertices. 
For any two distinct patches~$\Omega^{(i_0)}$ and $\Omega^{(i_1)}$, $i_0, i_1 \in \mathcal{I}_{\Omega}$, $i_0 \neq i_1$, the intersection of their closures, i.e.
$
\overline{\Omega^{(i_0)}} \cap \overline{\Omega^{(i_1)}},
$
is either empty, a vertex~$\bfm{\Xi}^{(i)}$ for some $i \in \mathcal{I}_{\Xi}$, or the closure~$\overline{\Gamma^{(i)}}$ of an inner edge~$\Gamma^{(i)}$ for some $i \in \mathcal{I}_{\Gamma}^{I}$.
We denote by $\nu_i$ the patch valency of a vertex~$\bfm{\Xi}^{(i)}$, $i \in \mathcal{I}_{\Xi}$, and  define $\mathcal{I}_{\Xi}^{B,\nu \geq 2} \subset \mathcal{I}_{\Xi}^{B}$ as the {set of indices}
of boundary vertices~$\bfm{\Xi}^{(i)}$, $i \in \mathcal{I}_{\Xi}^{B}$, 
with patch valency~$\nu_i \geq 2$.

Let $\mathcal{S}_h^{p,r}([0,1])$ be the univariate spline space of degree~$p$, regularity~$r$, and mesh size $h = \frac{1}{k+1}$ with $k$ different, uniformly distributed inner knots defined on the unit interval~$[0,1]$, and let $\mathcal{S}_{h}^{\ab{p},\ab{r}}([0,1]^2)$ with $\ab{p} = (p,p)$ and $\ab{r} = (r,r)$ be the tensor-product spline space
$
\mathcal{S}_h^{p,r}([0,1]) \otimes \mathcal{S}_h^{p,r}([0,1])
$
on the unit square~$[0,1]^2$. The B-splines of the spline spaces $\mathcal{S}_h^{p,r}([0,1])$ and $\mathcal{S}_{h}^{\ab{p},\ab{r}}([0,1]^2)$ are denoted by $N_{h,j}^{p,r}$, $j \in \mathcal{J}_h = \{0,1,\ldots,n_h-1\}$, and
$
N_{h,\ab{j}}^{\ab{p},\ab{r}} = N_{h,j_1}^{p,r} N_{h,j_2}^{p,r}
$, $\ab{j} = (j_1,j_2) \in \bb{\mathcal{J}}_h = \mathcal{J}_h \times \mathcal{J}_h$, respectively, 
with
\[
n_h = \dim \mathcal{S}_h^{p,r}([0,1]) = p+1+k(p-r) = p+1+\frac{1-h}{h}(p-r).
\] 
We assume that each patch $\overline{\Omega^{(i)}}$, $i \in \mathcal{I}_{\Omega}$, is 
defined by a bijective and regular geometry mapping 
$\ab{G}^{(i)}: [0,1]^2 \rightarrow \R^2$ with $\ab{G}^{(i)} \in \mathcal{S}_{h}^{\ab{p},\ab{r}}([0,1]^2) \times \mathcal{S}_{h}^{\ab{p},\ab{r}}([0,1]^2)$, and denote by $\ab{G}$ the resulting multi-patch parameterization composed of the single patch parameterizations~$\ab{G}^{(i)}$, $i \in \mathcal{I}_{\Omega}$. The schematic illustration of a planar three patch domain~$\overline{\Omega}$ with its geometry mappings~$\ab{G}^{(i)}$, $i \in \{0,1,2\}$, is visualized in Fig.~\ref{fig:three_patch_isogeometric_analysis}.

\begin{figure}
\begin{center}
\begin{tabular}{c}
\resizebox{0.9\textwidth}{!}{
 \begin{tikzpicture}

\definecolor{patchcolor}{HTML}{8B5CF6}      
\definecolor{vertexcolor}{HTML}{EF4444}    
\definecolor{edgecolor}{HTML}{22D3EE}    
\definecolor{labelcolor}{HTML}{EF4444}    

  \draw[-stealth,line width=0.1mm] (1.95,0.15) .. controls (1.35,0.45) ..(0.6, 0.3);
  \node[scale=0.3] at (1.4,0.49) {$\ab{G}^{(0)}$};  
    \draw[-stealth,line width=0.1mm] (-1.95,0.35) .. controls (-1.2,0.8) .. (-0.6, 0.3);
  \node[scale=0.3] at (-1.2,0.78) {$\ab{G}^{(1)}$};
    \draw[-stealth,line width=0.1mm] (-1.3,-1.25) .. controls (-0.8,-1.1) .. (-0.15, -0.7);
  \node[scale=0.3] at (-0.8,-1.18) {$\ab{G}^{(2)}$};
 
 \begin{scope}[cm={1., 0.10, -0.26, 0.9,(0.02,-0.04)}]
  \coordinate(A) at (0,0); \coordinate(B) at (1.8,2.5); \coordinate(C) at (0,4.2); \coordinate(D) at (-1.8,2.5);  \coordinate(E) at (-3.8,2.4); 
  \coordinate(F) at (-2.6,-0.4); \coordinate(G) at (-3,-2.9); \coordinate(H) at (-0.5,-2.4); \coordinate(I) at (2,-3.1); 
  \coordinate(J) at (2.5,-0.9); \coordinate(K) at (3.9,1.5);
  \draw[line width=0.12mm, patchcolor] (0.0830276, -0.0202994) .. controls (0.00843625, 0.447714) .. (-0.0283338, 0.928053);
  \draw[line width=0.12mm,patchcolor] (0.0830276, -0.0202994) .. controls (0.469181, -0.200027) .. (0.858765, -0.354949);
  \draw[line width=0.12mm,patchcolor] (0.895323, 0.593403) .. controls (0.889936, 0.109089) .. (0.858765, -0.354949);
  \draw[line width=0.12mm,patchcolor] (0.895323, 0.593403) .. controls (0.447681, 0.724473) .. (-0.0283338, 0.928053);
  \node[scale=0.3,patchcolor] at (0.45,0.3) {$\Omega^{(0)}$};

  \draw[line width=0.12mm,patchcolor] (-0.0463954, -0.0178091) .. controls (-0.464553, -0.232213) .. (-0.895251, -0.507768);
  \draw[line width=0.12mm,patchcolor] (-0.0463954, -0.0178091) .. controls (-0.129422, 0.445037) .. (-0.174048, 0.926239);
  \draw[line width=0.12mm,patchcolor] (-1.03215, 0.41779) .. controls (-0.590866, 0.670787) .. (-0.174048, 0.926239);
  \draw[line width=0.12mm,patchcolor] (-1.03215, 0.41779) .. controls (-0.94911, -0.00646313) .. (-0.895251, -0.507768);
  \node[scale=0.3,patchcolor] at (-0.5,0.2) {$\Omega^{(1)}$};

  \draw[line width=0.12mm,patchcolor] (0.0309216, -0.126322) .. controls (0.402589, -0.321429) .. (0.777823, -0.483336);
  \draw[line width=0.12mm,patchcolor] (0.0309216, -0.126322) .. controls (-0.392534, -0.352443) .. (-0.827679, -0.63994);
  \draw[line width=0.12mm,patchcolor] (-0.0807774, -1.23732) .. controls (-0.452036, -0.923429) .. (-0.827679, -0.63994);
  \draw[line width=0.12mm,patchcolor] (-0.0807774, -1.23732) .. controls (0.338466, -0.841568) .. (0.777823, -0.483336);
  \node[scale=0.3,patchcolor] at (-0.1,-0.6) {$\Omega^{(2)}$};

\draw[thick, vertexcolor, fill=vertexcolor] (0.0375, -0.0625) circle (0.14mm);

\draw[thick, vertexcolor, fill=vertexcolor] (-0.1125, 1.05) circle (0.14mm);

\draw[thick, vertexcolor, fill=vertexcolor] (0.9125, -0.45) circle (0.14mm);

\draw[thick, vertexcolor, fill=vertexcolor] (0.9625, 0.6625) circle (0.14mm);

\draw[thick, vertexcolor, fill=vertexcolor] (-1.125, 0.4375) circle (0.14mm);

\draw[thick, vertexcolor, fill=vertexcolor] (-0.9625, -0.65) circle (0.14mm);

\draw[thick, vertexcolor, fill=vertexcolor] (-0.0875, -1.3625) circle (0.14mm);

\draw[line width=0.25mm, edgecolor] (0.02525, -0.00984375) .. controls (-0.0625, 0.478125) .. (-0.10975, 0.991406);

\draw[line width=0.25mm, edgecolor] (0.08125, -0.0860312) .. controls (0.475, -0.278125) .. (0.86875, -0.434781);

\draw[line width=0.25mm, edgecolor] (0.96475, 0.602719) .. controls (0.9625, 0.084375) .. (0.91975, -0.398531);

\draw[line width=0.25mm, edgecolor] (0.9135, 0.670594) .. controls (0.45, 0.796875) .. (-0.054, 1.01934);

\draw[line width=0.25mm, edgecolor] (-0.010125, -0.0829688) .. controls (-0.45, -0.309375) .. (-0.910125, -0.611719);

\draw[line width=0.25mm, edgecolor] (-1.07022, 0.467531) .. controls (-0.596875, 0.740625) .. (-0.158969, 1.01878);

\draw[line width=0.25mm, edgecolor] (-1.11153, 0.394406) .. controls (-1.01563, -0.046875) .. (-0.965281, -0.584344);

\draw[line width=0.25mm, edgecolor] (-0.130063, -1.32153) .. controls (-0.51875, -0.978125) .. (-0.917563, -0.680281);

\draw[line width=0.25mm, edgecolor] (-0.039875, -1.31153) .. controls (0.4, -0.878125) .. (0.860125, -0.490281);

\node[scale=0.3,labelcolor] at (0.05,-0.25) {$\bfm{\Xi}^{(0)}$};
\node[scale=0.3,labelcolor] at (1.09,-0.45) {$\bfm{\Xi}^{(1)}$};
\node[scale=0.3,labelcolor] at (1.08,0.77) {$\bfm{\Xi}^{(2)}$};
\node[scale=0.3,labelcolor] at (-0.1,1.15) {$\bfm{\Xi}^{(3)}$};
\node[scale=0.3,labelcolor] at (-1.28,0.46) {$\bfm{\Xi}^{(4)}$};
\node[scale=0.3,labelcolor] at (-1.12,-0.65) {$\bfm{\Xi}^{(5)}$};
\node[scale=0.3,labelcolor] at (-0.05,-1.47) {$\bfm{\Xi}^{(6)}$};

\node[scale=0.3,edgecolor] at (0.45,-0.05) {$\Gamma^{(0)}$};
\node[scale=0.3,edgecolor] at (0.13,0.5) {$\Gamma^{(1)}$};
\node[scale=0.3,edgecolor] at (-0.5,-0.15) {$\Gamma^{(2)}$};
\node[scale=0.3,edgecolor] at (1.1,0.1) {$\Gamma^{(3)}$};
\node[scale=0.3,edgecolor] at (0.4,0.97) {$\Gamma^{(4)}$};
\node[scale=0.3,edgecolor] at (-0.6,0.88) {$\Gamma^{(5)}$};
\node[scale=0.3,edgecolor] at (-1.15,-0.1) {$\Gamma^{(6)}$};
\node[scale=0.3,edgecolor] at (-0.5,-1.15) {$\Gamma^{(7)}$};
\node[scale=0.3,edgecolor] at (0.48,-1) {$\Gamma^{(8)}$};

\end{scope}


\draw[line width=0.1mm] (-2.2, 0) -- (-1.7, 0);
\draw[line width=0.1mm] (-2.2, 0.5) -- (-1.7, 0.5);
\draw[line width=0.1mm] (-2.2, 0) -- (-2.2, 0.5);
\draw[line width=0.1mm] (-1.7, 0) -- (-1.7, 0.5);
\draw[-stealth,line width=0.1mm] (-2.24, -0.04) -- (-1.5, -0.04);
\draw[-stealth,line width=0.1mm] (-2.24, -0.04) -- (-2.24, 0.7);
\node[scale=0.3] at (-1.95,0.25) {$[0,1]^2$};
\node[scale=0.23] at (-1.5,0.05) {$\xi_1$};
\node[scale=0.23] at (-2.15,0.7) {$\xi_2$};


\draw[line width=0.1mm] (-1.7, -1.5) -- (-1.2, -1.5);
\draw[line width=0.1mm] (-1.7, -1) -- (-1.2, -1);
\draw[line width=0.1mm] (-1.7, -1.5) -- (-1.7, -1);
\draw[line width=0.1mm] (-1.2, -1.5) -- (-1.2, -1);
\draw[-stealth,line width=0.1mm] (-1.74, -1.54) -- (-1.0, -1.54);
\draw[-stealth,line width=0.1mm] (-1.74, -1.54) -- (-1.74, -0.8);
\node[scale=0.3] at (-1.45,-1.25) {$[0,1]^2$};
\node[scale=0.23] at (-1.0,-1.45) {$\xi_1$};
\node[scale=0.23] at (-1.65,-0.8) {$\xi_2$};


\draw[line width=0.1mm] (1.7, -0.2) -- (2.2, -0.2);
\draw[line width=0.1mm] (1.7, 0.3) -- (2.2, 0.3);
\draw[line width=0.1mm] (1.7, -0.2) -- (1.7, 0.3);
\draw[line width=0.1mm] (2.2, -0.2) -- (2.2, 0.3);
\draw[-stealth,line width=0.1mm] (1.66, -0.24) -- (2.4, -0.24);
\draw[-stealth,line width=0.1mm] (1.66, -0.24) -- (1.66, 0.5);
\node[scale=0.3] at (1.95,0.05) {$[0,1]^2$};
\node[scale=0.23] at (2.4,-0.15) {$\xi_1$};
\node[scale=0.23] at (1.75,0.5) {$\xi_2$};

\end{tikzpicture}
 }
\end{tabular}
\end{center}
\caption{ An example of a planar three patch domain~$\overline{\Omega}$ composed of the open patches~$\Omega^{(i)}$, $i\in \{0,1,2\}$, open edges~$\Gamma^{(i)}$, $i \in \{0, \ldots,8 \}$, and vertices~$\bfm{\Xi}^{(i)}$, $i \in \{0,\ldots,6 \}$, and with its geometry mappings~$\ab{G}^{(i)}$, $i \in \{0,1,2\}$.}
 \label{fig:three_patch_isogeometric_analysis}
\end{figure}


\subsection{Analysis-suitable $G^1$ multi-patch parameterization}

The planar multi-patch domain~$\overline{\Omega}$, introduced in Section~\ref{subsec:planar_multi-patch_domain}, is represented by a
particular multi-patch parameterization~$\ab{G}$ 
called analysis-suitable $G^1$ multi-patch parameterization~\cite{CoSaTa16}. A $C^0$-smooth multi-patch parameterization~$\ab{G}$ is called analysis-suitable $G^1$ if for each inner edge~$\Gamma^{(i)}$, $i \in \mathcal{I}^{I}_{\Gamma}$, with $\overline{\Gamma^{(i)}} = \overline{\Omega^{(i_0)}} \cap  \overline{\Omega^{(i_1)}}$, $i_0,i_1 \in \mathcal{I}_{\Omega}$, assuming without loss of generality that the two associated patch parameterizations~$\ab{G}^{(i_0)}$ and $\ab{G}^{(i_1)}$ are parameterized as
\begin{equation} \label{eq:standard_par}
\ab{G}^{(i_0)}(0,\xi) = \ab{G}^{(i_1)}(0,\xi) =:  \ab{G}^{(i)}_{0}(\xi), \quad \xi \in [0,1],
\end{equation}
there exist linear functions $\alpha^{(i,i_0)}:[0,1] \rightarrow \R$, $\alpha^{(i,i_1)}:[0,1] \rightarrow \R$, $\beta^{(i,i_0)}:[0,1] \rightarrow \R$ and $\beta^{(i,i_1)}: [0,1] \rightarrow \R$ such that
\begin{equation*}
\ab{G}^{(i,i_0)}_{1}(\xi) = \ab{G}^{(i,i_1)}_{1}(\xi), 
\quad \xi \in [0,1],   
\end{equation*}
with
\[
\ab{G}_{1}^{(i,\tau)}(\xi) = \frac{1}{\alpha^{(i,\tau)}(\xi)} \left(  \partial_1 \ab{G}^{(\tau)}(0,\xi) - \beta^{(i,\tau)}(\xi) \partial \ab{G}^{(i)}_{0}(\xi) 
\right), \;\; \tau\in \{i_0,i_1 \}.
\]
Parametrizing multi-patch domains by analysis-suitable $G^1$ multi-patch geometries allows to construct $C^1$-smooth isogeometric spline spaces with optimal approximation properties as proven in~\cite{HaTaTa26} and numerically demonstrated in e.g.~\cite{CoSaTa16, KaSaTa17b, KaSaTa19b, FaJuKaTa22, FaKaKoVi24, KaKoVi24b}, while the use of non-analysis-suitable $G^1$ multi-patch parameterizations can lead to $C^1$-smooth spline spaces with dramatically reduced approximation properties~\cite{KaSaTa17b}. Bilinear planar multi-patch geometries form an important subclass of the class of analysis-suitable $G^1$ planar multi-patch parameterizations~\cite{CoSaTa16}, but the class of analysis-suitable $G^1$ planar multi-patch parameterizations is much larger than the subclass of bilinear planar multi-patch geometries, cf.~\cite{KaSaTa17b, KaSaTa19b, FaJuKaTa22, FaKaKoVi24}. More precisely, any non-analysis-suitable planar $G^1$ multi-patch geometry can be approximated with high accuracy by an analysis-suitable $G^1$ planar multi-patch parameterization~\cite{KaSaTa17b, FaKaKoVi24}.


\subsection{Hierarchical B-splines and Truncated hierarchical B-splines} \label{subsec:HB_THB}

We recall 
the concept of hierarchical and truncated hierarchical B-spline basis constructions (e.g.~\cite{Kr97,GiJuSp2012}) {\new on} $[0,1]$ and $[0,1]^2$, which is used in this work. 
Let 
$$
\mathcal{D}^{(d)} = ( \mathcal{D}_{\ell}^{(d)} )_{\ell=0}^{N-1}, \quad 
\mathcal{D}^{(d)}_0 \supseteq \mathcal{D}^{(d)}_1 \supseteq \cdots \supseteq \mathcal{D}^{(d)}_{N-1}, \quad \mathcal{D}_0^{(d)} = [0,1]^d, 
$$
be a sequence of closed nested subdomains in $[0,1]^d$, $d \in \{1,2\}$, 
with $\mathcal{D}_{\ell}^{(d)}$, $\ell>0$, being a union of d-dimensional squares 
$$
 \prod_{i=1}^{d} \left[h_{\ell-1} \, j_i
, h_{\ell-1}(j_i+1) \right] \subseteq [0,1]^d, \quad {\rm for}\; {\rm some} \;\, j_i \in \{0,1,\ldots,\frac{1}{h_{\ell-1}}-1\},
$$
and let 
$$
\mathcal{S}_{h_0}^{{p},{r}}([0,1]) \subseteq 
\mathcal{S}_{h_1}^{{p},{r}}([0,1]) \subseteq \cdots \subseteq \mathcal{S}_{h_{N-1}}^{{p},{r}}([0,1]), 
$$
and
$$
\mathcal{S}_{h_0}^{\ab{p},\ab{r}}([0,1]^2) \subseteq 
\mathcal{S}_{h_1}^{\ab{p},\ab{r}}([0,1]^2) \subseteq \cdots \subseteq \mathcal{S}_{h_{N-1}}^{\ab{p},\ab{r}}([0,1]^2), 
$$
be sequences of nested spline
spaces defined on 
$[0,1]$ and $[0,1]^2$, respectively, for some initial $h_0$ and $h_\ell = h_0/2^\ell$.
Let further 
$$
\mathcal{N}^{p,r}_{\ell}(\mathcal{D}^{(1)}) = \{ N_{h_\ell,j}^{p,r} 
: \; j \in \mathcal{J}_{h_\ell} , \;   {\rm supp} \, N_{h_\ell,j}^{p,r} \subseteq {\mathcal{D}_{\ell}^{(1)}} \; {\rm and}\; {\rm supp} \, N_{h_\ell,j}^{p,r} \nsubseteq {\mathcal{D}_{\ell+1}^{(1)}} \}
$$
and 
$$
\mathcal{N}^{p,r}_{\ell}(\mathcal{D}^{(2)}) = \{ N_{h_\ell,\ab{j}}^{\ab{p},\ab{r}} 
: \; \ab{j} \in \bb{\mathcal{J}}_{h_\ell}, \;   {\rm supp} \, N_{h_\ell,\ab{j}}^{\ab{p},\ab{r}} \subseteq {\mathcal{D}_{\ell}^{(2)}} \; {\rm and}\; {\rm supp} \, N_{h_\ell,\ab{j}}^{\ab{p},\ab{r}} \nsubseteq {\mathcal{D}_{\ell+1}^{(2)}} \}
$$
denote the ``active'' B-splines from level $\ell$ with respect to nested subdomains $\mathcal{D}^{(d)}$, $d \in \{1,2\}$. 
Now, we can define the hierarchical spline space $\mathcal{H}^{p,r}(\mathcal{D}^{(d)})$  over $[0,1]^d$, $d\in \{1,2\}$, as
$$
\mathcal{H}^{p,r}(\mathcal{D}^{(d)}) = \Span 
\bigcup_{\ell=0}^{N-1} \mathcal{N}^{p,r}_{\ell}(\mathcal{D}^{(d)}).
$$
Note that all 
spanning functions 
of the hierarchical spline space $\mathcal{H}^{p,r}(\mathcal{D}^{(d)})$ are linearly independent, since all functions $\mathcal{N}^{p,r}_{\ell}(\mathcal{D}^{(d)})$ are linearly independent on $\mathcal{D}_{\ell}^{(d)}\backslash \mathcal{D}_{\ell+1}^{(d)}$, cf.~{\new \cite{GiJuSp2012}}. 

In order to define truncated hierarchical spline spaces, the concept of truncation at level $\ell+1$ has to be recalled. 
We first represent a spline function $\varphi \in \mathcal{N}^{p,r}_{\ell}(\mathcal{D}^{(d)})$ as a linear combination of basis functions in $\mathcal{N}^{p,r}_{\ell+1}(\mathcal{D}^{(d)})$, and then take into account only those basis functions having support partially also outside of $\mathcal{D}_{\ell+1}^{(d)}$. More precisely,
$$
{\rm tr_{\ell+1}} (\varphi) = {\rm tr_{\ell+1}} ( \hspace{-0.2cm}\sum_{\psi \in \mathcal{N}^{p,r}_{\ell+1}(\mathcal{D}^{(d)})} \hspace{-0.4cm} c_\psi \psi ) = \hspace{-0.4cm}
\sum_{\psi \in \mathcal{N}^{p,r}_{\ell+1}(\mathcal{D}^{(d)}) \atop {\rm supp}\;\psi \nsubseteq \mathcal{D}_{\ell+1}^{(d)}} \hspace{-0.6cm} c_\psi \psi.
$$
The truncated hierarchical spline space 
$\mathcal{T}^{p,r}(\mathcal{D}^{(d)})$  {\new on} $[0,1]^d$, $d\in \{1,2\}$, is then defined as
$$
\mathcal{T}^{p,r}(\mathcal{D}^{(d)}) = \Span \bigcup_{\ell=0}^{N-1} \left\{ tr_{N-1} (\cdots (tr_{\ell+1}(\varphi))\cdots) :\;\;  \varphi \in \mathcal{N}^{p,r}_{\ell}(\mathcal{D}^{(d)}) 
\right\}.
$$
Again, all 
spanning functions 
of the truncated hierarchical spline space $\mathcal{T}^{p,r}(\mathcal{D}^{(d)})$ are linearly independent, since all functions $\mathcal{N}^{p,r}_{\ell}(\mathcal{D}^{(d)})$ are linearly independent on $\mathcal{D}_{\ell}^{(d)}\backslash \mathcal{D}_{\ell+1}^{(d)}$, cf.~{\new \cite{GiJuSp2012}}. 

\subsection{The space of $C^1$-smooth isogeometric hierarchical spline functions} \label{subsec:Cs_space}

Let $\overline{\Omega}$ be a planar multi-patch domain parametrized by an analysis-suitable $G^1$ multi-patch geometry~$\ab{G}$ with individual patch parameterizations~$\ab{G}^{(i)} \in \mathcal{S}_{h_0}^{\ab{p},\ab{r}}([0,1]^2) \times \mathcal{S}_{h_0}^{\ab{p},\ab{r}}([0,1]^2)$, $i \in \mathcal{I}_{\Omega}$, where $h_0$ is the initial mesh size. In addition, let $\mathcal{D}^{i,(2)} = ( \mathcal{D}^{i,(2)}_\ell )_{\ell=0}^{N-1}$, $i \in \mathcal{I}_{\Omega}$, be sequences of nested subdomains in $[0,1]^2$ 
mapped 
to $\overline{\Omega^{(i)}}$ using $\ab{G}^{(i)}$. Then, we define by
$\V_h(\overline{\Omega})$ 
the space of isogeometric hierarchical spline functions {\new on} the multi-patch domain~$\overline{\Omega}$ 
given as 
\begin{equation*}
\V_h(\overline{\Omega}) =
\left\{
\phi \in L^2(\overline{\Omega}) :
\phi |_{\overline{\Omega^{(i)}}} \circ \ab{G}^{(i)} \in
\mathcal{H}^{p,r}(\mathcal{D}^{i,(2)}) ,
\; i \in \mathcal{I}_{\Omega}
\right\},
\end{equation*}
where $h$ is referred to the mesh size associated to the level~$\ell$ of refinement up to which it has been considered, i.e. $h=h_0/2^{N-1}$. 
Since
$\ab{G}^{(i)} \in \mathcal{S}_{h_0}^{\ab{p},\ab{r}}([0,1]^2) \times \mathcal{S}_{h_0}^{\ab{p},\ab{r}}([0,1]^2)$, $i \in \mathcal{I}_{\Omega}$, 
{we also have that}
$$
\ab{G}^{(i)} \in \mathcal{H}^{p,r}(\mathcal{D}^{i,(2)}) \times \mathcal{H}^{p,r}(\mathcal{D}^{i,(2)}). 
$$
We 
denote by $\V_h^1(\overline{\Omega})$ the corresponding space of $C^1$-smooth isogeometric hierarchical spline functions, i.e.,
\begin{equation*}
\V_h^1(\overline{\Omega}) = \V_h(\overline{\Omega}) \cap C^1(\overline{\Omega}).
\end{equation*}
An isogeometric spline function $\phi \in \V_h(\overline{\Omega})$ belongs to the space~$\mathcal{V}_h^1(\overline{\Omega})$ if and only if for any inner edge~$\Gamma^{(i)}$, $i \in \mathcal{I}^{I}_{\Gamma}$, with
$
\overline{\Gamma^{(i)}} = \overline{\Omega^{(i_0)}} \cap \overline{\Omega^{(i_1)}},\ i_0,i_1 \in \mathcal{I}_{\Omega},
$
assuming without loss of generality that the two associated patch parameterizations~$\ab{G}^{(i_0)}$ and $\ab{G}^{(i_1)}$ are parameterized as in~\eqref{eq:standard_par}, we have for all $\xi \in [0,1]$
\begin{equation}
\label{eq:gC}
\g_0^{(i,\LL)}[\phi](\xi) = \g_0^{(i,\RR)}[\phi](\xi) =: \gC^{(i)}_0[\phi](\xi), \quad
\g_1^{(i,\LL)}[\phi](\xi) = \g_1^{(i,\RR)}[\phi](\xi) =: \gC^{(i)}_1[\phi](\xi)
\end{equation}
 with
 \begin{align*}   \label{eq:gC2}
 \g_0^{(i,\Side)}[\phi](\xi) & = \left(\phi \circ \ab{G}^{(\tau)}\right)(0,\xi), \\
 \g_1^{(i,\Side)}[\phi](\xi) & = \frac{1}{\alpha^{(i,\Side)}(\xi)}\left( 
 \partial_1
 \left(\phi \circ \ab{G}^{(\tau)}\right)(0,\xi) -  
\beta^{(i,\Side)}(\xi) 
\partial \gC_0^{(i)}[\phi](\xi) \right),
 \end{align*}
for $\Side\in \{\LL,\RR\}$, see e.g. \cite{KaSaTa19a, KaVi20b}. 
Note that the 
function $\gC_{0}^{(i)}[\phi]$ 
describes the trace of the function $\phi$ along the inner edge~$\Gamma^{(i)}$, while the function $\gC_{1}^{(i)}[\phi]$ represents a specific 
transversal derivative of the function~$\phi$ across the inner edge~$\Gamma^{(i)}$, cf.~\cite{CoSaTa16}.
The $C^1$-smooth spline space~$\V_h^1(\overline{\Omega})$ can now be 
described as
\begin{equation}
\label{eq:spaceVs}
\V_h^1(\overline{\Omega})
=
\left\{
\phi \in \V_h(\overline{\Omega}) :
\g_0^{(i,\LL)}[\phi](\xi) = \g_0^{(i,\RR)}[\phi](\xi),\, 
\g_1^{(i,\LL)}[\phi](\xi) = \g_1^{(i,\RR)}[\phi](\xi),\, 
\xi \in [0,1], 
i \in \mathcal{I}_{\Gamma}^I
\right\}.
\end{equation}
Motivated by the one-level cases~\cite{KaSaTa19a, KaVi20b} to obtain a $C^1$-smooth spline space whose dimension is independent of the initial parameterizations~$\ab{G}^{(i)}$, $i \in \mathcal{I}_{\Omega}$, and of the valencies~$\nu_i$ of the vertices~$\bfm{\Xi}^{(i)}$, $i \in \mathcal{I}_{\Xi}$, we consider
instead of the space~$\V_h^1(\overline{\Omega})$ the subspace $\W_h^1(\overline{\Omega}) \subseteq \V_h^1(\overline{\Omega})$, which additionally requires the functions~$\phi \in \W_h^1(\overline{\Omega})$ being $C^2$-smooth at the vertices and having trace and derivative functions $\gC^{(i)}_0[\phi]$ and $\gC^{(i)}_1[\phi]$ from a specific space. More precisely, we choose the subspace 
\begin{equation}
\label{eq:spaceWs}
\W_h^1(\overline{\Omega})
=
\left\{
\begin{array}{ll}
\phi \in \V_h^1(\overline{\Omega}) : &
\gC^{(i)}_0[\phi] \in \mathcal{H}^{p,r+1}(\mathcal{D}^{i,(1)}), \;
\gC^{(i)}_1[\phi] \in \mathcal{H}^{p-1,r}(\mathcal{D}^{i,(1)}), \,
i \in \mathcal{I}_{\Gamma}^{I}, \mbox{ and} \\[0.5ex]
&
\phi \in C^{2}(\bfm{\Xi}^{(i)}), \;
i \in \mathcal{I}_{\Xi}^{I}
\end{array}
\right\},
\end{equation}
with 
$$
\mathcal{D}^{i,(1)} = \mathcal{D}^{i_0,(2)}\big|_{\ab{e}} \, \cap \, \mathcal{D}^{i_1,(2)}\big|_{\ab{e}},
$$
where
$
\overline{\Gamma^{(i)}} = \overline{\Omega^{(i_0)}} \cap \overline{\Omega^{(i_1)}},
$
$i \in \mathcal{I}^{I}_{\Gamma}$, and
$\ab{e} \subseteq [0,1]^2$ is the closed boundary edge of $[0,1]^2$ fulfilling
$
\ab{G}^{(i_0)}(\ab{e}) = \ab{G}^{(i_1)}(\ab{e}) = \overline{\Gamma^{(i)}}.
$
By further requiring $p \geq 3$, $1 \leq r \leq p-2$, and $h \leq \frac{p-r-1}{3-r+1}$, we ensure that 
for the spline space $\W_h^1(\overline{\Omega})$ the trace and derivative functions~$\gC_{0}^{(i)}[\phi]$ and $\gC_{1}^{(i)}[\phi]$ are h-refinable and that 
an independent $C^{2}$ interpolation at the single vertices can be done, cf.~\cite{KaSaTa19a, KaVi20b}.

\begin{rem}
In the definition of the spaces $\V_h^1(\overline{\Omega})$ and $\W_h^1(\overline{\Omega})$ in \eqref{eq:spaceVs} and \eqref{eq:spaceWs}, respectively, we can replace the hierarchical spline spaces $\mathcal{H}^{p,r}(\mathcal{D}^{i,(d)})$ with the truncated hierarchical spline spaces $\mathcal{T}^{p,r}(\mathcal{D}^{i,(d)})$.
\end{rem}

\begin{rem} \label{rem:comparison}
In~\cite{BrGiKaVa23}, a $C^1$-smooth hierarchical isogeometric spline space {\new on} an analysis-suitable $G^1$ multi-patch geometry has been generated as the span of a $C^1$-smooth hierarchical spline basis. The construction of this hierarchical spline basis is based on the one hand on the development and usage of a quite technical refinement mask for the $C^1$-smooth one-level basis functions and on the other hand on the use of a specific refinement of the mesh in the vicinity of a vertex. More precisely, whenever an element adjacent to a vertex has to be refined, then not just this one element needs to be refined but also some additional elements in its neighborhood. This specific refinement is needed since the $C^1$-smooth one-level basis functions are not necessarily locally linearly independent in the vicinity of a vertex. 

In contrast to~\cite{BrGiKaVa23}, our adaptive method for solving the 
{\new Kirchhoff--Love plate problem} will require neither the construction nor the use of $C^1$-smooth one-level basis functions or of a $C^1$-smooth hierarchical spline basis, but only the use of the standard B-splines and of their comparatively simple refinement mask. Similar to~\cite{BrGiKaVa23}, a specific refinement of the mesh in the vicinity of a vertex 
will be performed, by refining a slightly smaller neighborhood of elements compared to~\cite{BrGiKaVa23}, to ensure the imposition of linearly independent continuity conditions across the edges, cf. Section~\ref{subsec:adaptive_ref}. Note, however, that when using the same refinement as in~\cite{BrGiKaVa23}, the underlying $C^1$-smooth hierarchical discretization space of our proposed adaptive method, which, as mentioned, 
does not need to be explicitly computed, 
would be identical to the $C^1$-smooth hierarchical spline space~\cite{BrGiKaVa23} {\new in the special case of the Kirchhoff--Love plate problem given by the ``standard'' biharmonic equation with clamped boundary conditions}.
\end{rem}


\section{Dual-Primal Isogeometric Tearing and Interconnecting (IETI-DP)} \label{sec:IETI}

{\new
We describe the framework of a Dual-Primal Isogeometric Tearing and Interconnecting
(in short IETI-DP) method for solving the Kirchhoff--Love plate problem {\new on} planar multi-patch
domains $\Omega$, allowing adaptive refinement of the underlying mesh. The integration of
the adaptive refinement procedure into the IETI-DP framework and the implementation details
of the resulting adaptive method are presented in Section~4. Our presented IETI-DP method is a further development of the IETI-DP method~\cite{KaKoVi26} by simplifying the technicality of the IETI-DP algorithm, by proposing a preconditioner of the Schur complement formulation of the resulting dual-primal problem, and by providing the possibility to perform adaptive refinement.

\subsection{The problem statement}

We consider a homogeneous isotropic Kirchhoff--Love plate with Young's modulus $E>0$, material thickness $t>0$, and Poisson ratio of the plate material $\nu_{\mathrm p}\in(-1,1/2)$. {\new The mid-surface of the plate is assumed to be an open domain $\Omega\subset\mathbb{R}^2$ whose closure $\overline\Omega$ is a planar multi-patch domain possessing an analysis-suitable $G^1$ multi-patch parameterization
$\boldsymbol G$.} Let us recall the definition of the plate quantities that are needed to define the problem. 
The plate scalar bending stiffness 
and the plate bending stiffness matrix
are defined as
\[
 D_{\mathrm p} =\frac{Et^3}{12(1-\nu_{\mathrm p}^2)}>0, 
 \quad 
 {\rm and}
 \quad
  \boldsymbol D_{\mathrm p}
 =D_{\mathrm p}
 \begin{pmatrix}
  1 & \nu_{\mathrm p} & 0\\
  \nu_{\mathrm p} & 1 & 0\\
  0 & 0 & \dfrac{1-\nu_{\mathrm p}}{2}
 \end{pmatrix},
\]
respectively. 
Additionally, the bending-moment tensor is defined by
\[
 \boldsymbol{\mathcal M}(u)
 =D_{\mathrm p}\Big((1-\nu_{\mathrm p})\nabla^2u
                    +\nu_{\mathrm p}\Delta u\,\boldsymbol I\Big),
\]
where $\nabla^2$ denotes the Hessian
and $\boldsymbol I$ denotes the identity matrix. Furthermore, this gives the normal bending moment
as
\[
 \mathcal M_{nn}(u)
 =\boldsymbol n^T\boldsymbol{\mathcal M}(u)\boldsymbol n
 =D_{\mathrm p}\Big((1-\nu_{\mathrm p})\,\boldsymbol n^T\nabla^2u\,\boldsymbol n
                   +\nu_{\mathrm p}\,\Delta u\Big),
\]
where $\bfm{n}$ is the outward unit normal vector at the boundary~$\partial \Omega$.
For a free boundary segment, let $\boldsymbol t$ denote a unit tangent vector.  Then the twisting moment
and the effective Kirchhoff shear force are defined by
\[
 \mathcal M_{nt}(u)=\boldsymbol t^T\boldsymbol{\mathcal M}(u)\boldsymbol n,
 \qquad
 \mathcal Q_n(u)
 =-\boldsymbol n^T\operatorname{div}\boldsymbol{\mathcal M}(u)
  -\partial_{\ab{t}}\mathcal M_{nt}(u).
\]

We decompose the boundary $\partial \Omega$ into clamped, simply supported, and free parts,
\[
 \partial\Omega=\Gamma_C\cup\Gamma_S\cup\Gamma_F,
\]
where the three boundary parts are pairwise disjoint up to their common endpoints. The
Kirchhoff--Love plate problem aims to find a scalar unknown function $u$, representing the
transverse plate displacement, such that
\begin{equation}
 D_{\mathrm p}\,\Delta^2u(\boldsymbol x)=f(\boldsymbol x),\qquad \boldsymbol x\in\Omega,
  \label{eq:kirchhoff-strong}
\end{equation}
with the (nonhomogeneous) boundary conditions
\begin{equation}
 \begin{cases}
  u(\ab{x})=g_0(\ab{x}),\quad\partial_{\boldsymbol n}u(\ab{x})=g_1(\ab{x}),
       & \ab{x}\in\Gamma_C,\\[1mm]
  u(\ab{x})=g_0(\ab{x}),\quad\mathcal M_{nn}(u)(\ab{x})=g_2(\ab{x}),
       & \ab{x}\in\Gamma_S,\\[1mm]
  \mathcal M_{nn}(u)(\ab{x})=g_2(\ab{x}),\quad\mathcal Q_n(u)(\ab{x})=g_3(\ab{x}),
       & \ab{x}\in\Gamma_F.
 \end{cases}
 \label{eq:kirchhoff-boundary-conditions}
\end{equation}
Here, $f:\Omega\rightarrow\R$ and the boundary data $g_0$, $g_1$, $g_2$, and $g_3$ are assumed to
be sufficiently smooth on the boundary parts on which they are prescribed. Thus, on
$\Gamma_C$ both the displacement and its normal derivative are prescribed
essentially, on $\Gamma_S$ only the displacement is prescribed essentially while
the normal bending moment is natural, and on $\Gamma_F$ both the normal bending
moment and the effective Kirchhoff shear force are natural boundary conditions. If
$\partial\Omega$ is piecewise smooth, the natural corner condition associated with the
free boundary is imposed at boundary vertices at which the 
trace of the test function is not fixed by $\Gamma_C\cup\Gamma_S$; in the absence of concentrated
corner forces, the corresponding oriented jump of $\mathcal M_{nt}(u)$ vanishes.

Let 
$H_g^{2}(\Omega)$ and $H_0^{2}(\Omega)$ be the function spaces defined by 
\[
 H_g^2(\Omega)
 =\left\{\phi\in H^2(\Omega):
 \phi (\ab{x})=g_0 (\ab{x})\ \text{for } \ab{x} \in \Gamma_C\cup\Gamma_S,\;\,
 \partial_{\boldsymbol n}\phi (\ab{x}) =g_1 (\ab{x})\ \text{for } \ab{x} \in \Gamma_C\right\},
\]
and
\[
 H_0^2(\Omega)
 =\left\{\phi\in H^2(\Omega):
 \phi(\ab{x})=0\ \text{for } \ab{x} \in \Gamma_C\cup\Gamma_S,\;\,
 \partial_{\boldsymbol n}\phi (\ab{x}) =0\ \text{for } \ab{x} \in \Gamma_C\right\}.
\]
A possible weak form of problem~\eqref{eq:kirchhoff-strong} 
is given by finding $u \in H_g^{2}(\Omega)$ such that
\begin{equation}
\label{eq:weak_polyharmonic}
a(u,v) = F(v), \quad \mbox{for all } v \in H^{2}_0(\Omega),
\end{equation}
with the bilinear form
\begin{equation}  \label{eq:kirchhoff-bilinear}
 a(u,v) = \sum_{i \in \mathcal{I}_\Omega} a^{(i)}(u,v), \quad
a^{(i)}(u,v) =
 D_{\mathrm p}\int_{\Omega^{(i)}}
 \Big((1-\nu_{\mathrm p})\,\nabla^2u:\nabla^2v
      +\nu_{\mathrm p}\,\Delta u\,\Delta v\Big)\,\mathrm d\Omega^{(i)},
\end{equation}
where $\nabla^{2}u : \nabla^{2}v$ is the Frobenius product
\[
\nabla^{2}u : \nabla^{2}v
=
\partial_{11} u \, \partial_{11} v
+ 2 \partial_{12} u \, \partial_{12} v
+ \partial_{22} u \, \partial_{22} v,
\]
and the corresponding
linear functional 
\[
F(v) = \sum_{i \in \mathcal{I}_\Omega} F^{(i)}(v),
\]
with
\begin{equation}
 F^{(i)}(v)=
 \int_{\Omega^{(i)}} f(\ab{x})v(\ab{x})\,\mathrm d\Omega^{(i)}
 +\int_{(\Gamma_S\cup\Gamma_F)\cap\partial\Omega^{(i)}}
 g_2(\ab{x})\,\partial_{\boldsymbol n}v(\ab{x})\,\mathrm ds
 +\int_{\Gamma_F\cap\partial\Omega^{(i)}}
 g_3(\ab{x})\,v(\ab{x})\,\mathrm ds.
 \label{eq:kirchhoff-functional}
\end{equation}
Note that the local integrand in \eqref{eq:kirchhoff-bilinear} can equivalently be written as
$\boldsymbol\kappa(u)^T\boldsymbol D_{\mathrm p}\boldsymbol\kappa(v)$, with
$
\boldsymbol\kappa(u)
 =\begin{pmatrix}
   \partial_{11}u & \partial_{22}u & 2\partial_{12}u
 \end{pmatrix}^{T}.
$ 
 For homogeneous natural boundary data, $g_2=0$ on
$\Gamma_S\cup\Gamma_F$ and $g_3=0$ on $\Gamma_F$, the
boundary contribution to $F(v)$ vanishes.
Note also that setting $D_{\mathrm p}=1$ and
$\nu_{\mathrm p}=0$ gives the ``standard'' biharmonic equation. 
\begin{rem}
{\new In the fully free case, $\Gamma_F=\partial\Omega$ and
$\Gamma_C=\Gamma_S=\emptyset$, there are no essential boundary
conditions and $H_g^2(\Omega)=H_0^2(\Omega)=H^2(\Omega)$. Since the
Kirchhoff--Love plate bilinear form~\eqref{eq:kirchhoff-bilinear} then vanishes on the space of affine-linear functions $\mathbb P_1(\Omega)$, a fully free plate is determined only up
to an affine-linear function. Whenever a unique representative is required, we 
have to replace $
H^2(\Omega)$ by the normalized space
$
 \left\{\phi\in H^2(\Omega):
 \int_\Omega \phi(\boldsymbol x)\,p(\boldsymbol x)\,\mathrm d\Omega=0
 \;\text{for all }p\in\mathbb P_1(\Omega)\right\}.
$
In order to avoid this, we assume that $\Gamma_C \cup \Gamma_S$ is large 
enough that this additional normalization is not needed.}
\end{rem}
\begin{rem}
The alternative Laplacian bilinear form
\begin{equation} \label{eq:alternative_bilinear}
a_{\Delta}(u,v)
=
\sum_{i\in I_\Omega}
D_{\mathrm p}
\int_{\Omega^{(i)}}
\Delta u\,\Delta v
\,\mathrm d\Omega^{(i)},
\end{equation}
{\new as used in the IETI-DP method~\cite{KaKoVi26} for solving the ``standard'' biharmonic equation with clamped boundary conditions (i.e. $D_{\mathrm p}=1$, $\nu_{\mathrm p}=0$, $\Gamma_C=\partial \Omega$, $\Gamma_S=\Gamma_F=\emptyset$),} is equivalent to the Kirchhoff--Love plate bilinear form~\eqref{eq:kirchhoff-bilinear} in the purely clamped case
$\Gamma_C=\partial\Omega$. Indeed, then the test functions satisfy
$v=\partial_n v=0$ on $\partial\Omega$, and integration by parts gives
\[
\int_{\Omega}
\nabla^2 u:\nabla^2 v\,\mathrm d\Omega
=
\int_{\Omega}
\Delta u\,\Delta v\,\mathrm d\Omega.
\]
If $\Gamma_S$ or $\Gamma_F$ is nonempty, this equivalence does not hold in general: the
normal derivative of the test function is free on $\Gamma_S$, while both its trace and its
normal derivative are free on $\Gamma_F$. The Kirchhoff--Love plate bilinear bilinear form~\eqref{eq:kirchhoff-bilinear} therefore
has to be retained to recover the natural conditions involving $\mathcal M_{nn}(u)$ on
$\Gamma_S\cup\Gamma_F$ and $\mathcal Q_n(u)$ on $\Gamma_F$.
It is also convenient for the IETI-DP formulation,
since its kernel over a single unconstrained patch $\Omega^{(i)}$ consists only of
affine-linear functions, whereas the kernel of the Laplacian bilinear form~\eqref{eq:alternative_bilinear} over
$H^2(\Omega^{(i)})$ contains all harmonic functions.
\end{rem}
The weak problem~\eqref{eq:weak_polyharmonic} can equivalently be written as the minimization problem
\begin{equation}
 u=\underset{v\in H^2_g(\Omega)}{\operatorname{arg\,min}}
 \left(\frac12 a(v,v)-F(v)\right),
 \label{eq:KP-minimization}
\end{equation}
where the space $H_g^2(\Omega)$ incorporates the essential conditions on
$\Gamma_C\cup\Gamma_S$, 
while all natural
conditions on $\Gamma_S\cup\Gamma_F$ enter through the functional $F$. 

Let $W^1_{g,h}(\overline{\Omega})\subseteq W^1_h(\overline{\Omega})$ denote the common
discrete space with the essential boundary data imposed, i.e.,
\begin{equation}
 W^1_{g,h}(\overline{\Omega})=
 \left\{\phi\in W^1_h(\overline{\Omega}):
 \phi(\ab{x})=\widetilde g_0(\ab{x}) {\mbox{ for }}  \ab{x} \in \Gamma_C\cup\Gamma_S,\;
 \partial_{\boldsymbol n}\phi (\ab{x})=\widetilde g_1 (\ab{x}) { \mbox{ for }} \ab{x} \in \Gamma_C\right\}.
 \label{eq:KP-discrete-boundary-space}
\end{equation}
Here, $\widetilde{g}_{0}$ and $\widetilde{g}_{1}$ denote the $L^2$-orthogonal projections of the essential boundary data~$g_{0}$ and $g_{1}$ onto
$ \left\{ \phi|_{\bfm{x} \in \Gamma_C\cup\Gamma_S} : \right.$ $\left.\phi \in \W^{1}_h(\overline{\Omega}) \right\} $
and $ \left\{ \partial_{\ab{n}}\phi|_{\bfm{x} \in \Gamma_C} : \phi \in \W^{1}_h(\overline{\Omega}) \right\}$, respectively.
%
%
By applying a Galerkin projection with the finite-dimensional subspace~$\W_{g,h}^1(\overline{\Omega}) \subseteq H_{g}^{2}(\Omega)$, and by using the characterizations~\eqref{eq:spaceWs} and 
\eqref{eq:KP-discrete-boundary-space} for the spaces $\W_{h}^1(\overline{\Omega})$ and $\W_{g,h}^1(\overline{\Omega})$, respectively, the problem~\eqref{eq:KP-minimization} is transformed into the problem of finding $u_h \in \W^{1}_{g,h}(\overline{\Omega})$
as the solution of the constrained minimization problem
\begin{equation}
 u_h=\underset{v_h\in V_h(\Omega)}{\operatorname{arg\,min}}
 \left(\frac12 a(v_h,v_h)-F(v_h)\right)
 \label{eq:KP-discrete-minimization}
\end{equation}
subject to
\begin{align}
& \g_0^{(i,\LL)}[v_h](\xi) = \g_0^{(i,\RR)}[v_h](\xi), \;
\g_1^{(i,\LL)}[v_h](\xi) = \g_1^{(i,\RR)}[v_h](\xi), \;
\xi \in [0,1], \; i \in \mathcal{I}_{\Gamma}^I, \mbox{ and }
\label{eq:constraintsI1} \\
& \gC^{(i)}_0[v_h] \in \mathcal{H}^{p,r+1}(\mathcal{D}^{i,(1)}), \;
\gC^{(i)}_1[v_h] \in \mathcal{H}^{p-1,r}(\mathcal{D}^{i,(1)}), \;
i \in \mathcal{I}_{\Gamma}^{I}, \mbox{ and }
\label{eq:constraintsI2} \\
& v_h \in 
C^{2}(\bfm{\Xi}^{(i)}), \;
i \in \mathcal{I}_{\Xi}^{I}, \mbox{ and }
\label{eq:constraintsV} \\
& v_h(\bfm{x})=\widetilde g_0(\bfm{x}) {\mbox{ for }}
 \bfm{x}\in\Gamma_C\cup\Gamma_S,
\quad
\partial_{\boldsymbol n}v_h(\bfm{x})=\widetilde g_1(\bfm{x}),
{\mbox{ for }} \bfm{x}\in\Gamma_C.
\label{eq:constraintsB}
\end{align}
}

Due to the decomposition of
the isogeometric spline space~$\V_h(\overline{\Omega})$ into the direct sum 
\begin{equation*}
\label{eq:V}
\V_h(\overline{\Omega}) = \bigoplus_{i \in \mathcal{I}_\Omega} \V_h^{(i)}(\overline{\Omega}), \quad
\V_{h}^{(i)}(\overline{\Omega})
=
\Span \left\{
\phi_{{j}}^{(i)} :
j \in \JH
\right\}, \quad \JH=\{0,1,\ldots, n_{\mathcal{H}}^{i,(2)}-1 \},
\end{equation*}
with
the basis functions
\begin{equation*}
{\phi}^{(i)}_{{j}}(\bfm{x})  = 
\begin{cases}
   (\varphi_j \circ (\ab{G}^{(i)})^{-1})(\bfm{x}): \; \varphi_j \in \cup_{\ell=0}^{N-1}\mathcal{N}^{p,r}_{\ell} (\mathcal{D}^{i,(2)}); 
\mbox{ if }\f \; \bfm{x} \in \overline{\Omega^{(i)}},
\\ 0 \quad \mbox{ if }\f \, \bfm{x} \in \overline{\Omega} \backslash \overline{\Omega^{(i)}},
\end{cases} 
\end{equation*}
$n_{\mathcal{H}}^{i,(d)} = \dim \mathcal{H}^{p,r}(\mathcal{D}^{i,(d)})$, 
and $\mathcal{D}^{i,(2)}$ 
being the sequence of nested subdomains in $[0,1]^2$ mapped 
to $\overline{\Omega^{(i)}}$ using $\ab{G}^{(i)}$, 
each function $u_h \in \V_h(\overline{\Omega})$ can be represented as
\begin{equation} \label{eq:spline_representation}
u_h(\ab{x}) =\sum_{i \in \mathcal{I}_{\Omega}} 
\sum_{j \in \JH} u_{{j}}^{(i)} \phi_{{j}}^{(i)}(\ab{x}) = \sum_{i \in \mathcal{I}_{\Omega}} \left(\ab{u}^{(i)}\right)^T \ab{\phi}^{(i)}(\ab{x}) = 
\ab{u}^T \ab{\phi}(\ab{x}),
\end{equation}
with coefficients~$u_{{j}}^{(i)} \in \R$, column vectors of coefficients $\ab{u}^{(i)} = 
[u_{{j}}^{(i)}]_{j\in \JH}$ and $\ab{u}  = [\ab{u}^{(i)}]_{i \in \mathcal{I}_{\Omega}}$, and column vectors of basis functions $\ab{\phi}^{(i)} = 
 [\phi_{j}^{(i)}]_{j \in \JH} $ and $\ab{\phi} = [\ab{\phi}^{(i)}]_{i \in \mathcal{I}_{\Omega}}$.
For a function~$v_h \in \V_h(\overline{\Omega})$ with the spline representation~$v_h(\ab{x})=\ab{v}^T \ab{\phi}(\ab{x})$, the constraints~\eqref{eq:constraintsI1}--\eqref{eq:constraintsB} are linear with respect to the coefficients~$\ab{v}$, and can therefore also be expressed in matrix form
\[
\ab{C} \ab{v} = \ab{c}, \quad  
\ab{C} \in \R^{q \times  N_{\mathcal{H}}^{i,(2)}} , \quad
\ab{c} \in \R^{N_{\mathcal{H}}^{i,(2)}}, \quad N_{\mathcal{H}}^{i,(2)} = \sum_{i \in \mathcal{I}_{\Omega}} n_{\mathcal{H}}^{i,(2)},
\]
for some~$q \in \N$, assuming without loss of generality that 
$\Rank (\ab{C}) = q$.
The minimization problem~\eqref{eq:KP-discrete-minimization} 
with the constraints~\eqref{eq:constraintsI1}--\eqref{eq:constraintsB} can now be written as
\begin{equation*}
u_h =
\underset{\ab{C} \ab{v} = \ab{c}}{\underset{v_h \in \V_h(\overline{\Omega})}{\arg \min}}
\left( \frac{1}{2} a(v_h,v_h) - F(v_h) \right),
\end{equation*}
which is further equivalent to 
solving the
saddle point problem
\begin{equation}
\label{eq:large_problem}
\left(
\begin{array}{cc}
\ab{K} & \ab{C}^T \\
\ab{C} & \ab{0}
\end{array}
\right)
\left(
\begin{array}{c}
\ab{u} \\
\ab{\lambda}
\end{array}
\right)
=
\left(
\begin{array}{c}
\ab{f} \\
\ab{c}
\end{array}
\right)
\end{equation}
for the spline coefficients~$\ab{u}$ of a function $u_h \in \V_h(\overline{\Omega})$, i.e. $u_h(\ab{x}) = \ab{u}^T \ab{\phi}(\ab{x})$, 
and for Lagrange multipliers~$\ab{\lambda} \in \R^{q}$, 
where $\ab{K} \in \R^{N_{\mathcal{H}}^{i,(2)} \times N_{\mathcal{H}}^{i,(2)}}$ is the 
stiffness matrix of the form
\[
\ab{K} = \Diag (\{ \ab{K}^{(i)}\}_{i \in \mathcal{I}_{\Omega}}),
\quad
\ab{K}^{(i)} =
\left[
a^{(i)}(\phi^{(i)}_{{j}_1},\phi^{(i)}_{{j}_2})
\right]_{j_1,j_2 \in \JH},
\]
and $\ab{f}$ is the 
load vector of the form
\[
\ab{f} = [\ab{f}^{(i)}]_{i \in \mathcal{I}_\Omega}, \quad \ab{f}^{(i)} =
 \left[
F^{(i)}(\phi^{(i)}_{{j}})
\right]_{j \in \JH}. 
\]

\subsection{Formulation of a saddle point problem}


In case that the stiffness matrix~$\ab{K}$ is invertible, the saddle point problem~\eqref{eq:large_problem} is equivalent to the dual problem of finding first the Lagrange multipliers $\ab{\lambda}$ via the linear system
\begin{equation}
\label{eq:dual_problem}
\ab{C} \ab{K}^{-1} \ab{C}^T \ab{\lambda}
=
\ab{C} \ab{K}^{-1} \ab{f} - \ab{c},
\end{equation} 
and then by computing the coefficients~$\ab{u}$ via 
\begin{equation}
\label{eq:first_line}
\ab{u} = \ab{K}^{-1} \left( \ab{f} - \ab{C}^T \ab{\lambda} \right). 
\end{equation}
While Eq.~\eqref{eq:first_line} directly follows from the first row of~\eqref{eq:large_problem}, the linear system~\eqref{eq:dual_problem} is obtained by inserting~\eqref{eq:first_line} into the second row of~\eqref{eq:large_problem}. The advantage of solving the dual problem~\eqref{eq:dual_problem} and \eqref{eq:first_line} over the saddle point problem~\eqref{eq:large_problem} is that instead of solving the large problem~\eqref{eq:large_problem}, first the comparitively small problem~\eqref{eq:dual_problem} for the Lagrange multipliers~$\bfm{\lambda}$ is solved followed by problem~\eqref{eq:first_line} for the coefficients~$\ab{u}$ which can be further split into parallelizable, local problems for the coefficients~$\ab{u}^{(i)}$ of the single patches~$\Omega^{(i)}$, $i \in \mathcal{I}_{\Omega}$. Unfortunately, the dual problem~\eqref{eq:dual_problem} and \eqref{eq:first_line} cannot be employed in our case, since
the single matrices~$\ab{K}^{(i)}$, $i \in \mathcal{I}_{\Omega}$, and therefore also the matrix~$\ab{K}$, are 
not invertible. The idea is now to 
introduce an adaptation of the dual
problem~\eqref{eq:dual_problem} and \eqref{eq:first_line}, more precisely a dual-primal formulation (cf.~\cite{FarhatFETI-DP, KlPeSaJu12, SoglTakacs_IETI_Elasticity}), which allows us to solve the saddle point problem~\eqref{eq:large_problem} in a similar way as using the dual formulation~\eqref{eq:dual_problem} and \eqref{eq:first_line}. For this purpose, we first present an adaption of the saddle point problem~\eqref{eq:large_problem}.

We decompose the coefficient vector~$\ab{u}$ 
into 
\[
\ab{u} =
\left(
\begin{array}{c}
\ab{u}_B \\
\ab{u}_\Pi \\
\ab{u}_\Delta
\end{array}
\right), 
\]
where 
$\ab{u}_B = [\ab{u}_B^{(i)}]_{i \in \mathcal{I}_\Omega}$ 
consists of all coefficients of $\ab{u}$ that are involved in the {\new essential boundary constraints}~\eqref{eq:constraintsB}, 
$\ab{u}_{\Pi} = [\ab{u}_\Pi^{(i)}]_{i \in \mathcal{I}_\Omega}$ collects the coefficients of $\ab{u}$ involved in the $C^2$-smoothness conditions~\eqref{eq:constraintsV} at the inner vertices~$\ab{\Xi}^{(i)}$, $i \in \mathcal{I}_{\Xi}^I$,
and $\ab{u}_\Delta = [\ab{u}_\Delta^{(i)}]_{i \in \mathcal{I}_\Omega}$ contains the remaining coefficients of $\ab{u}$. 
Note that for each inner vertex~$\ab{\Xi}^{(i)}$, $i \in \mathcal{I}_{\Xi}^I$, we have exactly six coefficients of~$\ab{u}$ per patch, which belong to~$\ab{u}_{\Pi}$. 
Similarly, we can decompose 
the load vector~$\ab{f}$ into 
\[
\ab{f} =
\left(
\begin{array}{c}
\ab{f}_B \\
\ab{f}_\Pi \\
\ab{f}_\Delta
\end{array}
\right), \qquad
\ab{f}_B = [\ab{f}_B^{(i)}]_{i \in \mathcal{I}_\Omega},
\qquad
\ab{f}_\Pi = [\ab{f}_\Pi^{(i)}]_{i \in \mathcal{I}_\Omega},
\quad\ab{f}_\Delta = [\ab{f}_\Delta^{(i)}]_{i \in \mathcal{I}_\Omega},
\]
and the stiffness matrix~$\ab{K}$ into
\begin{equation*}
\ab{K} =
\left(
\begin{array}{ccc}
\ab{K}_{BB} & \ab{K}_{\Pi B}^{T} &  \ab{K}_{\Delta B}^{T} \\
\ab{K}_{\Pi B} & \ab{K}_{\Pi \Pi} & \ab{K}_{\Delta \Pi}^{T} \\
\ab{K}_{\Delta B} & \ab{K}_{\Delta \Pi} & \ab{K}_{\Delta \Delta} 
\end{array}
\right), 
\end{equation*}
with
$
\ab{K}_{BB} = \Diag (\{ \ab{K}_{BB}^{(i)}\}_{i \in \mathcal{I}_{\Omega}})$, 
$\ab{K}_{\Pi B} = \Diag (\{ \ab{K}_{\Pi B}^{(i)}\}_{i \in \mathcal{I}_{\Omega}})$,
$\ab{K}_{\Delta B} = \Diag (\{ \ab{K}_{\Delta B}^{(i)}\}_{i \in \mathcal{I}_{\Omega}})$,
$\ab{K}_{\Pi \Pi} = \Diag (\{ \ab{K}_{\Pi \Pi}^{(i)}\}_{i \in \mathcal{I}_{\Omega}})$, 
$\ab{K}_{\Delta \Pi} = \Diag (\{ \ab{K}_{\Delta \Pi}^{(i)}\}_{i \in \mathcal{I}_{\Omega}})$,
$\ab{K}_{\Delta \Delta} = \Diag (\{ \ab{K}_{\Delta \Delta}^{(i)}\}_{i \in \mathcal{I}_{\Omega}}).
$
In addition, we decompose the matrix~$\ab{C}$ into matrices $\ab{C}_B$, $\ab{C}_{\Xi}$, and $\ab{C}_{\Gamma}$. 
{\new Thereby, $\ab{C}_{B}$ is a partial permutation matrix that fixes the coefficients collected in
$\ab u_B$ to a prescribed boundary coefficient vector~$\ab g$ via
}
\begin{equation} \label{eq:constraints_CB}
\ab{C}_B \, \ab{u} = \ab{g},
\end{equation}
{\new
where $\ab g$ is not the vector of values of the boundary functions {$\widetilde{g}_0$} and {$\widetilde{g}_1$},
rather, it denotes the coefficient vector obtained from the prescribed essential boundary
data by the separate constrained boundary approximation described in
Section~\ref{subsec:coefficients_ub}. In that construction, the required $C^1$-compatibility
conditions at boundary vertices are also imposed. Thus,
\eqref{eq:constraints_CB} realizes the essential boundary conditions~\eqref{eq:constraintsB}
through these already computed boundary coefficients.
}
$\ab{C}_{\Xi} \, \ab{u} = \ab{0}$ represents the $C^{2}$-smoothness conditions~\eqref{eq:constraintsV} at all inner vertices $\bfm{\Xi}^{(i)}$, $i \in \mathcal{I}_{\Xi}^I$, and 
$\ab{C}_{\Gamma} \, \ab{u} = \ab{0}$ represents the $C^1$-smoothness conditions~\eqref{eq:constraintsI1} and~\eqref{eq:constraintsI2} at all inner edges~$\Gamma^{(i)}$, $i \in \mathcal{I}_\Gamma^{I}$, reduced by redundant information already given in the matrices~$\ab{C}_{B}$ and $\ab{C}_{\Xi}$. 
{\new The} matrices $\ab{C}_{B}$, $\ab{C}_{\Xi}$ and $\ab{C}_{\Gamma}$ can 
further be decomposed into
\[
\ab{C}_{B} =
\left(
\begin{array}{ccc}
\ab{C}_{B,B} & \ab{0} & \ab{0}
\end{array}
\right), \quad
\ab{C}_{\Xi} =
\left(
\begin{array}{ccc}
\ab{0} & \ab{C}_{\Xi,\Pi} & \ab{0}
\end{array}
\right) \quad
\mbox{and} \quad
\ab{C}_{\Gamma} =
\left(
\begin{array}{ccc}
\ab{C}_{\Gamma,B} & \ab{C}_{\Gamma,\Pi} & \ab{C}_{\Gamma,\Delta}
\end{array}
\right),
\]
respectively. 
The system~\eqref{eq:large_problem} can then be written as
\begin{equation}
\label{eq:decomposed_problem}
\left(
\begin{array}{cccccc}
\ab{K}_{BB} & \ab{K}_{\Pi B}^T &  \ab{K}_{\Delta B}^T & \ab{C}_{B,B}^T & \ab{0} &  \ab{C}_{\Gamma,B}^T \\
\ab{K}_{\Pi B} & \ab{K}_{\Pi \Pi} & \ab{K}_{\Delta \Pi}^T &\ab{0} & \ab{C}_{\Xi,\Pi}^T & \ab{C}_{\Gamma,\Pi}^T \\
\ab{K}_{\Delta B} & \ab{K}_{\Delta \Pi} & \ab{K}_{\Delta \Delta}&\ab{0} & \ab{0} & \ab{C}_{\Gamma,\Delta}^T \\
\ab{C}_{B,B} & \ab{0} & \ab{0} & \ab{0} & \ab{0} & \ab{0} \\
\ab{0} &\ab{C}_{\Xi,\Pi} & \ab{0} & \ab{0} & \ab{0} & \ab{0} \\
\ab{C}_{\Gamma,B} & \ab{C}_{\Gamma,\Pi} & \ab{C}_{\Gamma,\Delta} & \ab{0} & \ab{0} & \ab{0}
\end{array}
\right)
\left(
\begin{array}{c}
\ab{u}_B \\
\ab{u}_\Pi \\
\ab{u}_\Delta \\
\ab{\lambda}_B \\
\ab{\lambda}_{\Xi} \\
\ab{\lambda}_{\Gamma}
\end{array}
\right)
=
\left(
\begin{array}{c}
\ab{f}_B \\
\ab{f}_\Pi \\
\ab{f}_\Delta \\
\ab{g} \\
\ab{0} \\
\ab{0}
\end{array}
\right).
\end{equation}
{\new
With the above definition of $\ab g$ and with the coefficients ordered such that $\ab u_B$
forms the first block of $\ab u$, the matrix $\ab C_B$ merely selects these already
identified boundary coefficients. Therefore,
$\ab{C}_{B,B}=\ab I$ and consequently $\ab u_B=\ab g$.
}
Using some straightforward manipulations, the saddle point problem~\eqref{eq:decomposed_problem} simplifies to
\begin{equation} \label{eq:decomposed_problem_NoCB_PiDelta}
\left( \begin{array}{cccc}
{\ab{K}}_{\Pi \Pi} & {\ab{K}}_{\Delta \Pi }^T &  {\ab{C}}_{\Xi,\Pi}^T & {\ab{C}}_{\Gamma,\Pi}^T \\
{\ab{K}}_{\Delta \Pi } & \ab{K}_{\Delta \Delta} &  \ab{0} & {\ab{C}}_{\Gamma,\Delta}^T  \\
{\ab{C}}_{\Xi,\Pi} & \ab{0}  & \ab{0} & \ab{0}  \\
{\ab{C}}_{\Gamma,\Pi} & {\ab{C}}_{\Gamma,\Delta} & \ab{0} & \ab{0} 
\end{array} \right)
\left( 
\begin{array}{c}
\ab{u}_\Pi \\ 
\ab{u}_\Delta \\ 
\ab{\lambda}_{\Xi}\\
\ab{\lambda}_{\Gamma}
\end{array}
\right) = 
\left(
\begin{array}{c}
\overline{\ab{f}}_\Pi \\ 
\overline{\ab{f}}_\Delta\\
\ab{0}\\
\overline{\ab{g}}
\end{array}
\right),
\end{equation}
where $\overline{\ab{f}}_{\Pi} = \ab{f}_{\Pi} - \ab{K}_{\Pi B} \,{\ab{u}_{B}}$, $\overline{\ab{f}}_{\Delta} = \ab{f}_{\Delta} - \ab{K}_{\Delta B} \,{\ab{u}_{B}} $ and $\overline{\ab{g}} = - \ab{C}_{\Gamma,B} \ab{u}_B$. 
While the solving of the saddle point problem~\eqref{eq:decomposed_problem_NoCB_PiDelta} is explained in the next subsection and is done via a dual-primal problem, the detailed construction of the matrices $\ab{C}_{\Xi}$ and $\ab{C}_{\Gamma}$ is presented in  Section~\ref{subsec:constraint_CXi} and \ref{subsec:constraint_CGamma}, respectively. 

\subsection{Formulation and solving of a dual-primal problem}

{\new
In the following, let us assume that every patch $\Omega^{(i)}$ contains at least one inner vertex, or
has sufficient essential boundary conditions.
This means that we assume that the coefficients contained in \(\mathbf u_B^{(i)}\) and \(\mathbf u_\Pi^{(i)}\) eliminate the local affine-linear kernel.
}

Let us now further decompose the coefficients of $\ab{u}_\Delta$ into
\[
\ab{u}_\Delta =
\left(
\begin{array}{c}
\ab{u}_E \\
\ab{u}_R 
\end{array}
\right), 
\quad
\ab{u}_E = [\ab{u}_E^{(i)}]_{i \in \mathcal{I}_\Omega}, \quad
\ab{u}_R = [\ab{u}_R^{(i)}]_{i \in \mathcal{I}_\Omega},
\]
where the subindex $E$ refers to coefficients of $\ab{u}_\Delta$ associated to the edges and the first neighboring columns of coefficients, while $\ab{u}_R$ represent the remaining coefficients of $\ab{u}_\Delta$. 
Then we can write $\ab{C}_{\Gamma,\Delta} = \left(\begin{array}{cc}
\ab{C}_{\Gamma,E} &  \ab{0} 
\end{array}\right)$ with respect to the coefficients~$\ab{u}_{E}$ and $\ab{u}_{R}$.

In order to improve the conditioning of the system~\eqref{eq:decomposed_problem_NoCB_PiDelta}, let us first multiply the last row in~\eqref{eq:decomposed_problem_NoCB_PiDelta} by a particular orthogonal matrix.
More precisely, we replace the equations
\[
{\ab{C}}_{\Gamma,\Pi} \ab{u}_{\Pi} + {\ab{C}}_{\Gamma,\Delta} \ab{u}_{\Delta} = \overline{\ab{g}}
\]
by the equations
\[
\widetilde{{\ab{C}}}_{\Gamma,\Pi} \ab{u}_{\Pi} + \widetilde{{\ab{C}}}_{\Gamma,\Delta} \ab{u}_{\Delta} = \widetilde{\ab{g}},
\]
where
\[
\widetilde{{\ab{C}}}_{\Gamma,\Pi} = \ab{\Sigma}^{-1/2} \ab{Q} \, {\ab{C}}_{\Gamma,\Pi}, \qquad
\widetilde{{\ab{C}}}_{\Gamma,\Delta} = \ab{\Sigma}^{-1/2} \ab{Q} \, {\ab{C}}_{\Gamma,\Delta}, \qquad
\widetilde{\ab{g}} = \ab{\Sigma}^{-1/2} \ab{Q} \, \overline{\ab{g}},
\]
with $\ab{\Sigma}$ being a diagonal matrix with positive entries and $\ab{Q}$ being an orthogonal matrix such that
\[
{\ab{C}}_{\Gamma,E} \, {\ab{C}}_{\Gamma,E}^T = \ab{Q}^T \ab{\Sigma} \ab{Q}
\]
is a diagonalization of the symmetric positive definite matrix
$
{\ab{C}}_{\Gamma,E} \, {\ab{C}}_{\Gamma,E}^T.
$
As a consequence, we obtain
\begin{equation}
\label{eq:conditionCGammaTilda}
\widetilde{\ab{C}}_{\Gamma,E} \, \widetilde{\ab{C}}_{\Gamma,E}^T
=
\ab{\Sigma}^{-1/2} \ab{Q} \, {\ab{C}}_{\Gamma,E} \, {\ab{C}}_{\Gamma,E}^T \, \ab{Q}^T \ab{\Sigma}^{-1/2}
=
\ab{\Sigma}^{-1/2} \ab{Q} \, \ab{Q}^T \, \ab{\Sigma} \, \ab{Q} \, \ab{Q}^T \ab{\Sigma}^{-1/2}
=
\ab{I}.
\end{equation}

Furthermore, instead of imposing the continuity conditions ${\ab{C}}_{\Xi,\Pi}$ on the coefficients $\ab{u}_\Pi$, we introduce global primal degrees of freedom $\ab{w}_\Pi \in \R^{6 |\mathcal{I}_\Xi^I|}$ and replace the conditions
$
{\ab{C}}_{\Xi,\Pi} \, \ab{u}_\Pi = \ab{0}
$
by
\begin{equation}
\label{eq:localToglobal}
\ab{u}_\Pi = {\ab{N}}_{\Xi,\Pi} \, \ab{w}_\Pi, \quad {\ab{N}}_{\Xi,\Pi} = \ker \left( {\ab{C}}_{\Xi,\Pi} \right),
\end{equation}
where the construction of the matrix~${\ab{N}}_{\Xi,\Pi}$ is explained in Section~\ref{subsec:constraint_CXi}.
Decomposing ${\ab{N}}_{\Xi,\Pi}$ into
\[
{\ab{N}}_{\Xi,\Pi}
=
\left(
\begin{array}{cccc}
{{\ab{N}}^{(1)}_{\Xi,\Pi}}{}^T &
{{\ab{N}}^{(2)}_{\Xi,\Pi}}{}^T &
\cdots &
{{\ab{N}}^{(|\mathcal{I}_\Omega|)}_{\Xi,\Pi}}{}^T
\end{array}
\right)^T,
\qquad
{\ab{N}}^{(i)}_{\Xi,\Pi} \in \R^{6 \#_{\Xi}^{(i)} \times 6 |\mathcal{I}_\Xi^I|},
\]
where $\#_{\Xi}^{(i)}$ denotes the number of inner vertices of the patch $\Omega^{(i)}$, we can write
\begin{equation*}
\label{eq:global_local-Primal}
\ab{u}^{(i)}_\Pi = {\ab{N}}^{(i)}_{\Xi,\Pi} \, \ab{w}_\Pi, \quad i \in \mathcal{I}_\Omega,
\end{equation*}
which associates the global primal degrees $\ab{w}_\Pi$ with the patch-local degrees $\ab{u}^{(i)}_\Pi$, $i \in \mathcal{I}_\Omega$.
{\new
Multiplying the first block row in \eqref{eq:decomposed_problem_NoCB_PiDelta} (after replacing ${{\ab{C}}}_{\Gamma,\Pi}$ with $\widetilde{{\ab{C}}}_{\Gamma,\Pi}$) by matrix ${\ab{N}}_{\Xi,\Pi}^T$ from the left, and by
eliminating $\ab{u}_\Pi$ by means of~\eqref{eq:localToglobal}
}, 
we obtain
\begin{equation} \label{eq:GlobalPrimal_Short}
\left( \begin{array}{ccc}
{\ab{K}}_{\Delta \Delta}  &  {\ab{K}}_{\Delta \Pi } \,{\ab{N}}_{\Xi,\Pi} &  \widetilde{\ab{C}}_{\Gamma,\Delta}^T \\
{\ab{N}}_{\Xi,\Pi}^T {\ab{K}}_{\Delta \Pi }^T & {\ab{N}}_{\Xi,\Pi}^T \ab{K}_{\Pi \Pi} {\ab{N}}_{\Xi,\Pi} & {\ab{N}}_{\Xi,\Pi}^T \widetilde{\ab{C}}_{\Gamma,\Pi}^T  \\[0.1cm]
\widetilde{\ab{C}}_{\Gamma,\Delta} &  \widetilde{\ab{C}}_{\Gamma,\Pi} {\ab{N}}_{\Xi,\Pi} & \ab{0} 
\end{array} \right)
\left( 
\begin{array}{c}
\ab{u}_\Delta \\ 
\ab{w}_\Pi \\ 
\ab{\lambda}_{\Gamma}
\end{array}
\right) = 
\left(
\begin{array}{c} 
\overline{\ab{f}}_\Delta\\
{\ab{N}}_{\Xi,\Pi}^T \overline{\ab{f}}_\Pi\\
\widetilde{\ab{g}}
\end{array}
\right).
\end{equation}
Using a basis transformation which guarantees the orthogonality of
primal and dual basis functions, i.e.~replacing $\ab{u}_\Delta$ with 
\begin{equation} \label{eq:orthogonal_basis}
\ab{w}_\Delta = \ab{u}_\Delta + {\ab{K}}_{\Delta \Delta}^{-1} {\ab{K}}_{\Pi \Delta}^T {\ab{N}}_{\Xi,\Pi} \ab{w}_\Pi,
\end{equation}
system \eqref{eq:GlobalPrimal_Short} transforms into 
\begin{equation} \label{eq:GlobalPrimal_Short_Orthogonal}
\left( \begin{array}{ccc}
{\ab{K}}_{\Delta \Delta}  &  \ab{0} &  \widetilde{\ab{C}}_{\Gamma,\Delta}^T \\
\ab{0} & \ab{S}_{\Pi \Pi} & {\ab{B}_\Pi}^T \\
\widetilde{\ab{C}}_{\Gamma,\Delta} &  {\ab{B}_\Pi} & \ab{0} 
\end{array} \right)
\left( 
\begin{array}{c}
\ab{w}_\Delta \\ 
\ab{w}_\Pi \\ 
\ab{\lambda}_{\Gamma}
\end{array}
\right) = 
\left(
\begin{array}{c} 
\overline{\ab{f}}_\Delta\\
\widetilde{\ab{f}}_\Pi\\
\widetilde{\ab{g}}
\end{array}
\right),
\end{equation}
where
$$
\ab{S}_{\Pi \Pi} = 
{{\ab{N}}^T_{\Xi,\Pi}} 
 \left( 
\ab{K}_{\Pi \Pi} - \ab{K}_{\Pi \Delta} \ab{K}^{-1}_{\Delta \Delta} \ab{K}_{\Delta \Pi}
\right)
{\ab{N}}_{\Xi,\Pi} = 
\sum_{i \in \mathcal{I}_\Omega} 
{{\ab{N}}^{(i)}_{\Xi,\Pi}}{}^T 
 ( 
\ab{K}^{(i)}_{\Pi \Pi} - \ab{K}^{(i)}_{\Pi \Delta} {\ab{K}^{(i)}_{\Delta \Delta}}{\hspace{-0.1cm}}^{-1} \ab{K}^{(i)}_{\Delta \Pi}
)
{\ab{N}}^{(i)}_{\Xi,\Pi},
$$
$$
\ab{B}_\Pi = \left( 
\widetilde{\ab{C}}_{\Gamma, \Pi} - \widetilde{\ab{C}}_{\Gamma, \Delta} \ab{K}^{-1}_{\Delta \Delta} \ab{K}_{\Delta \Pi}
\right) {\ab{N}}_{\Xi,\Pi} = 
\sum_{i \in \mathcal{I}_\Omega} \left( 
\widetilde{\ab{C}}^{(i)}_{\Gamma, \Pi} - \widetilde{\ab{C}}^{(i)}_{\Gamma, \Delta} {\ab{K}^{(i)}_{\Delta \Delta}}{\hspace{-0.1cm}}^{-1} \ab{K}^{(i)}_{\Delta \Pi}
\right) {\ab{N}}^{(i)}_{\Xi,\Pi},
$$
and
$$
\widetilde{\ab{f}}_{\Pi} = {\ab{N}}_{\Xi,\Pi}^T \left( 
\overline{\ab{f}}_{\Pi}-
\ab{K}_{\Pi \Delta} \ab{K}^{-1}_{\Delta \Delta} \overline{\ab{f}}_{\Delta}
\right) = \sum_{i \in \mathcal{I}_\Omega}
{{\ab{N}}^{(i)}_{\Xi,\Pi}}{}^T \left( 
\overline{\ab{f}}^{(i)}_{\Pi}-
\ab{K}^{(i)}_{\Pi \Delta} {\ab{K}^{(i)}_{\Delta \Delta}}{\hspace{-0.1cm}}^{-1} \overline{\ab{f}}^{(i)}_{\Delta}
\right).
$$
From the first two rows of \eqref{eq:GlobalPrimal_Short_Orthogonal}, we now obtain 
\begin{equation} \label{eq:wPi_wDelta}
\ab{w}^{(i)}_{\Delta} = {{\ab{K}}^{(i)}_{\Delta \Delta}}{\hspace{-0.1cm}}^{-1}
\left( \overline{\ab{f}}^{(i)}_\Delta - \widetilde{\ab{C}}^{(i)}_{\Gamma,\Delta}{\hspace{-0.1cm}}^T \ab{\lambda}_{\Gamma} \right), \;\,
i \in \mathcal{I}_\Omega,\quad {\rm and} \quad
\ab{w}_{\Pi} = {\ab{S}}_{\Pi \Pi}^{-1}
\left( \widetilde{\ab{f}}_\Pi - \ab{B}_\Pi^T \ab{\lambda}_{\Gamma} \right). 
\end{equation}
Inserting \eqref{eq:wPi_wDelta} into the third row of \eqref{eq:GlobalPrimal_Short_Orthogonal} implies the symmetric positive definite linear system 
\begin{equation} \label{eq:spdFr}
 \ab{F} \ab{\lambda}_{\Gamma} = \ab{r},
\end{equation}
where
\begin{equation} \label{eq:spdMatrixF}
\ab{F} = \widetilde{\ab{C}}_{\Gamma,\Delta} \ab{K}^{-1}_{\Delta \Delta} \widetilde{\ab{C}}_{\Gamma,\Delta}^T + \ab{B}_\Pi \ab{S}^{-1}_{\Pi \Pi} \ab{B}_\Pi^T = \sum_{i \in \mathcal{I}_\Omega}  
\widetilde{\ab{C}}^{(i)}_{\Gamma,\Delta} {\ab{K}^{(i)}_{\Delta \Delta}}{\hspace{-0.1cm}}^{-1} \widetilde{\ab{C}}^{(i)}_{\Gamma,\Delta}{\hspace{-0.1cm}}^T + \ab{B}_\Pi \ab{S}^{-1}_{\Pi \Pi} \ab{B}_\Pi^T,
\end{equation}
and
$$
\ab{r} = \widetilde{\ab{C}}_{\Gamma,\Delta} \ab{K}^{-1}_{\Delta \Delta} \overline{\ab{f}}_{\Delta} + \ab{B}_\Pi \ab{S}^{-1}_{\Pi \Pi} \widetilde{\ab{f}}_{\Pi} - \widetilde{\ab{g}} = 
\sum_{i \in \mathcal{I}_\Omega}  
\widetilde{\ab{C}}^{(i)}_{\Gamma,\Delta} {\ab{K}^{(i)}_{\Delta \Delta}}{\hspace{-0.1cm}}^{-1} \overline{\ab{f}}^{(i)}_{\Delta} + \ab{B}_\Pi \ab{S}^{-1}_{\Pi \Pi} \widetilde{\ab{f}}_{\Pi} - \widetilde{\ab{g}}.
$$
To summarize, we first solve the system~\eqref{eq:spdFr}. Then, we compute $\ab{w}_\Pi$ and $\ab{w}_\Delta$ using~\eqref{eq:wPi_wDelta} and finally obtain the solution vectors $\ab{u}_\Pi$ and $\ab{u}_\Delta$ by means of~\eqref{eq:localToglobal} and~\eqref{eq:orthogonal_basis}.

When solving the system~\eqref{eq:spdFr}, it is desirable to use a good preconditioner. Recall that
\[
\ab{K}_{\Delta \Delta} =
\left(
\begin{array}{cc}
\ab{K}_{EE} & \ab{K}_{RE}^T \\
\ab{K}_{RE} & \ab{K}_{RR}
\end{array}
\right),
\qquad
\widetilde{\ab{C}}_{\Gamma,\Delta} =
\left(
\begin{array}{cc}
\widetilde{\ab{C}}_{\Gamma,E} & \ab{0}
\end{array}
\right).
\]
Then, the matrix $\ab{F}$ in~\eqref{eq:spdMatrixF} simplifies to
\[
\ab{F}
=
\sum_{i \in \mathcal{I}_\Omega}
\widetilde{\ab{C}}^{(i)}_{\Gamma,E}
{\ab{S}^{(i)}_{EE}}^{-1}
\widetilde{\ab{C}}^{(i)}_{\Gamma,E}{}^T
+
\ab{B}_\Pi \ab{S}^{-1}_{\Pi \Pi} \ab{B}_\Pi^T,
\]
where
$
{\ab{S}^{(i)}_{EE}}
=
{\ab{K}^{(i)}_{EE}}
-
{\ab{K}^{(i)}_{ER}}
{\ab{K}^{(i)}_{RR}}^{-1}
{\ab{K}^{(i)}_{RE}}.
$
The property~\eqref{eq:conditionCGammaTilda} indicates that the Dirichlet preconditioner
\begin{equation}\label{eq:Dirichlet_prec}
\ab{M}
=
\sum_{i \in \mathcal{I}_\Omega}
\widetilde{\ab{C}}^{(i)}_{\Gamma,E}
{\ab{S}^{(i)}_{EE}}
\widetilde{\ab{C}}^{(i)}_{\Gamma,E}{}^T
\end{equation}
may already be a good choice, cf.~\cite{Ta25}.
Therefore, instead of solving the system~\eqref{eq:spdFr}, we solve the preconditioned system
\begin{equation}  \label{eq:preconditionedSystem}
\ab{M} \ab{F} \ab{\lambda}_{\Gamma} = \ab{M}\ab{r}.
\end{equation}
{\new In the following example, the behavior of the proposed Dirichlet preconditioner under {\em uniform refinement} is numerically investigated and compared with standard diagonal scaling. Furthermore, also the good behavior of the proposed Dirichlet preconditioner under {\em adaptive refinement} is numerically demonstrated on the basis of several examples in Section~\ref{sec:numerical_examples}.}
\begin{ex} \label{ex:Example5}
{\new  
We consider 
the four-patch domain with two inner extraordinary vertices shown in Fig.~\ref{fig:example5_b}, 
and use the harmonic load $f(\ab{x})= f(x,y)=\cos(\frac{\pi}{2} x) \cos(\frac{\pi}{2} y)$. We set 
$$
D_{\mathrm p}=1, \quad \nu_{\mathrm p}=0, \quad \Gamma_C=\partial \Omega,\quad \Gamma_S=\Gamma_F=\emptyset,
$$
and all boundary data active on $\Gamma_C$ 
to be homogeneous, i.e.,
$g_0=g_1=0$. We compare the condition numbers $\kappa$ of the linear system~\eqref{eq:spdFr} obtained 
using standard diagonal scaling and 
the proposed Dirichlet preconditioner~\eqref{eq:Dirichlet_prec}, cf.~Fig.~\ref{fig:example5} and Table~\ref{tab:ConditionAndStepsUniform}. Table~\ref{tab:ConditionAndStepsUniform} also reports for the two preconditioners the required number of iterations using the conjugate gradient method for solving the linear system~\eqref{eq:spdFr}. The iteration is terminated once the residual has been reduced by a factor of $10^{-8}$ relative to the initial residual. 
In Table~\ref{tab:ConditionAndStepsUniform}, we additionally compare two successive solutions $u_h$ and $u_{h/2}$ by computing the corresponding $h$--$h/2$ type relative errors 
(see e.g.~\cite{ErGaPr20}) with respect to the $L^2$ norm and a norm equivalent to the $H^2$-seminorm. The concept of calculating these $h$--$h/2$ type relative errors
is explained in detail in Section~\ref{sec:numerical_examples}.
}

\begin{figure}[htb!]
    \centering
    \begin{tabular}{cc}
  \includegraphics[scale=0.23]{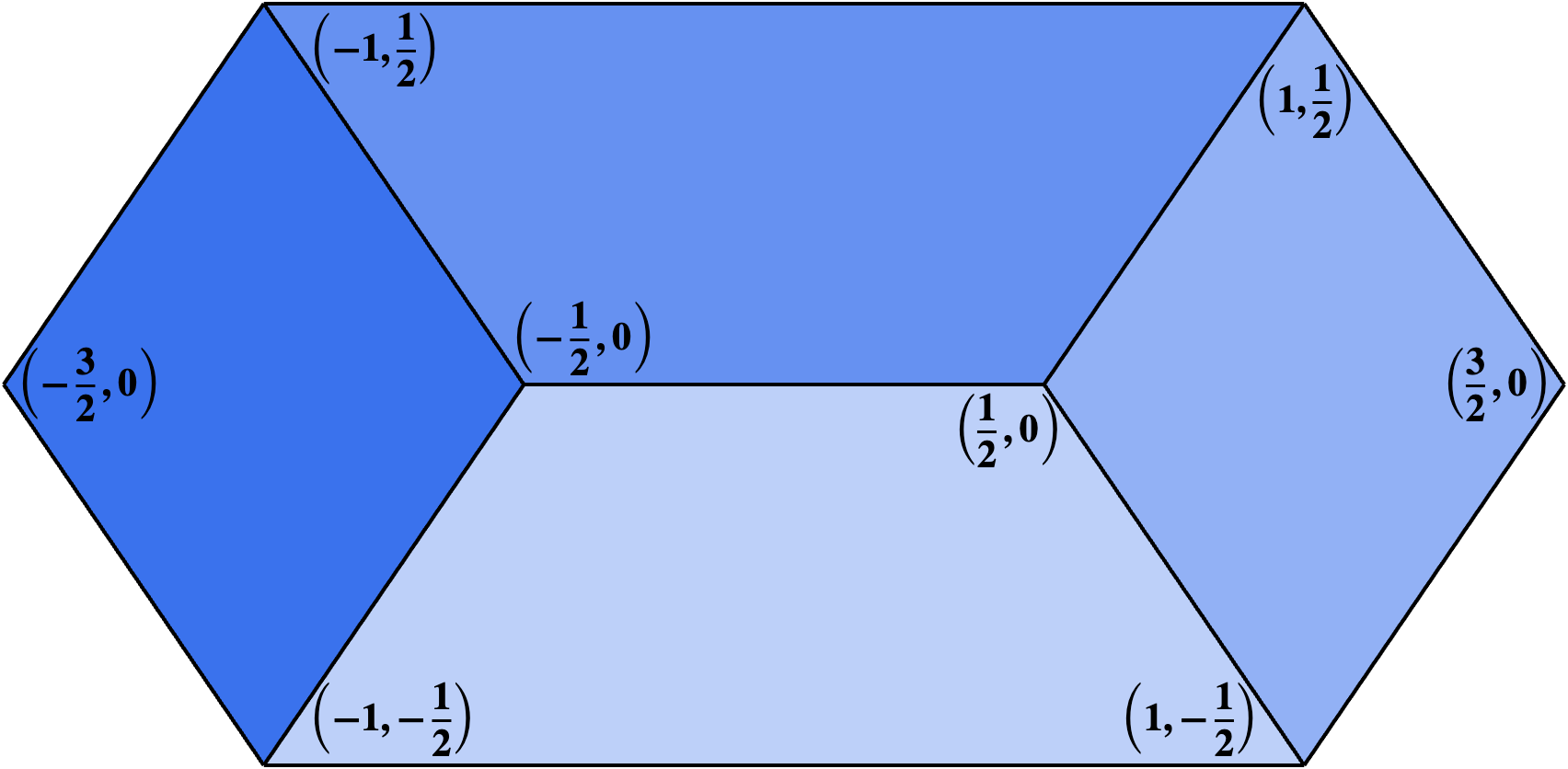}   & \hskip1.em
  \includegraphics[scale=0.33]{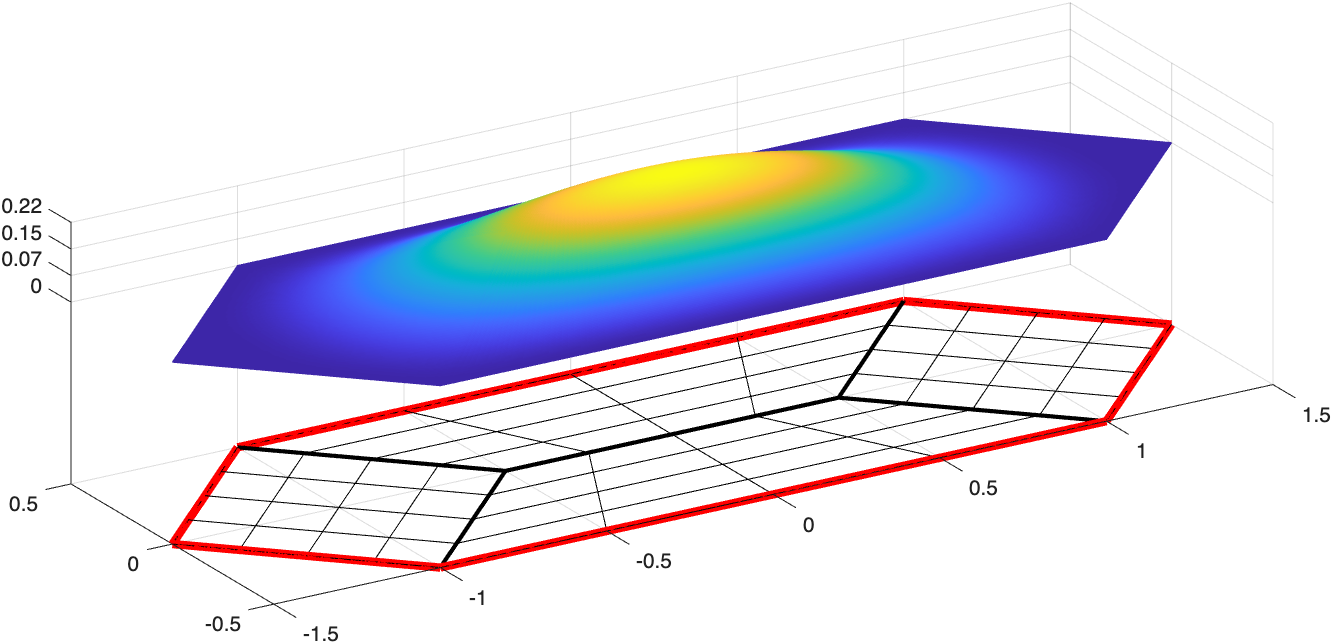} 
  \end{tabular}
\caption{\new Example~\ref{ex:Example5}.
The bilinear four-patch domain (left) and the associated numerical solution, scaled with factor 100 (right).
}
    \label{fig:example5_b}
\end{figure}

\begin{figure}[htb!]
    \centering
    \begin{tabular}{cc}
  \includegraphics[scale=0.3]{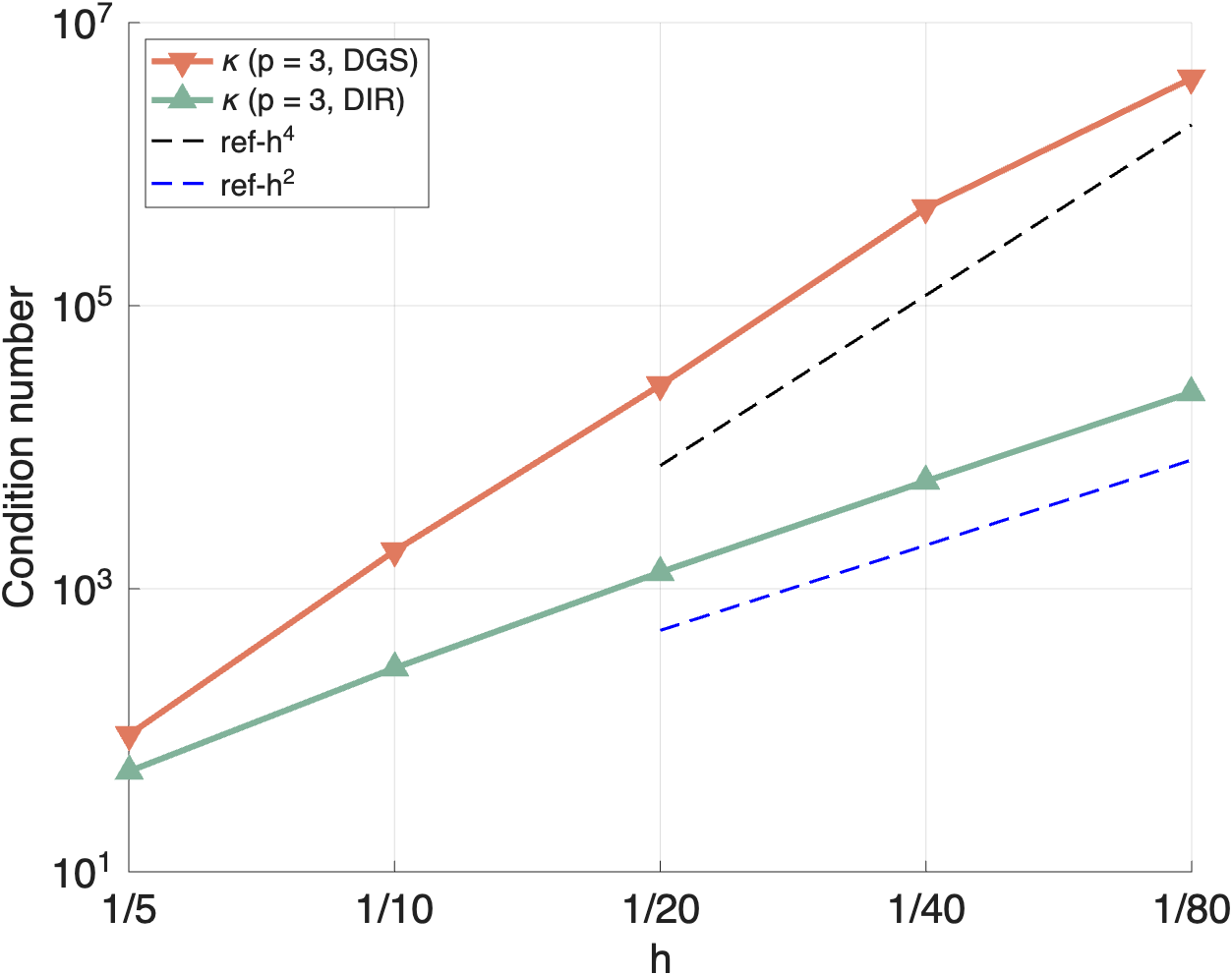}   & \hskip4em
  \includegraphics[scale=0.3]{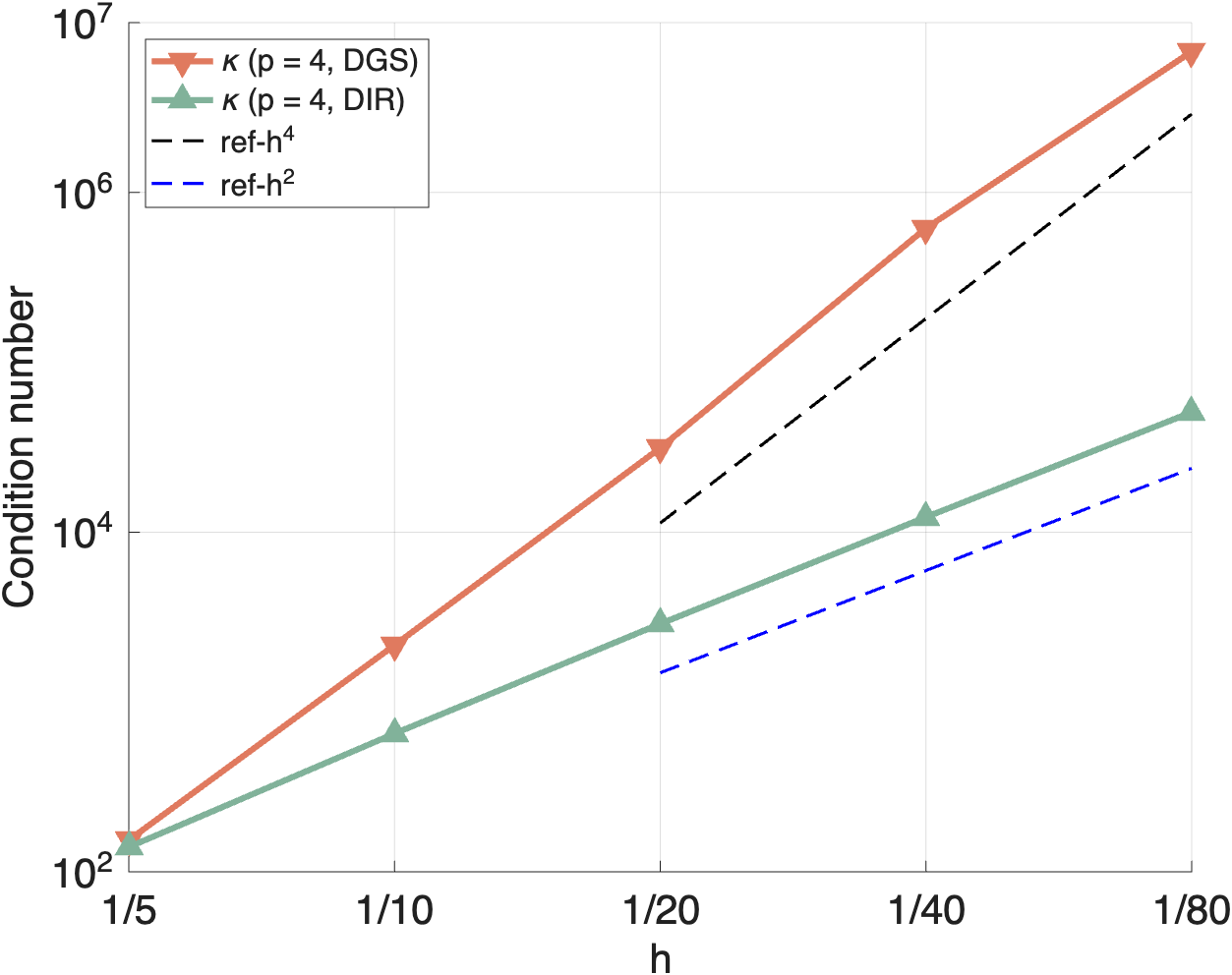} \\[0.2cm]
   Condition numbers $\kappa$ for $p=3$  &  Condition numbers $\kappa$ for $p=4$ \\[0.1cm]
  \end{tabular}
\caption{\new Example~\ref{ex:Example5}. 
Condition numbers $\kappa$ with respect to the number of Lagrange multipliers for $p=3$ (left) and $p=4$ (right), visualized in 
red lines 
for diagonal scaling (DGS) and in green lines 
for the Dirichlet preconditioner (DIR).}
    \label{fig:example5}
\end{figure}

\begin{table}[htb!]
\centering
\begin{tabular}{|c||c|c|c|c||c|c|c|c||}
\hline
& \multicolumn{4}{c||}{Dirichlet preconditioner $\ab{M}$} &
\multicolumn{4}{c||}{Diagonal scaling} \\
\hline
\hline
Level~$\ell$ & $\kappa$ & it & $L^2$ & $H^2$ & $\kappa$ & it & $L^2$ & $H^2$ \\
 \hline         
0 & 51.20 & 13 & - & - & 92.72 & 36 & - & -\\ 
1 & 2.75e+02 & 24 & 4.35e-03 & 3.61e-02  & 1.85e+03 & 74  & 4.35e-03 & 3.61e-02 \\
2 & 1.31e+03 & 45 & 3.49e-04 & 9.33e-03  & 2.74e+04 & 143 & 3.49e-04 & 9.33e-03 \\ 
3 & 5.76e+03 & 90 & 2.42e-05 & 2.27e-03  & 4.89e+05 & 334 & 2.42e-05 & 2.27e-03 \\ 
4 & 2.44e+04 & 189 & 1.39e-06 & 4.98e-04 & 4.04e+06 & 786 & 1.39e-06 & 4.98e-04 \\ 
\hline
\end{tabular}
\vskip2em
\centering
\begin{tabular}{|c||c|c|c|c||c|c|c|c||}
\hline
& \multicolumn{4}{c||}{Dirichlet preconditioner $\ab{M}$} &
\multicolumn{4}{c||}{Diagonal scaling} \\
\hline
\hline
Level~$\ell$ & $\kappa$ & it & $L^2$ & $H^2$ & $\kappa$ & it & $L^2$ & $H^2$ \\
 \hline         
0 & 1.40e+02 & 18  &    -     &    -     & 1.54e+02 & 42   &    -     &    - \\ 
1 & 6.51e+02 & 30  & 4.20e-04 & 9.36e-03 & 2.16e+03 & 88   & 4.20e-04 & 9.36e-03 \\
2 & 2.88e+03 & 62  & 2.74e-05 & 2.23e-03 & 3.13e+04 & 184  & 2.74e-05 & 2.23e-03 \\ 
3 & 1.22e+04 & 126 & 1.39e-06 & 4.38e-04 & 6.16e+05 & 415  & 1.39e-06 & 4.38e-04 \\ 
4 & 5.08e+04 & 255 & 4.89e-08 & 6.42e-05 & 6.75e+06 & 1001 & 5.41e-08 & 6.42e-05 \\ 
\hline
\end{tabular}
\caption{\new Example~\ref{ex:Example5}. 
Condition numbers~$\kappa$, iteration counts, and the reported error measures with respect to the $L^2$ and to a norm equivalent to the $H^2$-seminorm for five levels of uniform refinement for $p=3$ (top) and $p=4$ (bottom).}
\label{tab:ConditionAndStepsUniform}
\end{table}

\end{ex}

\section{Implementation details of the adaptive IETI-DP method} \label{sec:implementation_details}

The presented IETI-DP method from Section~\ref{sec:IETI}, which is based on the use of (truncated) hierarchical B-splines for representing the solution~$u_h$ of the 
{\new Kirchhoff--Love plate problem~\eqref{eq:kirchhoff-strong}}, allows to perform adaptive refinement of the underlying mesh to compute the solution~$u_h$. In this section,  we present the details of the resulting adaptive IETI-DP framework, which are needed to implement the method. This includes the introduction of the selected adaptive refinement strategy, the computation of the boundary coefficient vector~$\ab{u}_{B}$ as well as the assembly of the constraint matrices $\ab{C}_{\Xi,\Pi}$ and $\ab{C}_{\Gamma}$ and of the transformation matrix~$\ab{N}_{\Xi,\Pi}$. For the calculation of the coefficient vector~$\ab{u}_{B}$ and of the constraint matrices $\ab{C}_{\Xi,\Pi}$ and $\ab{C}_{\Gamma}$, we follow the approach described in~\cite{KaKoVi26}, and adapt it accordingly to the IETI-DP technique proposed in Section~\ref{sec:IETI}, particularly with regard to the case of adaptive refinement.

\subsection{The adaptive refinement procedure} \label{subsec:adaptive_ref}

The selected refinement strategy relies on the commonly used adaptive loop (cf.~\cite{BuGaGiPrVa22})
\begin{equation} \label{eq:modules}
\mbox{Solve} \rightarrow \mbox{Estimate} \rightarrow \mbox{Mark} \rightarrow \mbox{Refine}
\end{equation}
to pass from one 
level of refinement to the next one. The aim is to construct from a current hierarchical mesh a new, refined hierarchical mesh by refining exactly those regions of the mesh where the current approximation should be improved while preserving the structural properties of the mesh required by the (truncated) hierarchical construction. Below, we explain the integration of the four modules~\eqref{eq:modules} into our adaptive IETI-DP method to achieve this goal.

In the module {\em Solve}, we compute for a given hierarchical mesh a solution~$u_h$ of the 
{\new Kirchhoff--Love plate problem~\eqref{eq:kirchhoff-strong}} using the IETI-DP method presented in Section~\ref{sec:IETI}. The resulting solution~$u_h$ is then employed in the module {\em Estimate} to calculate for each mesh element a posteriori error estimate to decide which regions have to be refined. Instead of using a standard residual error estimator, which involves 
the computation of third and fourth order derivatives for the 
{\new Kirchhoff--Love plate problem}, we rather 
take the bubble estimator approach~\cite{AnBuCo20, BaSm93, CoAnVaBu20} by 
constructing a space with so-called $C^1$ bubble functions of degree $p+1$, each supported on a single mesh element.  
Denoting 
the resulting bubble space by $\mathcal{B}_h$, the estimator is obtained as the unique function $e_h \in \mathcal{B}_h$, satisfying
{\new 
\[
    a(e_h,b_h)
    =
    F(b_h)-a(u_h,b_h),
    \qquad\text{for all } b_h\in B_h.
\]
}
Since each bubble function is supported on only one element, the resulting estimator is localized, and the contribution of each element is measured by the corresponding local energy norm. In this way, the estimator identifies those regions which still have to be refined. Thereby, the selection of the corresponding mesh elements is done in the module {\em Mark} using D\"orfler's strategy~\cite{Do96}. The marked elements are then refined in the module {\em Refine} to get a refined hierarchical mesh. In general, however, it is not enough just to refine the marked elements. Instead, the set of marked, and thus of refined, elements has to be slightly enlarged in order to satisfy two additional requirements as explained in the following two paragraphs. 

First, to preserve structural requirements of the mesh by preventing the construction of very irregular hierarchical meshes, we generate hierarchical meshes that are admissible~\cite{BuGi16}. Note that a hierarchical mesh is called admissible of class $\mu > 1$ if for every active element $Q \in \mathcal{Q}^\ell$ of the hierarchical mesh $\mathcal{Q}$ (consisting of different levels $\mathcal{Q}^\ell$) all basis functions that have a nonzero support on $Q$ belong to at most $\mu$ successive refinement levels. That means the admissibility condition controls how many refinement levels may interact on each active element and therefore limits the level gap among basis functions overlapping the same element. Whenever we refine an element $Q$, we 
must also refine its neighborhood
$$
N(Q,\mu) = \{ Q' \in \mathcal{Q}^{\ell-\mu+1}: \; Q' \cap \mathcal{E}(Q,\ell-\mu+1+\chi) \neq \emptyset \},
$$
where $\chi=1$ 
for the case of using the truncated hierarchical 
construction
and $\chi=0$ otherwise, and where
$$
\mathcal{E}(Q,k) = \bigcup 
\{ \, \mbox{supp } \psi:\; \psi \in \mathcal{N}^{p,r}_k(\mathcal{D}^{(2)}) \mbox{ and } \mbox{supp }\psi \cap Q  \neq \emptyset\}
$$
denotes the union of supports of all level $k$ basis functions that do not vanish on $Q$. 
Since 
the truncated hierarchical construction is based on a smaller support extension $\mathcal{E}$ compared to the non-truncated construction, the use of a truncated hierarchical basis generally requires the refinement of a smaller surrounding region to preserve the same admissibility class~$\mu$ as in case of a non-truncated hierarchical basis.   

Second, to assure in case of a non-truncated hierarchical construction that the imposition of the $C^2$-continuity conditions at the inner vertices~$\bfm{\Xi}^{(i)}$, $i \in \mathcal{I}_{\Xi}^{I}$, in the IETI-DP method works as described in Section~\ref{sec:IETI}, and hence that the assembly of the constraint matrix~$\ab{C}_{\Xi,\Pi}$ can be done as presented in Section~\ref{subsec:constraint_CXi}, we have to refine (but only) in case of $p-r=2$ some additional elements close to an inner vertex~$\bfm{\Xi}^{(i)}$, $i \in \mathcal{I}_{\Xi}^{I}$. More precisely, whenever an element adjacent to an inner vertex~$\bfm{\Xi}^{(i)}$, $i \in \mathcal{I}_{\Xi}^{I}$, has to be
refined, then not just this one element needs to be refined but also an additional element in both directions 
within the same patch, namely the neighboring element in each direction. This is due to the fact that when refining the element next to an inner vertex~$\bfm{\Xi}^{(i)}$, $i \in \mathcal{I}_{\Xi}^{I}$, the third coarse function in each direction associated to the boundary ring of basis functions of a patch has already a support partially outside of this (refined) element, and hence stays active. After adding finer functions we all together have eight functions in the patch, namely two from the coarse level and six from the finer level, which are involved in the $C^2$ condition at the vertex~$\bfm{\Xi}^{(i)}$. Since in each patch adjacent to an inner vertex~$\bfm{\Xi}^{(i)}$, $i \in \mathcal{I}_{\Xi}^{I}$, always only exactly six basis functions should be involved in the $C^2$ condition at the vertex~$\bfm{\Xi}^{(i)}$, we have to refine one further element in each direction to fulfill this property. However, this does not happen, i.e. no additional refinement is needed, if we consider the truncated case. Note that the potentially additional refinement of elements in the vicinity of an inner vertex in the non-truncated case is similar to the $C^1$ adaptive method~\cite{BrGiKaVa23}, where even up to two additional elements in both directions have to be refined but regardless of the truncated or non-truncated case. 

\subsection{Computing the coefficients of $\ab{u}_B$} \label{subsec:coefficients_ub}

{\new
The coefficient vector~$\ab{u}_B$ contains only those
coefficients that are involved in the essential boundary conditions~\eqref{eq:constraintsB}.
Hence, coefficients controlling the trace of $u_h$ on $\Gamma_C\cup\Gamma_S$ and coefficients required to prescribe
$\partial_{\boldsymbol n}u_h$ on $\Gamma_C$ belong to
$\ab{u}_B$. In contrast, coefficients associated solely
with the free boundary~$\Gamma_F$ are not contained in $\ab{u}_B$. 
We 
denote by $\mathcal I_{\Omega}^{B,E}$ the index set of boundary patches containing an essential boundary segment, i.e.
\[
 \mathcal I_{\Omega}^{B,E}
 =\left\{i\in\mathcal I_{\Omega}^{B}:
 \partial\Omega^{(i)}\cap {\rm int}_{\partial \Omega}(\Gamma_C\cup\Gamma_S) \neq \emptyset \right\},
\]
and introduce for every $i\in\mathcal I_{\Omega}^{B,E}$ the local index set
\[
 \mathcal J_{B}^{(i)}
 =\left\{j\in\mathcal J_{\mathcal H}^{(i)}:
 \phi_j^{(i)}\big|_{(\Gamma_C\cup\Gamma_S)\cap\partial\Omega^{(i)}}\not\equiv0
 \ \; \text{or}\ \; 
 \partial_{\boldsymbol n}\phi_j^{(i)}\big|_{\Gamma_C\cap\partial\Omega^{(i)}}\not\equiv0
 \right\}.
\]
Thus, along a clamped boundary segment, $\mathcal J_B^{(i)}$ contains the two boundary
rings needed to prescribe the trace and the normal derivative there, whereas along a simply
supported boundary segment only the coefficients needed to prescribe the trace ({i.e. one boundary ring}) are
included. 
}
Following the method in~\cite{KaKoVi26}, we compute the coefficient vector~$\ab{u}_B$
by solving an IETI-type saddle point problem
\begin{equation}
\label{eq:Boundary_IETI_problem}
\left(
\begin{array}{cc}
{\new \MM} & \ab{B}^T \\
\ab{B} & \ab{0}
\end{array}
\right)
\left(
\begin{array}{c}
\ab{u}_B \\
\ab{\mu}
\end{array}
\right)
=
\left(
\begin{array}{c}
\ab{b} \\
\ab{0}
\end{array}
\right),
\end{equation}
with the matrix ${\new \MM} $ 
of the form
$
{\new \MM}  = \Diag (\{ {\new \MM^{(i)}}\}_{i \in \mathcal{I}^{B,E}_{\Omega}}),
$
where
{\new 
\begin{align*}
{\new \MM^{(i)}}
= & \; 
\omega_0
\left[
\int_{(\Gamma_C\cup\Gamma_S)\cap\partial\Omega^{(i)}}
\phi_{j_1}^{(i)}(\boldsymbol x)\,\phi_{j_2}^{(i)}(\boldsymbol x)\,\mathrm ds
\right]_{j_1,j_2\in\mathcal J_B^{(i)}}
\\
& +  \omega_1
\left[
\int_{\Gamma_C \cap\partial\Omega^{(i)}}
\partial_{\boldsymbol n}\phi_{j_1}^{(i)}(\boldsymbol x)\,
\partial_{\boldsymbol n}\phi_{j_2}^{(i)}(\boldsymbol x)\,\mathrm ds
\right]_{j_1,j_2\in\mathcal J_B^{(i)}},
\end{align*}
}
with the 
load vector $\ab{b}$ of the form
$
 \ab b=[\ab b^{(i)}]_{i\in\mathcal I_{\Omega}^{B,E}},
$
where
{\new
\begin{align*}
\ab b^{(i)}
=
\omega_0
\left[
\int_{(\Gamma_C\cup\Gamma_S)\cap\partial\Omega^{(i)}}
 g_0(\boldsymbol x)\,\phi_j^{(i)}(\boldsymbol x)\,\mathrm ds
\right]_{j\in\mathcal J_B^{(i)}}
+
\omega_1
\left[
\int_{\Gamma_C \cap\partial\Omega^{(i)}}
 g_1(\boldsymbol x)\,\partial_{\boldsymbol n}\phi_j^{(i)}(\boldsymbol x)\,\mathrm ds
\right]_{j\in\mathcal J_B^{(i)}},
\end{align*}
}
and with the full row rank constraint matrix~$\ab B$ of the form
\[
 \ab B=[\ab B^{(i)}]_{i\in\mathcal I_{\Xi}^{B,E,\nu\geq2}},
\]
where each submatrix~$\ab B^{(i)}$ represents in its rows the smoothness conditions at a
boundary vertex~$\boldsymbol\Xi^{(i)}$ 
{\new that are required for the prescribed essential
boundary coefficients in $\ab u_B$ to be compatible with the $C^1$ coupling across the
adjacent inner edges.} 
In the matrices~${\new \MM^{(i)}}$ and vectors~$\ab{b}^{(i)}$ above, 
$\omega_{0}$ and $\omega_1$ are suitable non-negative weights,
and 
{\new 
\[
 \mathcal I_{\Xi}^{B,E,\nu\geq2}
 =\left\{i\in\mathcal I_{\Xi}^{B,\nu\geq2}:
 \boldsymbol\Xi^{(i)}\in \overline{\Gamma_C\cup\Gamma_S}\right\}.
\]
}

Since the matrix ${\new \MM}$ is invertible, the linear system~\eqref{eq:Boundary_IETI_problem}
can then be solved by first determining the Lagrange multipliers~$\ab{\mu}$ from the symmetric
positive definite system
\[
\left( \ab{B} {\new \MM^{-1}} \ab{B}^T \right) \ab{\mu}
=
\ab{B} {\new \MM^{-1}} \ab{b},
\]
followed by computing $\ab{u}_B$ via
\[
\ab{u}_B =
{\new \MM^{-1}} \left( \ab{b} - \ab{B}^T \ab{\mu} \right).
\]
{\new
The coefficient vector obtained in this way is precisely the prescribed boundary
coefficient vector denoted by $\ab g$ in~\eqref{eq:constraints_CB}, i.e., $\ab g=\ab u_B$
for the subsequent IETI-DP system.
}

It remains to explain the construction of the submatrix~$\ab{B}^{(i)}$ for a boundary vertex~$\bfm{\Xi}^{(i)}$ with patch valency~$\nu_i \geq 2$, i.e., for $i \in {\new \mathcal{I}_{\Xi}^{B,E,\nu \geq 2}}$. We first generate a matrix~$\widetilde{\ab{B}}^{(i)}$, which collects in its rows certain smoothness conditions at the boundary vertex~$\bfm{\Xi}^{(i)}$. More precisely, for any two neighboring patches~$\Omega^{(i_j)}$ and $\Omega^{(i_k)}$, $i_j,i_k \in \mathcal{I}_{\Omega}^B$, $i_j \neq i_k$, with $ \overline{\Omega^{(i_j)}} \cap  \overline{\Omega^{(i_k)}}= \overline{\Gamma^{(i_m)}}$ for some $i_m \in \mathcal{I}_{\Gamma}^{I}$ such that $\bfm{\Xi}^{(i)} \in \overline{\Gamma^{(i_m)}}$, assuming that the inner edge~$\overline{\Gamma^{(i_m)}}$ is parameterized as in~\eqref{eq:gC}, i.e.,
$
\ab{G}^{(i_j)}(0,\xi) = \ab{G}^{(i_{k})}(0,\xi), \ \xi \in [0,1],
$
and that the boundary vertex~$\bfm{\Xi}^{(i)}$ is given by
$
\bfm{\Xi}^{(i)} = \ab{G}^{(i_j)}(0,0) = \ab{G}^{(i_{k})}(0,0)$, we add to the rows of the matrix~$\widetilde{\ab{B}}^{(i)}$ the smoothness conditions 
\begin{align} \label{eq:boundary_vertex}
  & \g_0^{(i_m,i_j)}[u_h](0)   =  \g_0^{(i_m,i_{k})}[u_h](0), \quad
  {\partial \g_0^{(i_m,i_j)}[u_h](0)  =  \partial \g_0^{(i_m,i_{k})}[u_h](0)}, \nonumber \\
  & {\partial^2 \g_0^{(i_m,i_j)}[u_h](0)   =  \partial^2 \g_0^{(i_m,i_{k})}[u_h](0)},\\
  & \g_1^{(i_m,i_j)}[u_h](0)   =  \g_1^{(i_m,i_{k})}[u_h](0), \quad
  {\partial \g_1^{(i_m,i_j)}[u_h](0)  =  \partial \g_1^{(i_m,i_{k})}[u_h](0)}, \nonumber
\end{align}
with respect to the coefficient vector~$\left( \widetilde{\ab{u}}_E^T, \, \ab{u}_B^T\right)^T$, where  $\widetilde{\ab{u}}_E$ collects the coefficients of $\ab{u}_E$ which are involved in the smoothness conditions~\eqref{eq:boundary_vertex}. The obtained matrix~$\widetilde{\ab{B}}^{(i)}$ is then used to construct the matrix~$\ab{B}^{(i)}$ via
\begin{equation}
\label{eq:boundary_vertexOutside}
\left(
\begin{array}{cc}
\ab{0} & \ab{B}^{(i)}
\end{array}
\right)
=
(\ab{N}^{(i)})^T \widetilde{\ab{B}}^{(i)},
\end{equation}
with the matrix $\ab{N}^{(i)}$ describing a basis of the null space of the matrix $(\overline{\ab{B}}^{(i)})^T$, where $\overline{\ab{B}}^{(i)}$ denotes the matrix obtained from $\widetilde{\ab{B}}^{(i)}$ by setting all columns corresponding to the coefficient vector~${\ab{u}}_B$ to zero, {and with} the zero matrix 
representing the resulting zero columns with respect to the coefficient vector~$\widetilde{\ab{u}}_E$.

\subsection{Constructing matrices $\ab{C}_{\Xi,\Pi}$ and $\ab{N}_{\Xi,\Pi}$} \label{subsec:constraint_CXi}

We construct the matrix ${\ab{C}}_{\Xi,\Pi}$ as 
\begin{equation}
\label{eq:C_XiF}
\ab{C}_{\Xi,\Pi} = [\ab{C}^{(i)}_{\Xi,\Pi}]_{i \in \mathcal{I}_{\Xi}^{I}},
\end{equation}
where each submatrix~$\ab{C}^{(i)}_{\Xi,\Pi}$ collects 
in its rows the $C^{2}$-smoothness conditions~\eqref{eq:constraintsV} at the inner vertex~$\bfm{\Xi}^{(i)}$, $i \in \mathcal{I}_{\Xi}^{I}$, with respect to the coefficient vector~$\ab{u}_{\Pi}$.
Assuming that the neighboring patches~$\Omega^{(i_0)}, \ldots, \Omega^{(i_{\nu_i-1})}$ 
with
$
\bfm{\Xi}^{(i)} = \cap_{\ell=0}^{\nu_i-1} \overline{\Omega^{(i_{\ell})}}
$
are 
arranged in counterclockwise order, then the rows of the matrix~$\ab{C}^{(i)}_{\Xi,\Pi}$ are simply determined 
by the
$
6(\nu_i-1)
$
linearly independent $C^{2}$-smoothness conditions 
\begin{equation}
\label{eq:inner_vertex}
\left(\ab{u}^{(i_{\ell})}\right)^T \partial_1^{j_1} \partial_2^{j_2} \ab{\phi}^{(i_\ell)}(\ab{\Xi}^{(i)})
=
\left(\ab{u}^{(i_{\ell+1})}\right)^T \partial_1^{j_1} \partial_2^{j_2} \ab{\phi}^{(i_{\ell+1})}(\ab{\Xi}^{(i)}),
\quad
j_1 + j_2 \leq 2,
\quad
\ell = 0,\ldots,\nu_i-2,
\end{equation}
at the inner vertex $\ab{\Xi}^{(i)}$, cf.~\cite{KaKoVi26}.

Our next goal is to efficiently construct the matrix ${\ab{N}}_{\Xi,\Pi}$ from ${\ab{C}}_{\Xi,\Pi}$. To this end, we first define for each $i \in \mathcal{I}_{\Xi}^{I}$ the matrix 
$
\widehat{\ab{C}}^{(i)}_{\Xi,\Pi} \in \R^{6(\nu_i-1) \times 6\nu_i},
$
which 
is obtained from the 
matrix ${\ab{C}}^{(i)}_{\Xi,\Pi}$ 
by removing all zero columns.
Next, let $\ab{N}^{(i)}_{\Xi,\Pi}  \in \R^{6 \nu_i \times 6}$ 
be the matrix representing a basis of the 
null space 
of the matrix~$\widehat{\ab{C}}^{(i)}_{\Xi,\Pi}$. Furthermore, let ${\ab{N}}^{(i)}_{\Xi,\Pi} \in \R^{|\ab{u}_\Pi| \times 6}$ be the matrix obtained from $\widehat{\ab{N}}^{(i)}_{\Xi,\Pi}$ by extending it with zero rows at the positions corresponding to the zero columns removed when passing from ${\ab{C}}^{(i)}_{\Xi,\Pi}$ to $\widehat{\ab{C}}^{(i)}_{\Xi,\Pi}$.
Finally, we obtain ${\ab{N}}_{\Xi,\Pi}$ as
\[
{\ab{N}}_{\Xi,\Pi}
=
\left(
\begin{array}{cccc}
{{\ab{N}}^{(1)}_{\Xi,\Pi}} &
{{\ab{N}}^{(2)}_{\Xi,\Pi}} &
\cdots &
{{\ab{N}}^{(|\mathcal{I}_\Xi^I|)}_{\Xi,\Pi}}
\end{array}
\right)
\in
\R^{|\ab{u}_\Pi| \times 6 |\mathcal{I}_\Xi^I|}.
\]

\subsection{Constructing matrix $\ab{C}_{\Gamma}$} \label{subsec:constraint_CGamma}

We adapt the construction from~\cite{KaKoVi26} to the case of {\new the Kirchhoff--Love plate problem and} to adaptive refinement. Recall that
${\ab{C}}_{\Gamma} = \left( \begin{array}{cccc}
{\ab{C}}_{\Gamma,B} & {\ab{C}}_{\Gamma,\Pi} & {\ab{C}}_{\Gamma,E} & \ab{0} \end{array}\right)$ with respect to the coefficient vector~$\ab{u}=\left( \ab{u}_B^T, \, \ab{u}_\Pi^T , \, \ab{u}_E^T, \, \ab{u}_R^T \right)^T$.
We generate the matrix $\ab{C}_{\Gamma}$ as 
\[
[\ab{C}_{\Gamma}^{(i)}]_{i \in \mathcal{I}_{\Gamma}^{I}},
\]
where the rows of the 
submatrix $\ab{C}_{\Gamma}^{(i)}$,  $i \in \mathcal{I}_\Gamma^I$, specify 
the $C^1$-smoothness conditions~\eqref{eq:constraintsI1} and~\eqref{eq:constraintsI2} for the inner edge~$\Gamma^{(i)}$ with respect to the coefficient vector~
$\ab{u}$. Thereby, these conditions have to be linearly independent with the smoothness conditions~\eqref{eq:constraintsV} and/or with the boundary conditions~\eqref{eq:constraintsB} at the two endpoints of $\overline{\Gamma^{(i)}}$, which are already incorporated in the matrices $\ab{C}_{\Xi,\Pi}$ and/or $\ab{C}_B$. 
{\new 
Note that $\ab{C}_B$ contains endpoint conditions only for boundary vertices in $\overline{\Gamma_C\cup\Gamma_S}$.
} 
Below, we
assume that for the inner edge~$\Gamma^{(i)}$, $i \in \mathcal{I}_{\Gamma}^{I}$, the two neighboring patches~$\Omega^{(i_0)}$ and $\Omega^{(i_1)}$, $i_0,i_1 \in \mathcal{I}_{\Omega}$, satisfying
$
\overline{\Gamma^{(i)}} = \overline{\Omega^{(i_0)}} \cap \overline{\Omega^{(i_1)}},
$
are parameterized as in~\eqref{eq:standard_par}. Before we explain the construction of the matrix~$\ab{C}_{\Gamma}^{(i)}$, $i \in \mathcal{I}_\Gamma^I$, we introduce some required tools. 

Let $\mathcal{D}^{\tau,(1)} = ( \mathcal{D}^{\tau,(1)}_\ell )_{\ell=0}^{N-1}$, $\tau \in \{i_0,i_1\}$, be two sequences of nested subdomains in $[0,1]$ which are mapped 
to
$
\overline{\Gamma^{(i)}} \subset \overline{\Omega^{(\tau)}}, 
$
using $\ab{G}^{(\tau)}$. 
In addition, we introduce another sequence of nested subdomains in $[0,1]$ 
given as 
\[
\mathcal{D}^{i,(1)} = ( \mathcal{D}^{i,(1)}_\ell )_{\ell=0}^{N-1},
\qquad
\mathcal{D}^{i,(1)}_\ell = \mathcal{D}^{i_0,(1)}_\ell \cap \mathcal{D}^{i_1,(1)}_\ell,
\]
to determine via Eq.~\eqref{eq:gC} a well-defined
trace function $\g_0^{(i)}[u_h] \in \mathcal{H}^{p,r+1}(\mathcal{D}^{i,(1)})$ and 
specific first derivative function $\g_1^{(i)}[u_h] \in \mathcal{H}^{p-1,r}(\mathcal{D}^{i,(1)})$  possessing the spline representations 
$$
\g_0^{(i)}[u_h] = \hspace{-1cm} \sum_{\varphi_j \in 
\cup_{\ell=0}^{N-1}\mathcal{N}^{p,r+1}_{\ell} (\mathcal{D}^{i,(1)})}
 \hspace{-1cm} d_{j,0}^{(i)} \, \varphi_{j} \quad \mbox{and} \quad
 \g_1^{(i)}[u_h] = \hspace{-1cm} \sum_{\varphi_j \in 
\cup_{\ell=0}^{N-1}\mathcal{N}^{p-1,r}_{\ell} (\mathcal{D}^{i,(1)})}
 \hspace{-1cm} d_{j,1}^{(i)} \, \varphi_{j},
$$
respectively, with the coefficients~$d_{j,0}^{(i)}$ and $d_{j,1}^{(i)}$, collected in the coefficient vector
\[
\ab{d}^{(i)} = \big[d_{j,\ell}^{(i)}\big]_{j=0, \ldots, \widetilde{n}_{\mathcal{H},\ell}^{i,(d)}-1  ;\ell=0,1}, \quad \widetilde{n}_{\mathcal{H},\ell}^{i,(d)}  = \dim \mathcal{H}^{p-\ell,r+1-\ell}(\mathcal{D}^{i,(1)}).
\]

We are now ready to explain the construction of the 
matrix~$\ab{C}_{\Gamma}^{(i)}$ for the inner edge~$\Gamma^{(i)}$, $i \in \mathcal{I}_{\Gamma}^{I}$, for which we have to distinguish two cases. In the first case, we assume 
that $\overline{\Gamma^{(i)}} \cap \partial \Omega = \emptyset$, i.e., both endpoints of $\overline{\Gamma^{(i)}}$ are inner vertices, and let us construct linearly independent smoothness conditions~\eqref{eq:constraintsI1} and~\eqref{eq:constraintsI2}, reduced by the ones which are already included in the matrix~$\ab{C}_{\Xi,\Pi}$. 
Let $\xi^{p,r}_{\ell,j}$, $j=0,\ldots,n_{h_\ell}-1$, denote the Greville abscissae for the space $\mathcal{S}_{h_\ell}^{p,r}([0,1])$ corresponding to the basis functions $N_{h_\ell,j}^{p,r}$.
We aim to construct mixed hierarchical Greville abscissae $\zeta_{\tau,\tilde{j}}^{p,r}$ for the hierarchical spaces $\mathcal{H}^{p,r}(\mathcal{D}^{\tau,(1)})$, 
$\tau \in \{i_0,i_1\}$. We associate to every active hierarchical B-spline function its corresponding Greville abscissa from the level where that function is defined, i.e.,
\begin{equation*}  \label{eq:set1Greville}
\left\{ \zeta_{\tau,\tilde{j}}^{p,r}: \; \tilde{j}=0,1,\ldots, n_\mathcal{H}^{\tau,(1)}-1 \right\} = \bigcup_{\ell=0}^{N-1} \left\{ 
\xi^{p,r}_{\ell,j}: \; j\in \{0,\ldots,n_{h_\ell}-1\}, \; N^{p,r}_{h_\ell, j} \in \mathcal{N}_\ell^{p,r}(\mathcal{D}^{\tau,(1)}) 
\right\}.
\end{equation*}
The smoothness conditions can now be expressed as  
\begin{equation}
\label{eq:interfaceSmoothnessConditions_innerVertex}
\begin{array}{ll}
\g_\rho^{(i,\LL)}[u_h](\zeta^{p,r}_{i_0,\tilde{j}}) = \gC^{(i)}_\rho[u_h](\zeta^{p,r}_{i_0,\tilde{j}}), &
\tilde{j} = 0,\ldots,
n_\mathcal{H}^{i_0,(1)} - 1,
\\[0.5ex]
\g_\rho^{(i,\RR)}[u_h](\zeta^{p,r}_{i_1,\tilde{j}}) = \gC^{(i)}_\rho[u_h](\zeta^{p,r}_{i_1,\tilde{j}}), &
\tilde{j} = 3-\rho,\ldots, 
n_\mathcal{H}^{i_1,(1)} - 4 + \rho,
\end{array}
\quad \rho = 0,1.
\end{equation}
However, active functions from different levels may share the same Greville abscissa. Additionally, in the transition region between coarse and fine levels the corresponding Greville abscisae may be badly graded leading to instability and badly conditioned smoothness conditions \eqref{eq:interfaceSmoothnessConditions_innerVertex}. 
Consequently we have to replace some equations in \eqref{eq:interfaceSmoothnessConditions_innerVertex} with another  ones to ensure well conditioned linearly independent smoothness condition. We follow here in the next paragraph the idea of weighted collocation presented in \cite{WeightedCollocation2013}.

We say that a function $ N_{h_\ell,j}^{p,r} \in \mathcal{N}^{p,r}_{\ell}(\mathcal{D}^{\tau,(1)})$ is in a transition region if its support has a nonempty intersection with the refinement subdomain $\mathcal{D}_{\ell+1}^{\tau,(1)}$. Then we replace the equations in \eqref{eq:interfaceSmoothnessConditions_innerVertex}, which are evaluated at the corresponding Greville abscissa 
$\zeta_{\tau,\tilde{j}}^{p,r} = \xi_{\ell,j}^{p,r}$, with new ones. For this, let us first represent the function $ N_{h_\ell,j}^{p,r} \in \mathcal{N}^{p,r}_{\ell}(\mathcal{D}^{\tau,(1)})$ with respect to finer functions from the space $\mathcal{S}_{h_{\ell+1}}^{p,r}([0,1])$ as
$$
 N_{h_\ell,j}^{p,r}  = \sum_{\iota_j} c_{\iota_j}^{(\ell+1)}  N_{h_{\ell+1},\iota_j}^{p,r},
$$
with the corresponding weights $c_{\iota_j}^{(\ell+1)}$,
and let $\xi_{\ell+1,\iota_j}^{p,r}$ be the corresponding Greville absissae of the functions $N_{h_{\ell+1},\iota_j}^{p,r}$. Then, the equations 
\begin{equation*}
\g_\rho^{(i,\tau)}[u_h](\xi_{\ell,j}^{p,r}) =  
\gC^{(i)}_\rho[u_h](\xi_{\ell,j}^{p,r}), \quad \rho=0,1,
\end{equation*}
are replaced with the equations
\begin{equation*}
  \sum_{\iota_j} c_{\iota_j} 
\g_\rho^{(i,\tau)}[u_h](\xi_{\ell+1,\iota_j}^{p,r}) =  \sum_{\iota_j} c_{\iota_j}  \gC^{(i)}_\rho[u_h](\xi_{\ell+1,\iota_j}^{p,r}), \quad \rho=0,1.
\end{equation*}
The obtained linearly independent smoothness conditions can now be 
collected in matrix form as
\begin{equation}  \label{eq:ui0ui1di_coefficients}
\left( \begin{array}{ccc}
\ab{C}^{(i,\LL)} &  \ab{0} & \ab{D}^{(i,\LL)} \\
\ab{0} & \ab{C}^{(i,\RR)} & \ab{D}^{(i,\RR)}
\end{array} \right)
\left( 
\begin{array}{c}
\ab{u}_P^{(\LL)} \\  
\ab{u}_P^{(\RR)}\\
\ab{d}^{(i)}
\end{array}
\right) = 
\left(
\begin{array}{c}
\ab{0}\\
\ab{0}
\end{array}
\right),
\qquad
\ab{u}_P^{(\tau)} = \left( 
\begin{array}{c}
\ab{u}_B^{(\tau)} \\  
\ab{u}_\Pi^{(\tau)}\\
\ab{u}_E^{(\tau)}
\end{array}
\right) ,
\quad 
\tau \in \{\LL,\RR\}.
\end{equation}
By 
combining the two coefficient vectors $\ab{u}_P^{(\LL)}$ and $\ab{u}_P^{(\RR)}$ into the vector
$
\ab{u}_P^{(\LL,\RR)} = \left(\begin{array}{c} \ab{u}_P^{(\LL)}\\ \ab{u}_P^{(\RR)} \end{array}\right),
$
and by appropriately reordering the conditions in \eqref{eq:ui0ui1di_coefficients}, we get a linear system 
\begin{equation}  
\label{eq:systemC1C2D1D2}
\left( \begin{array}{cc}
\ab{C}^{(i)}_{1} &  \ab{D}^{(i)}_1 \\
\ab{C}^{(i)}_2 & \ab{D}^{(i)}_2
\end{array} \right)
\left( 
\begin{array}{c}
\ab{u}_P^{(\LL,\RR)} \\  
\ab{d}^{(i)}
\end{array}
\right) = 
\left(
\begin{array}{c}
\ab{0}\\
\ab{0}
\end{array}
\right),    
\end{equation} 
where $\ab{D}^{(i)}_2$ is a full rank square 
matrix. From the second row of~\eqref{eq:systemC1C2D1D2}, we obtain
\begin{equation} \label{eq:vec_d}
\ab{d}^{(i)} = - (\ab{D}^{(i)}_2)^{-1} \ab{C}^{(i)}_2 \, \ab{u}_P^{(\LL,\RR)}.
\end{equation}
Inserting 
the expression~\eqref{eq:vec_d} into the first row of~\eqref{eq:systemC1C2D1D2} 
leads to
\[
\left(
\ab{C}_1^{(i)} - \ab{D}_1^{(i)} (\ab{D}^{(i)}_2)^{-1} \ab{C}^{(i)}_2
\right)
\ab{u}_P^{(\LL,\RR)} = \ab{0},
\]
which hence gives us conditions involving only the coefficients $\ab{u}_P^{(\LL)}$ and $\ab{u}_P^{(\RR)}$.
The matrix $\ab{C}^{(i)}_\Gamma$ is then given by the matrix 
$
\ab{C}_1^{(i)} - \ab{D}_1^{(i)} (\ab{D}^{(i)}_2)^{-1} \ab{C}_2^{(i)},
$
extended by zero columns at the positions corresponding to coefficients in 
 $\ab{u}$ that are not contained in $\ab{u}_P^{(\LL)}$ and $\ab{u}_P^{(\RR)}$. 

In the second case, we assume that one or both endpoints of $\overline{\Gamma^{(i)}}$ are boundary vertices of $\Omega$, and
consider in detail the case in which one endpoint is a boundary vertex, denoted by the vertex~$\ab{\Xi}^{(j)}$, $j \in \mathcal{I}_{\Xi}^B$, which satisfies 
$
\ab{\Xi}^{(j)} = \ab{G}^{(i_0)}(0,0) = \ab{G}^{(i_1)}(0,0).
$
Note that the remaining cases can be 
handled in a similar manner. 
We first construct as in the previous paragraph an initial matrix~$\ab{C}_{\Gamma}^{(i)}$. 
{\new 
If $j\in\mathcal{I}_{\Xi}^{B,E,\nu\geq 2}$ and $\overline{\Gamma^{(i)}}$ is the edge with the smallest index having $\ab{\Xi}^{(j)}$ as endpoint, then we add further smoothness conditions as additional rows to the initial matrix~$\ab{C}^{(i)}_{\Gamma}$, otherwise, no additional rows associated with this vertex are added to this matrix. These additional smoothness conditions represent the direct complement of the conditions~\eqref{eq:boundary_vertexOutside} with respect to the conditions~\eqref{eq:boundary_vertex}. 
}
The direct complement is computed iteratively as follows. We start with an empty set and successively add one of the conditions~\eqref{eq:boundary_vertex} to this set, but only if the current set with the newly added condition from~\eqref{eq:boundary_vertex} is linearly independent from the conditions~\eqref{eq:boundary_vertexOutside}, and only until enough conditions have been added, namely the number of conditions in~\eqref{eq:boundary_vertex} minus the number of conditions in~\eqref{eq:boundary_vertexOutside}. The resulting set of conditions specifies then the desired additional smoothness conditions, which are added as rows to the matrix~$\ab{C}_\Gamma^{(i)}$, extended with zeros at the positions corresponding to the remaining coefficients in 
$\ab{u}$.
{\new
If, instead,
$
\ab{\Xi}^{(j)}\in \overline{\Gamma_F}\setminus\overline{\Gamma_C\cup\Gamma_S},
$
then no corresponding conditions are contained in $\ab{C}_B$. In this case, if $\overline{\Gamma^{(i)}}$ is the edge with the smallest index having $\ab{\Xi}^{(j)}$ as endpoint, we add the full set of smoothness conditions~\eqref{eq:boundary_vertex} as additional rows to $\ab{C}^{(i)}_{\Gamma}$, otherwise, no additional rows associated with this vertex are added to this matrix. 
}

\section{Numerical examples} \label{sec:numerical_examples}

In this section, we will present some illustrative examples for solving the 
{\new Kirchhoff--Love plate problem} {\new on} several {\new planar} analysis-suitable $G^{1}$ multi-patch geometries, including bilinear and non-bilinear multi-patch domains, using the proposed adaptive IETI-DP method. We will analyze the convergence behavior for the numerical solution $u_h$ obtained by solving the saddle point problem \eqref{eq:large_problem} 
using the adaptive 
refinement procedure presented in Section~\ref{subsec:adaptive_ref}. 
In Examples~\ref{ex:Example1}, \ref{ex:Example2}, and \ref{ex:Example4}, we use degrees $p \in \{3,4,5\}$, and in Example~\ref{ex:Example3}, we restrict ourselves to degree 
$p=3$. 
In all examples we set the regularity to $r=p-2$, {\new which is the highest possible regularity ensuring $h$-refinability of the traces and first derivatives of the computed numerical solution~$u_h$ across the patch interfaces, cf. Section~\ref{subsec:Cs_space},} and use the D\"orfler marking parameter equal to $0.8$. For the adaptive computations, we employ both hierarchical and truncated hierarchical B-spline spaces and consider admissibility classes 
{\new $\mu=2$ and/or $\mu=3$}. Moreover, the functions {\new $f$, $g_0$, $g_1$, $g_2$ and $g_3$}, which define {\new the right-hand side and the boundary data} of the 
{\new Kirchhoff--Love plate problem~\eqref{eq:kirchhoff-strong}, are either derived from a prescribed exact solution~$u$ or explicitly specified. In the former case,} the resulting approximations $u_h$ are compared with the exact solution~$u$ by computing the relative errors in {\new a quantity equivalent to} the $H^2$-seminorm, {\new given by} 
{\new
\begin{equation} \label{eq:eqiuv2seminorms}
 e(u,u_h):= 
 \frac{\| \nabla^2 u- \nabla^2 u_h\|_{F,L^2}}{\| \nabla^2 u \|_{F,L^2}},
\end{equation}
where $\| \nabla^2 u \|_{F,L^2} := \| \nabla^2 u : \nabla^2 u  \|_{L^2}$,
}
with respect to the number of degrees of freedom ($N_{d}$), {\new whereas, when no exact solution is prescribed, we compare two successive discrete solutions {\new $u_h$ and $u_{h/2}$}, and use the corresponding $h$--$h/2$-type relative error} {\new given by $e(u_h,u_{h/2})$.}
{\new For brevity, we refer to 
$e(u,u_h)$ and $e(u_h,u_{h/2})$ 
simply as the $H^2$-seminorm error and the $H^2$ error estimate, respectively.} To compute the convergence order with respect to the number of degrees of freedom ($N_{d}$), we will use the formula 
$$
 \frac{\log \left(e^{(\ell-1)}/ 
 e^{(\ell)}\right)}{\log \left(N_{d}^{(\ell)} / 
 N_{d}^{(\ell-1)}\right)} ,
$$
where $e^{(\ell)}$ denotes {\new the corresponding relative quantity $e(u,u_h)$ 
or~$e(u_h,u_{h/2})$} with respect to refinement level $\ell$.
Note that the optimal convergence order in the $H^2$-seminorm (or regarding $H^2$-error estimates) with respect to the number of degrees of freedom $N_{d}$ is $(p-1)/2$, which means that the errors decay with an optimal rate of $\mathcal{O}(N_d^{-(p-1)/2})$. 

{\new In the following, we present four numerical examples illustrating different aspects of the proposed IETI-DP method.}

\begin{figure}[htb!]
    \centering
    \begin{tabular}{cc}
  \hspace{-0.8cm}  \includegraphics[scale=0.25]{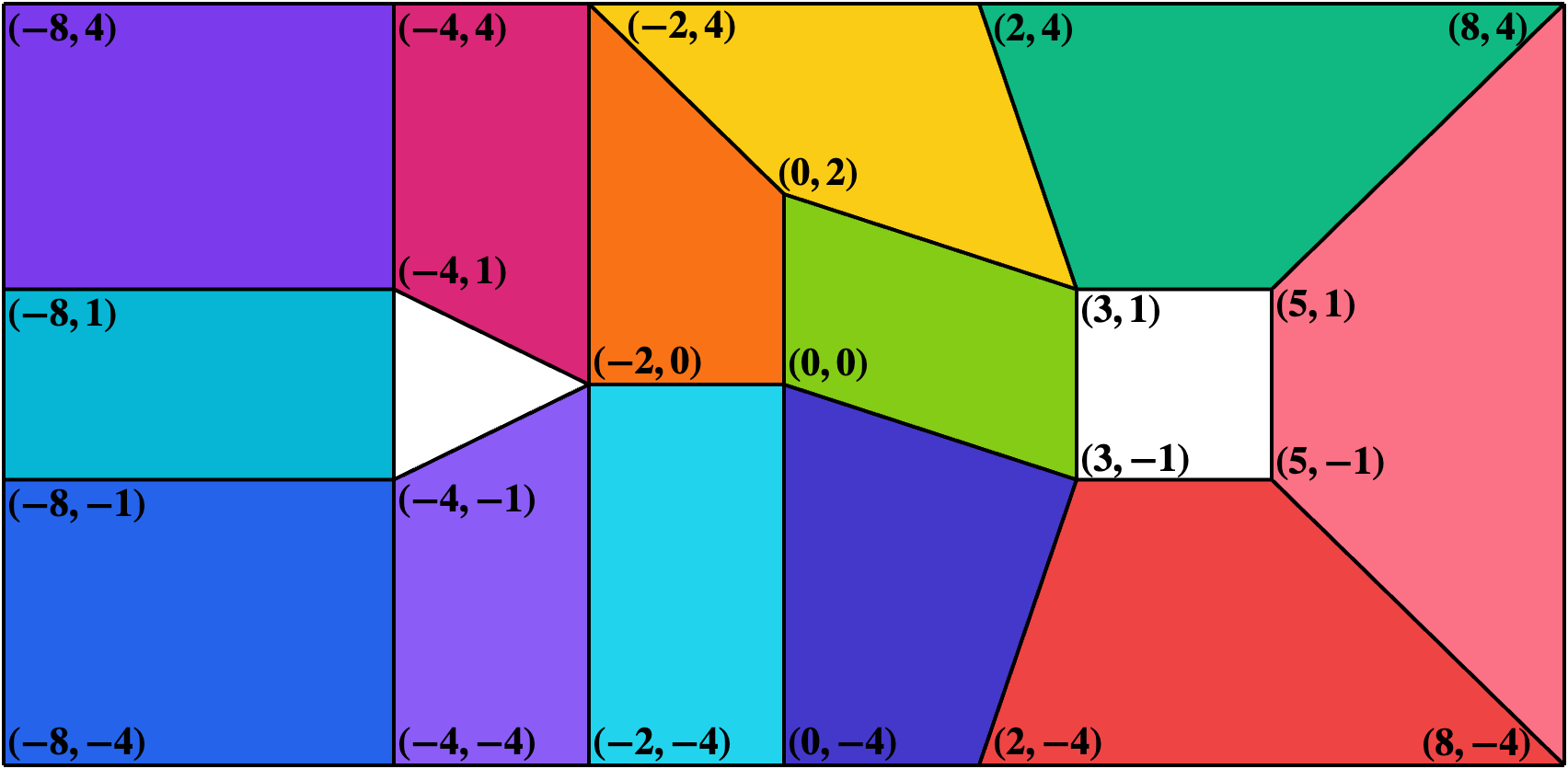}   & 
  \hspace{-1cm}
  \includegraphics[scale=0.14]{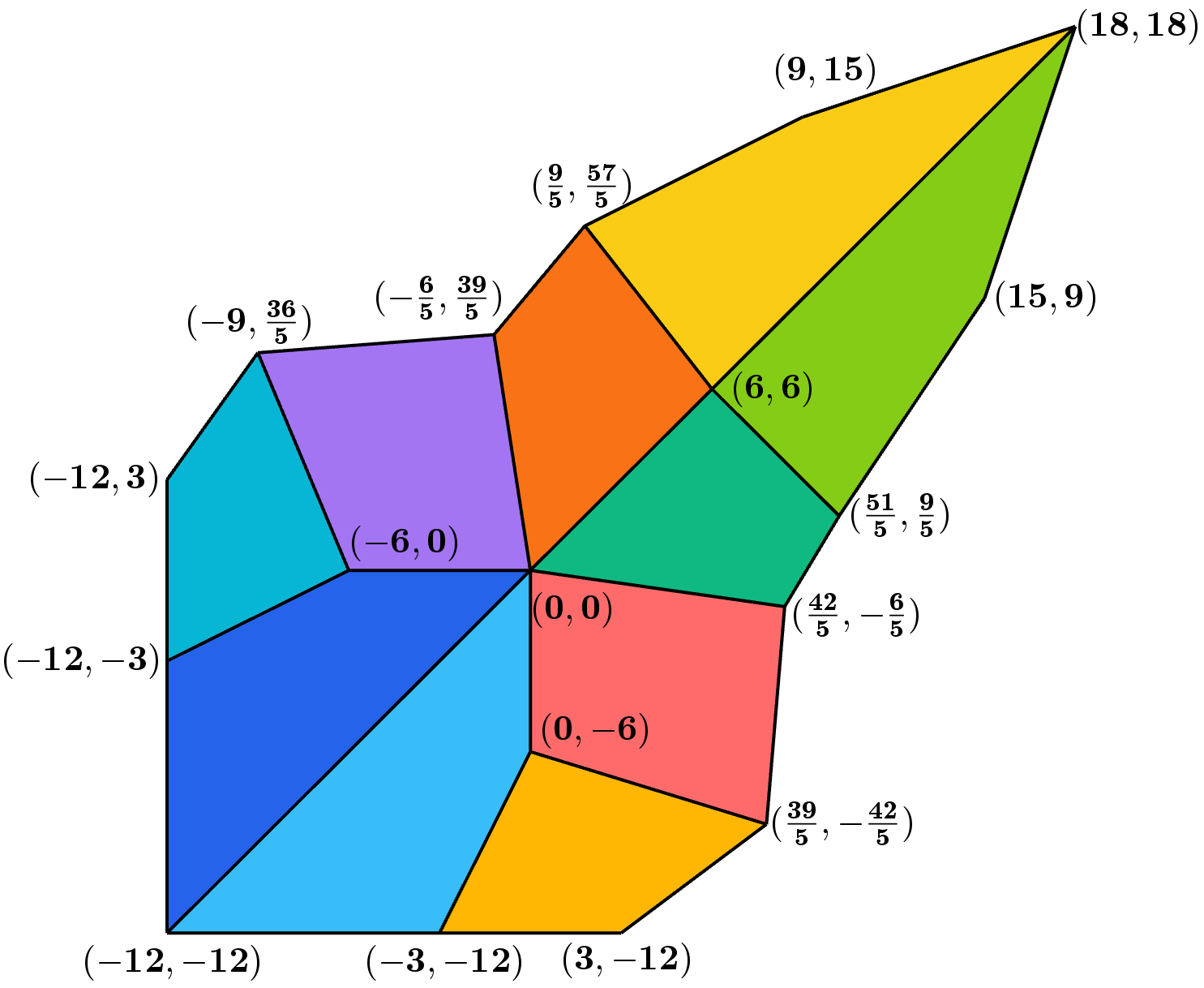} \\[0.2cm]
   Example~\ref{ex:Example1}  &  Example~\ref{ex:Example2} \\[0.8cm]
  \includegraphics[scale=0.115]{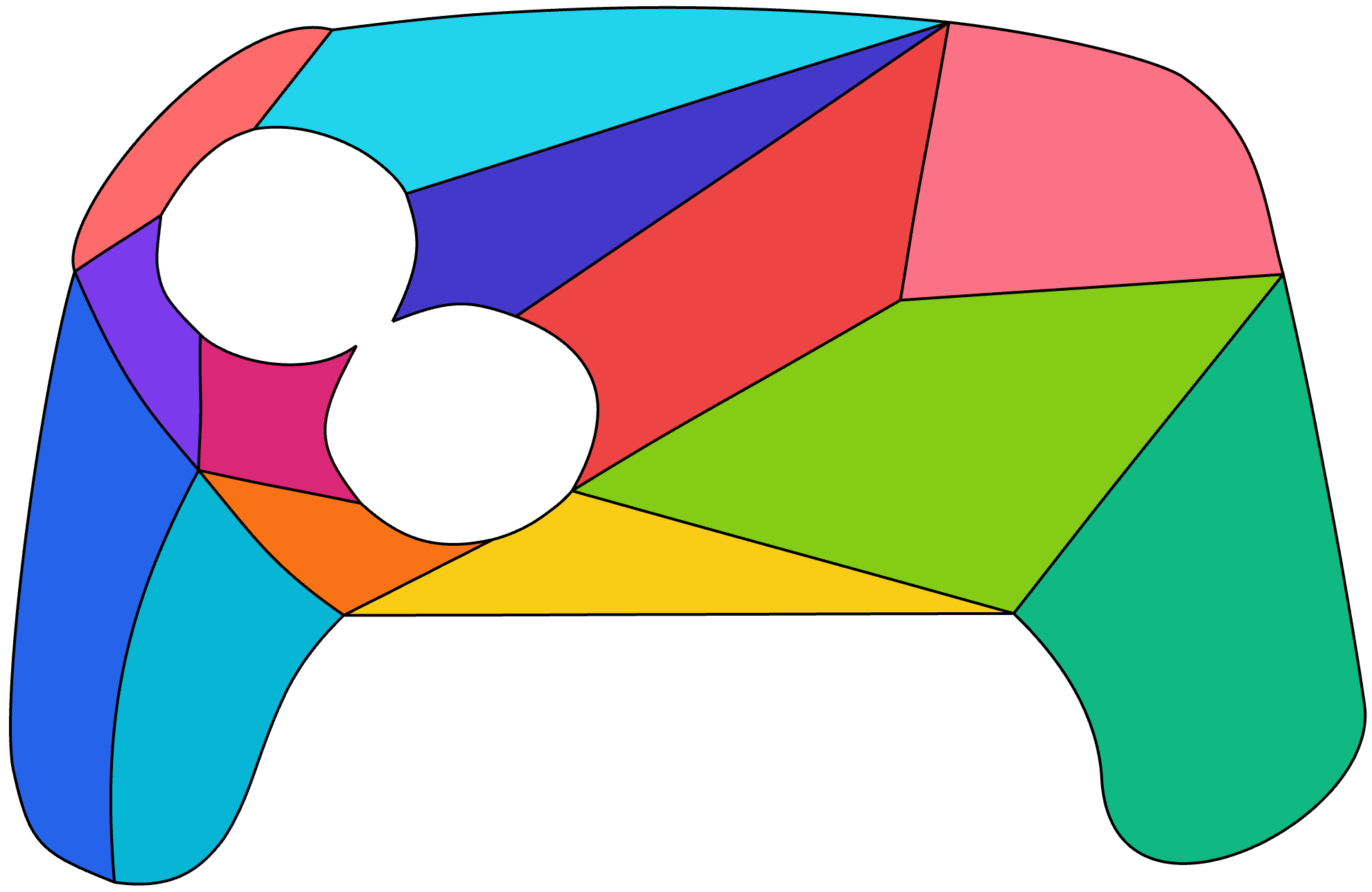}  \hskip3em &
  \includegraphics[scale=0.13]{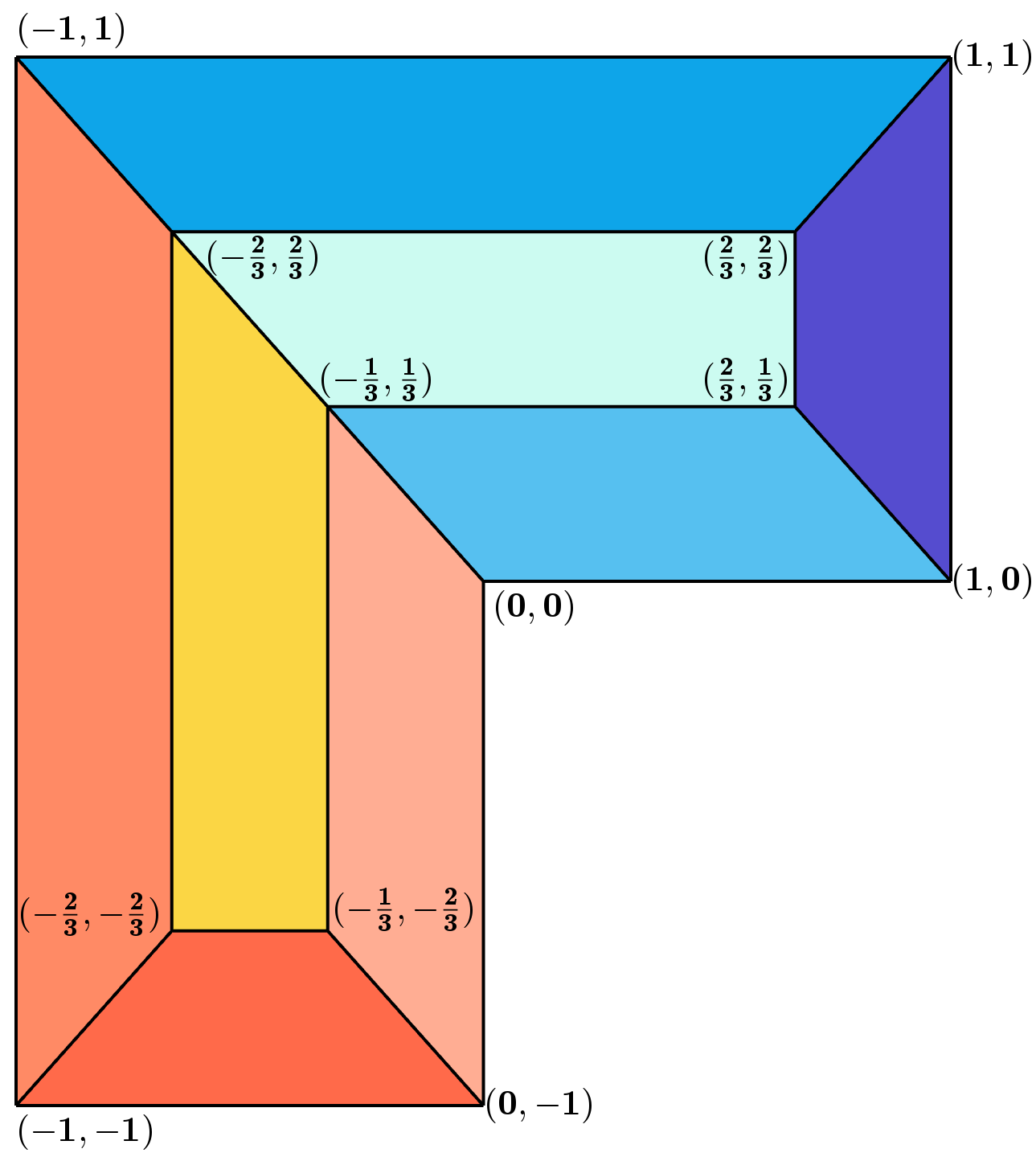}\\[0.4cm]
   {Example~\ref{ex:Example3}}  &  { Example~\ref{ex:Example4}} \\[0.4cm]
  \end{tabular}
\caption{The three bilinear multi-patch domains (Example~\ref{ex:Example1} (top left), Example~\ref{ex:Example2} (top right), and Example~\ref{ex:Example4} (bottom right)), and the {non-bilinear} thirteen-patch controller domain (Example~\ref{ex:Example3} (bottom left)).}
    \label{fig:domains}
\end{figure}

\begin{figure}[htb!]
    \centering
    \begin{tabular}{cc}
  \hspace{-1cm} \includegraphics[scale=0.33]{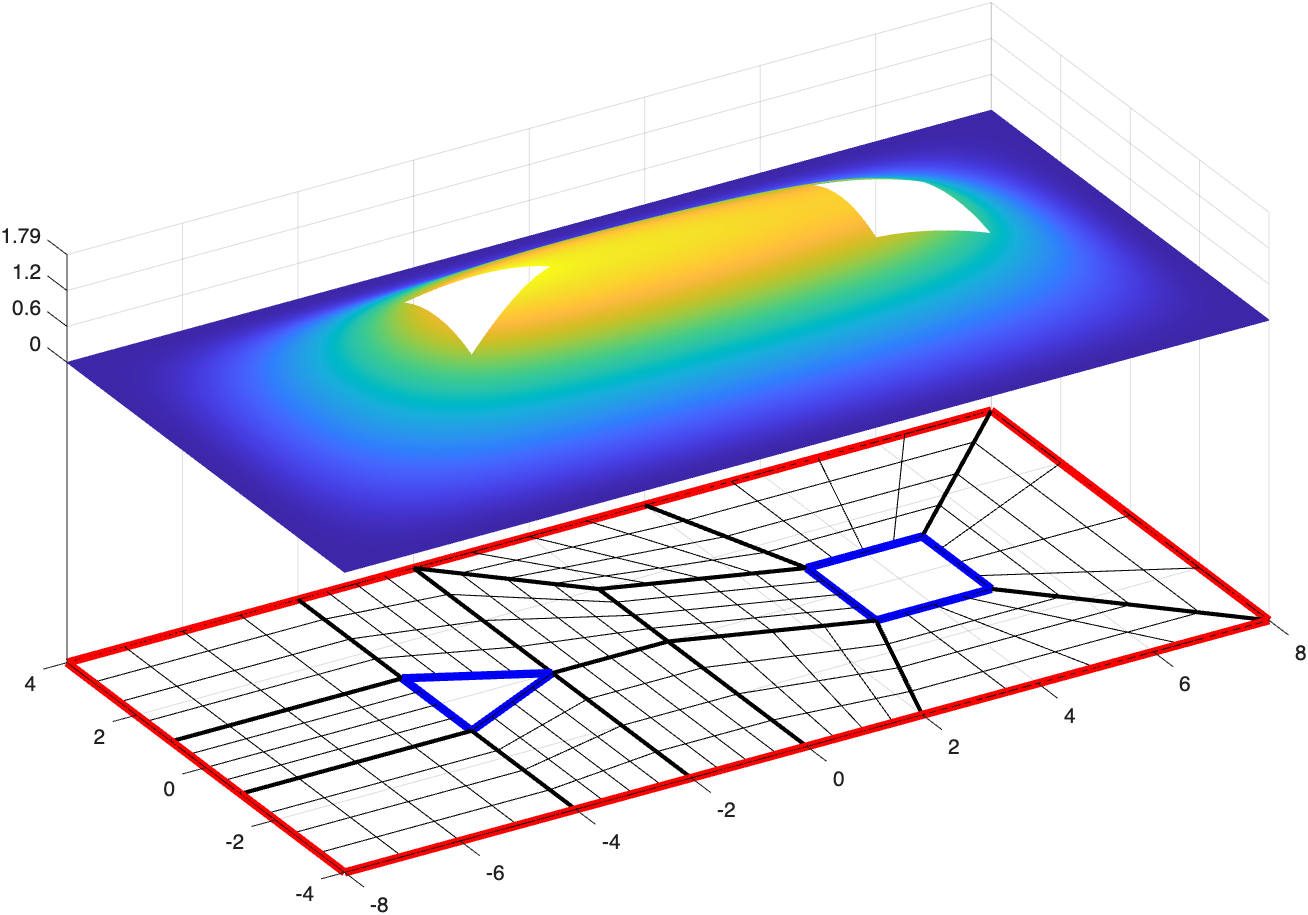}   & \hspace{-0.cm}
  \includegraphics[scale=0.31]{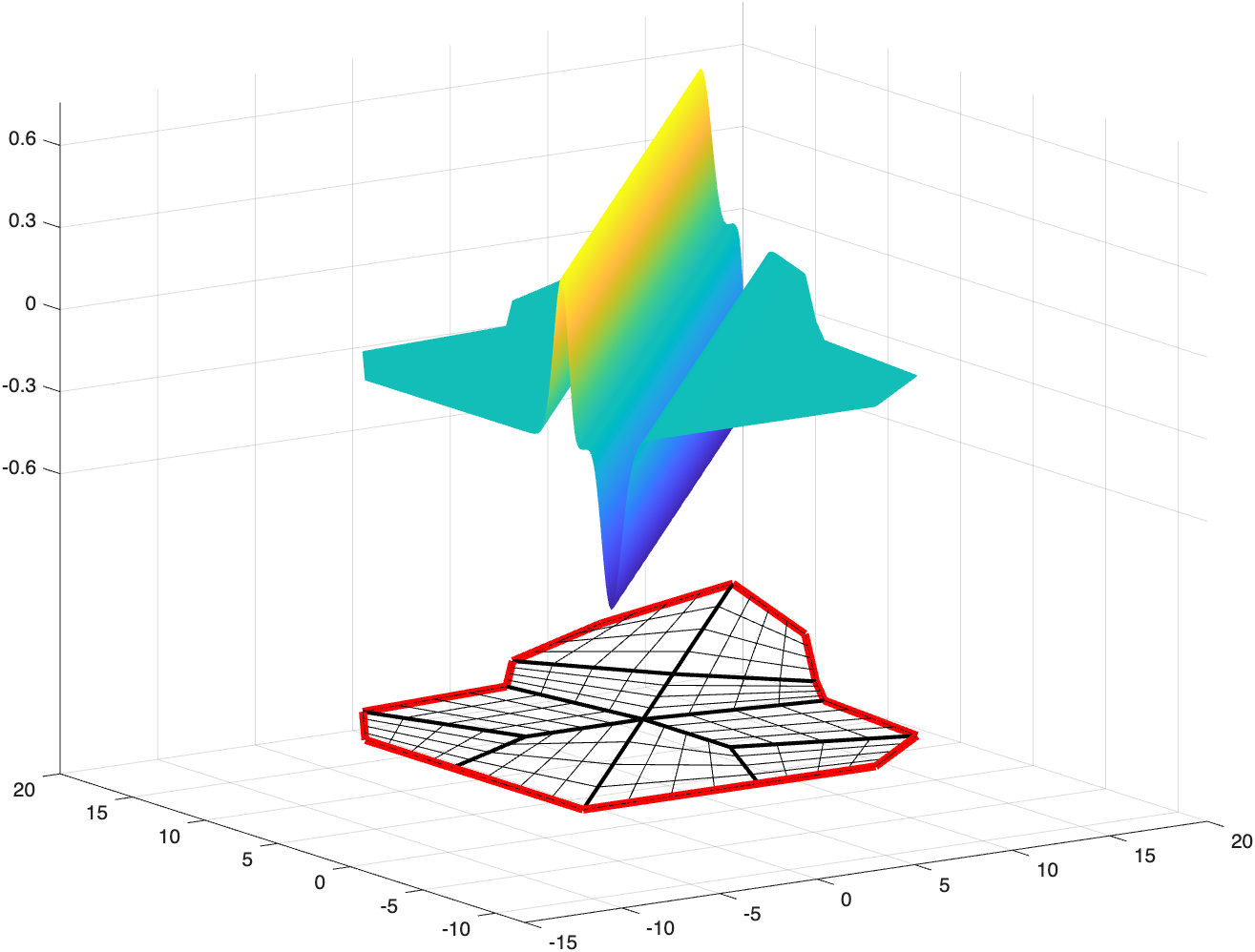} \\[0.2cm]
   {\new Example~\ref{ex:Example1}}   &  {\new Example~\ref{ex:Example2}} \\[0.8cm]
  \includegraphics[scale=0.353]{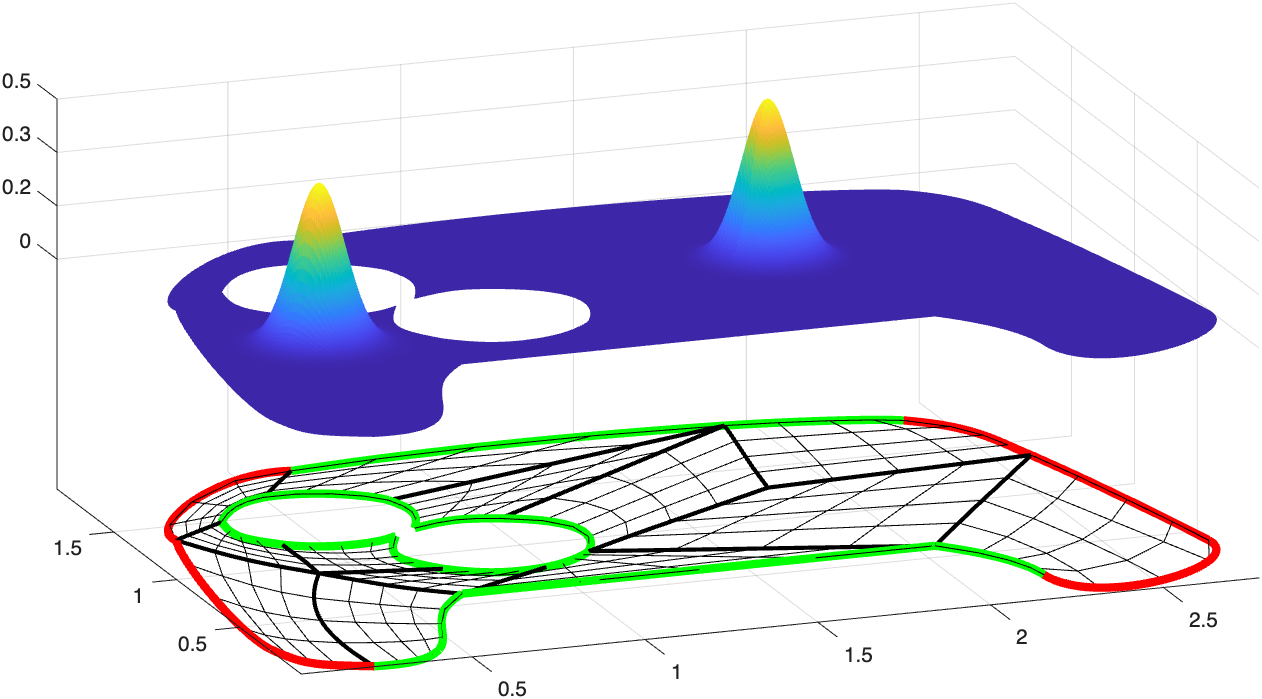}  \hskip1em &
  \includegraphics[scale=0.3]{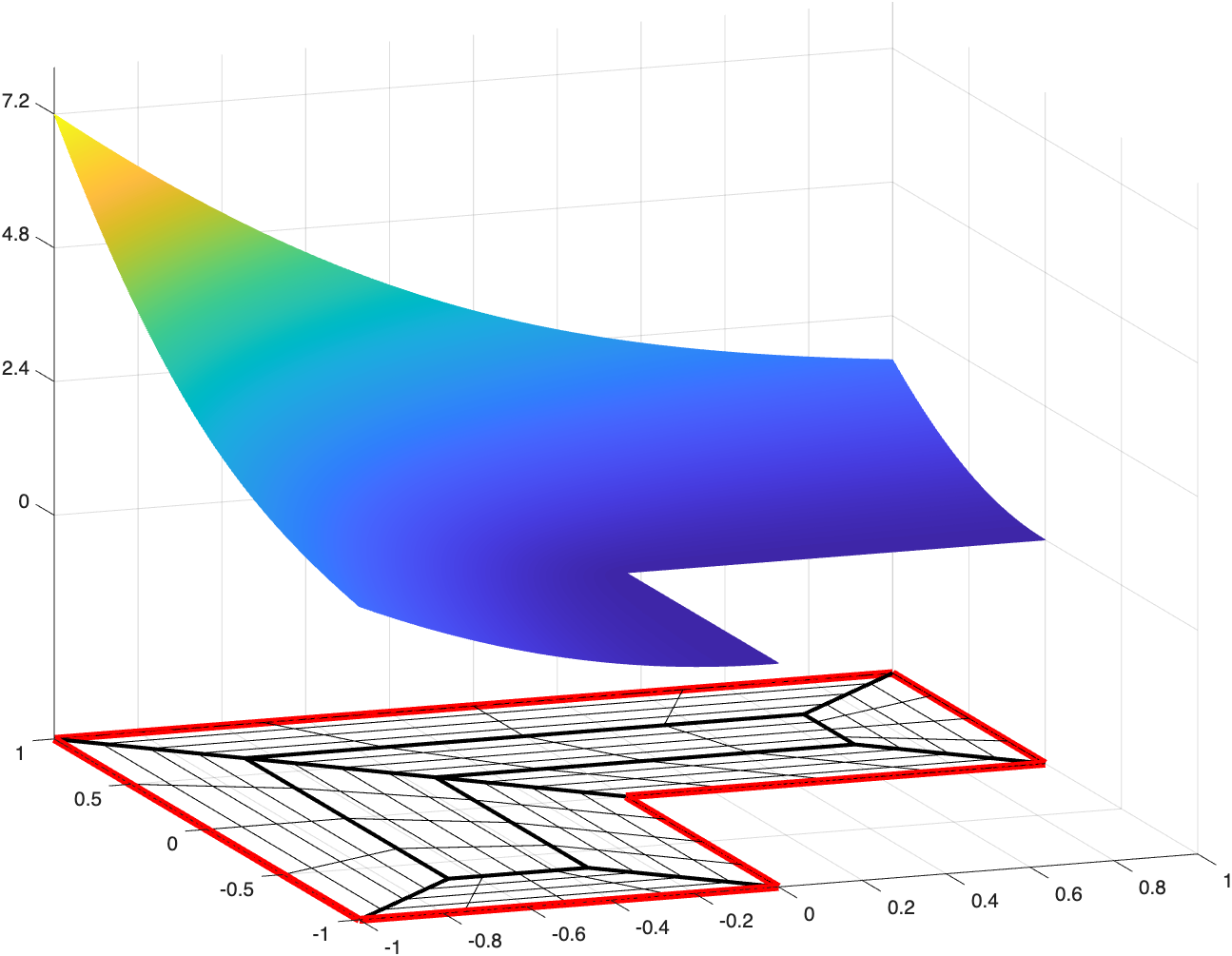}\\[0.4cm]
    {\new Example~\ref{ex:Example3}}  &   {\new Example~\ref{ex:Example4}} \\[0.4cm]
  \end{tabular}
\caption{{\new The numerical solutions {corresponding to} 
Example~$2$ (top left), Example~$3$ (top right), Example~$4$ (bottom left) and Example~$5$ (bottom right). The 
{solutions for Example~\ref{ex:Example1} and \ref{ex:Example3} }
are scaled by factors $3000$ and $100$, respectively.  
}
}
\label{fig:exactSolutions}
\end{figure}

\begin{ex} \label{ex:Example1}
We consider as computational domain the bilinear {\new thirteen-patch} domain, visualized in Fig.~\ref{fig:domains} (top left) possessing vertices with different patch valencies as well as one {\new triangular  and one rectangular} hole. 
{\new We consider the Kirchhoff--Love plate problem subject to the constant distributed load $f(\ab{x})=1$, which may represent, for example, a uniform pressure or a uniformly distributed transverse weight. We set $\Gamma_S=\emptyset$, while $\Gamma_C$ (red) and $\Gamma_F$ (blue) together constitute the entire boundary $\partial\Omega$, as indicated in Fig.~\ref{fig:exactSolutions} (top left). We choose the 
material and thickness parameters { (steel)} as
$$
E = 2 \cdot 10^{11} \;\mathrm{N/m^2}, \quad \nu_{\mathrm p}= 0.3 ,\quad t = 0.01\, \mathrm{m},
$$
and all boundary data 
acting on $\Gamma_C$ and $\Gamma_F$ are homogeneous, i.e.,
$g_0=g_1=g_2=g_3=0$.
}
{\new We start from a coarse mesh with initial mesh size $h_0=1/5$. Since no exact solution is prescribed in this example, the convergence is assessed by means of the  error estimate
$e(u_h,u_{h/2})$.}
We perform nine adaptive refinement steps using both hierarchical and truncated hierarchical B-splines and compare the resulting convergence behavior with that obtained from four levels of uniform refinement; see Fig.~\ref{fig:example1}. 

{\new We select admissible class~$\mu=2$.} In case of adaptive refinement we obtain optimal convergence rates of order~$(p-1)/2$ for both the hierarchical and truncated hierarchical construction and further for all considered degrees~$p \in \{3,4,5 \}$, while in case of uniform refinement the rates are clearly reduced. For 
degree $p=4$, we visualize in Fig.~\ref{fig:example1} (top left) the computed condition numbers for the linear system~\eqref{eq:spdFr} using standard diagonally scaling and the proposed Dirichlet preconditioner~\eqref{eq:Dirichlet_prec} 
and in Fig.~\ref{fig:example1_meshes} the resulting admissible meshes after five levels of adaptive refinement.

\begin{figure}[htb!]
    \centering
    \begin{tabular}{cc}
  \includegraphics[scale=0.36]{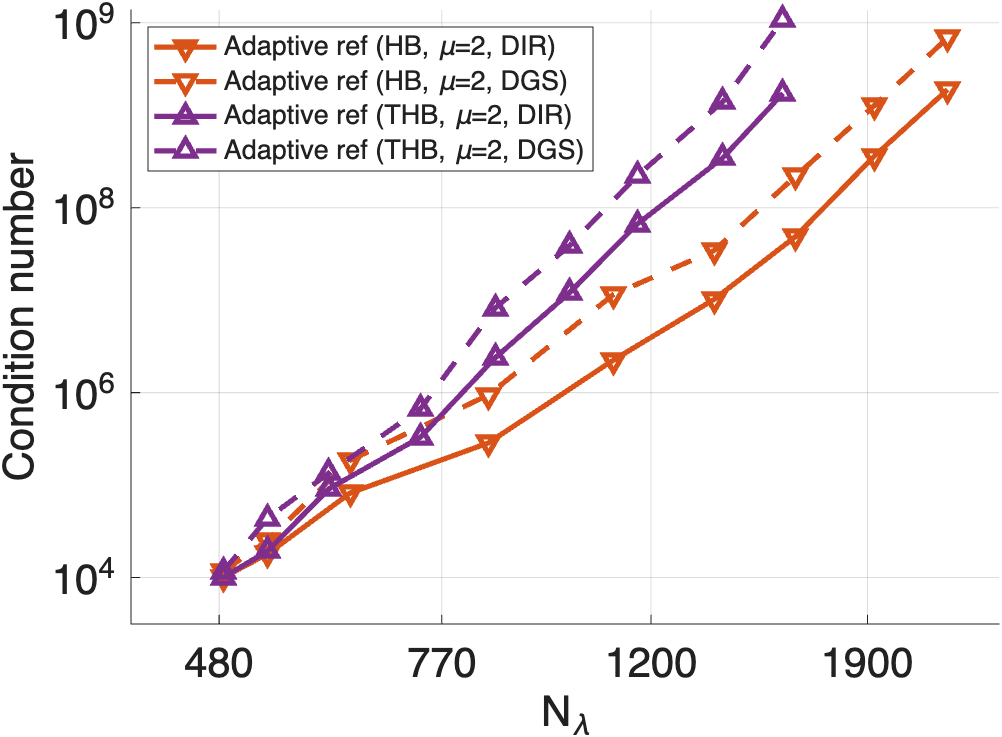}   & \hskip-1.em
  \includegraphics[scale=0.35]{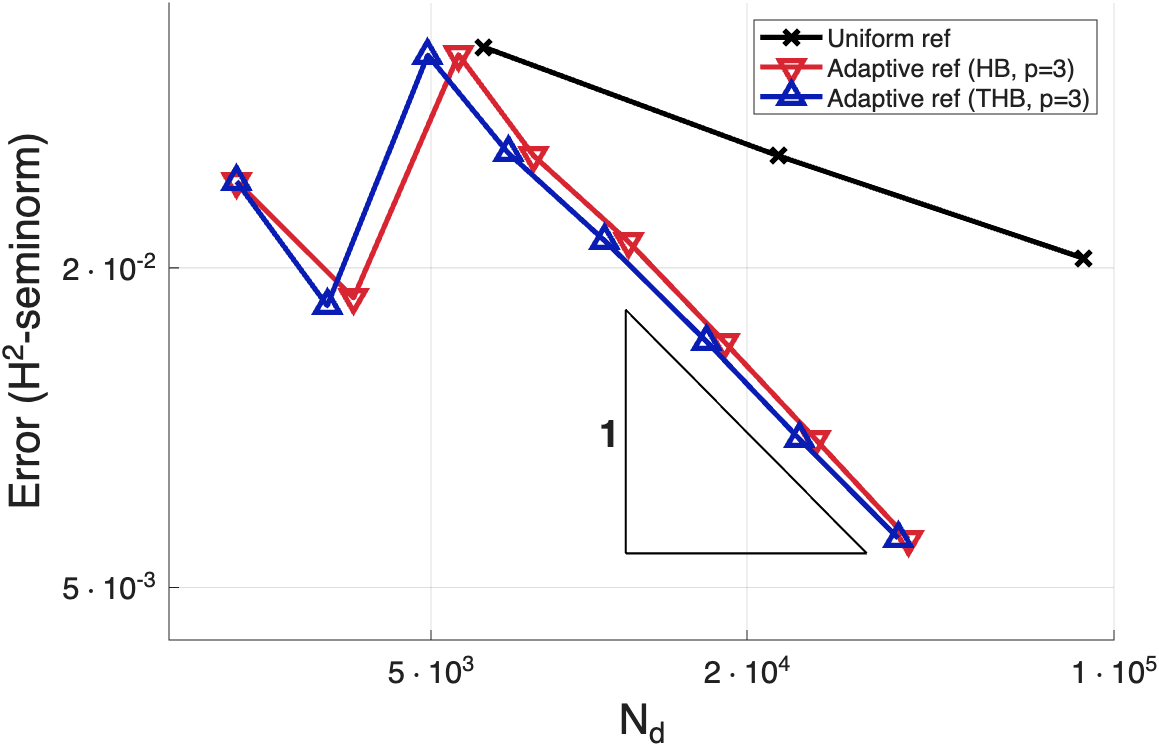} \\[0.2cm]
   Condition numbers $\kappa$ for $p=4$  & $p=3$ \\[0.8cm]
  \includegraphics[scale=0.35]{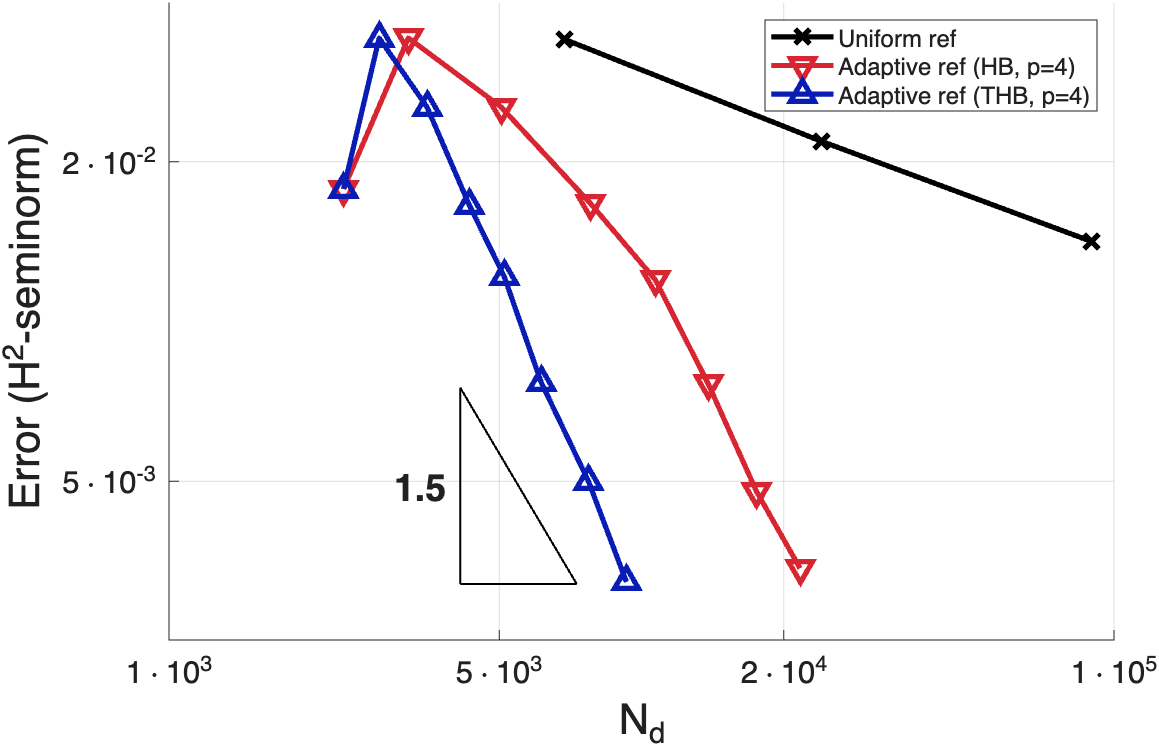}   & \hskip-1em
  \includegraphics[scale=0.35]{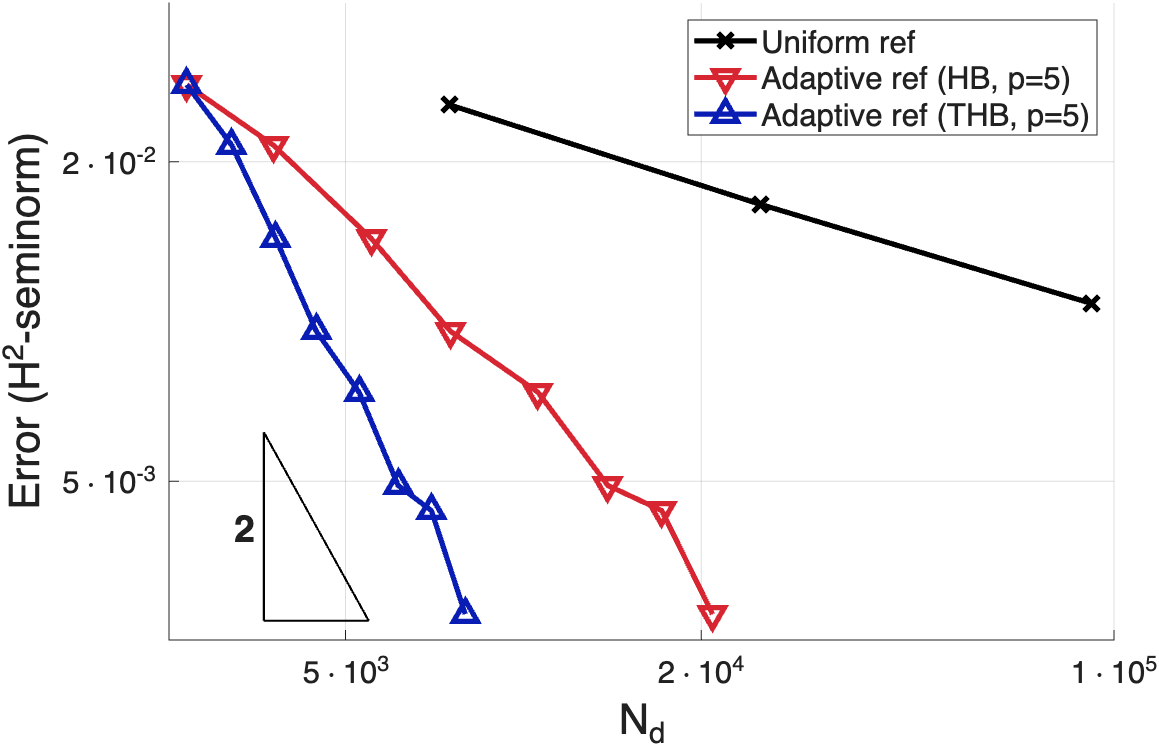}\\[0.4cm]
   $p=4$  &  $p=5$ \\[0.4cm]
  \end{tabular}
\caption{Example~\ref{ex:Example1}. {\new Condition numbers $\kappa$ (top left) for the linear system~\eqref{eq:spdFr} for $p=4$ versus the number of Lagrange multipliers $N_{\ab{\lambda}}$, 
where the dashed lines 
correspond to diagonal scaling (DGS) and the solid lines to the Dirichlet preconditioner (DIR). The remaining panels show the $h$--$h/2$ convergence plots for $p=3$ (top right), $p=4$ (bottom left), and $p=5$ (bottom right).}}
\label{fig:example1}
\end{figure}

\begin{figure}[htb!]
    \centering
    \begin{tabular}{cc}
  \includegraphics[scale=0.095]{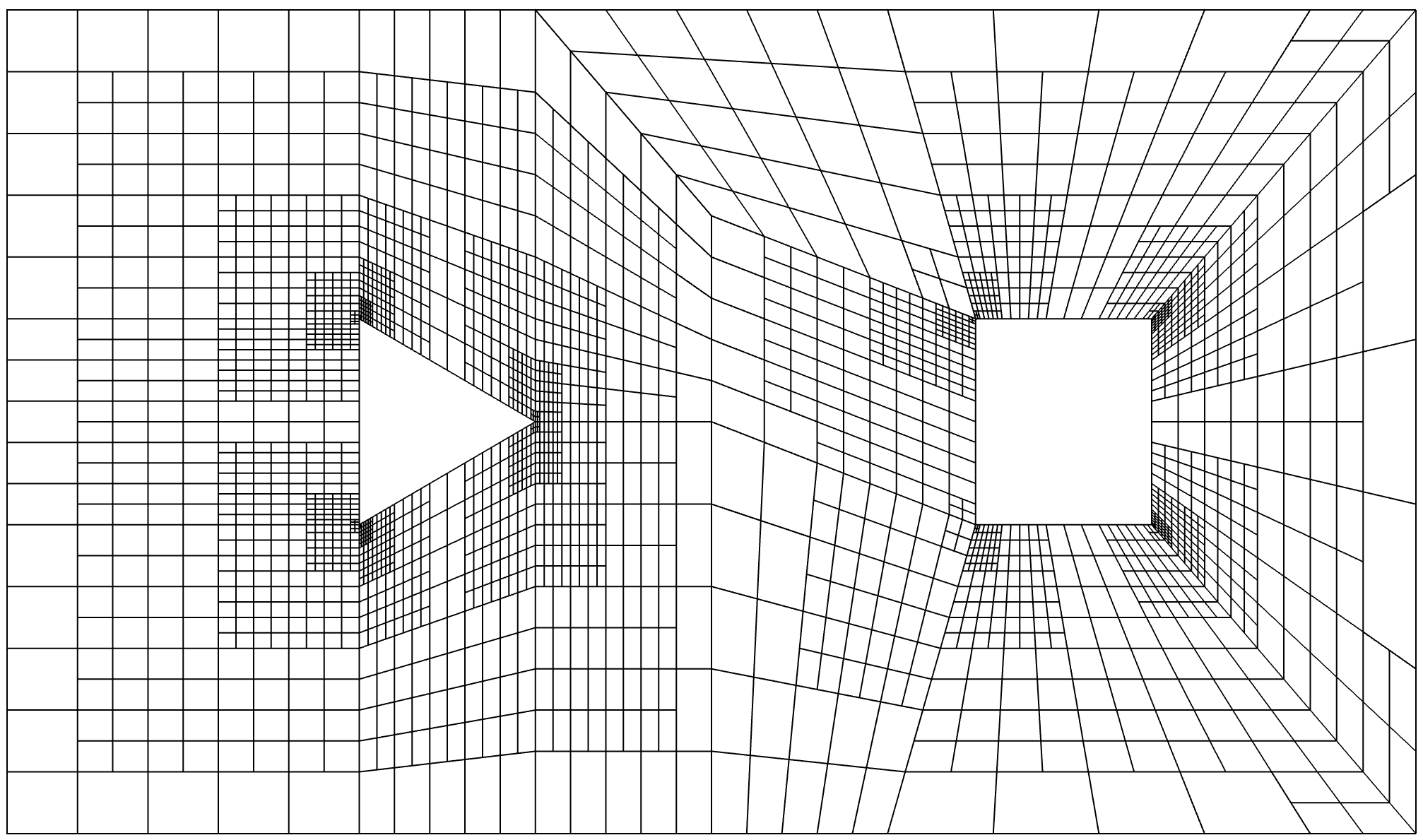}   & 
  \hspace{1cm}
  \includegraphics[scale=0.095]{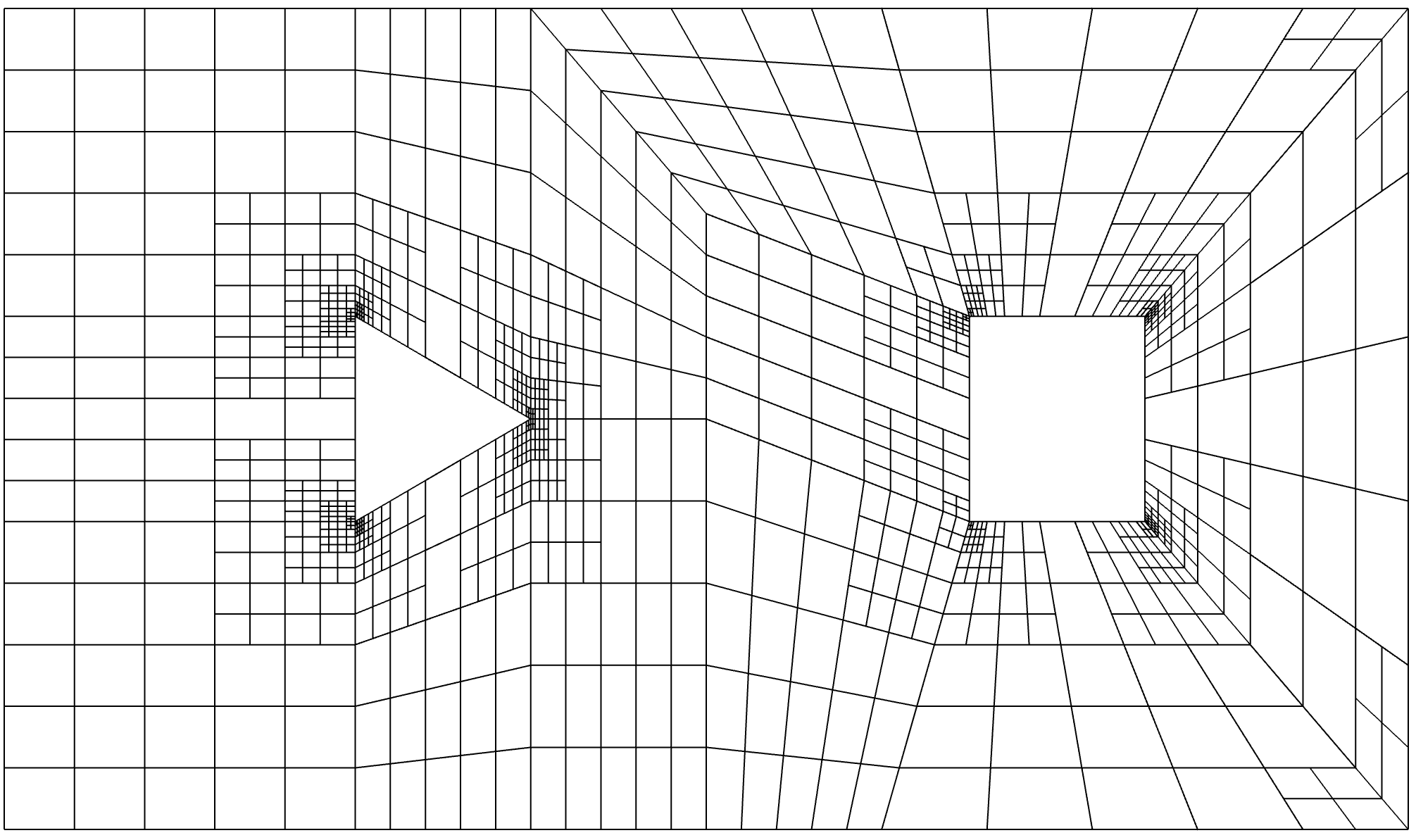} \\[0.2cm]
  \end{tabular}
\caption{Example~\ref{ex:Example1}. The admissible meshes for $\mu=2$ after five levels of adaptive refinement when using hierarchical B-splines for $p=4$ (left), and truncated hierarchical B-splines for $p=4$ (right).}
\label{fig:example1_meshes}
\end{figure}
\end{ex}

\begin{ex} \label{ex:Example2}
{\new The aim of this example is to study a line singularity for the ``standard'' biharmonic equation with clamped boundary conditions.
For this reason, we select 
$$
D_{\mathrm p}=1, \quad \nu_{\mathrm p}=0, \quad \Gamma_C=\partial \Omega,\quad \Gamma_S=\Gamma_F=\emptyset.
$$
}
We take the bilinear ten-patch domain shown in Fig.~\ref{fig:domains} (top right), and the exact solution
\begin{equation} \label{eq:exact_sol_example2}
u(x,y) = (y-x)^{\frac{13}{3}} e^{-(y-x)^2},
\end{equation}
see Fig.~\ref{fig:exactSolutions} (top right). {\new This solution has a singularity along the straight line $y=x$, which coincides with the patch interfaces in the central part of the domain. In addition, two smooth peaks occur near the singularity line and also require local refinement.} We again start with a coarse mesh with an initial mesh size $h_0=1/5$, and perform ten steps with our adaptive IETI-DP method. Fig.~\ref{fig:example2} presents the resulting errors in the $H^2$-seminorm with respect to the number of degrees of freedom $N_{d}$ for the hierarchical as well as for the truncated hierarchical construction, and shows additionally for comparison the errors obtained by performing uniform refinement for five levels.

Following the observation for the Poisson equation in \cite{BuGaGiPrVa22, BrGiKaVa23} and generalizing it to the biharmonic 
{\new equation}, the expected convergence rate in the $H^2$-seminorm for isotropic meshes is $\min ( \frac{p-1}{2}, \frac{17}{12})$ for uniform refinement and $\min ( \frac{p-1}{2}, \frac{17}{6})$ for adaptive refinement. This is due to the fact that the exact solution~\eqref{eq:exact_sol_example2} belongs to the space $H^{29/6-\epsilon}(\Omega)$ for any $\epsilon>0$. As numerically demonstrated in Fig.~\ref{fig:example2}, the obtained orders for the local refinement are indeed again advantageous in comparison to the uniform refinement. 

For 
admissible class~${\new \mu=3}$ and 
degree $p=4$, Fig.~\ref{fig:example2} (top left) shows further the computed condition numbers for the linear system~\eqref{eq:spdFr} using standard diagonal scaling and the proposed Dirichlet preconditioner~\eqref{eq:Dirichlet_prec}, 
and Fig.~\ref{fig:example2_meshes} presents the resulting admissible meshes after nine levels of adaptive refinement.

\begin{figure}[htb!]
    \centering
    \begin{tabular}{cc}
  \includegraphics[scale=0.19]{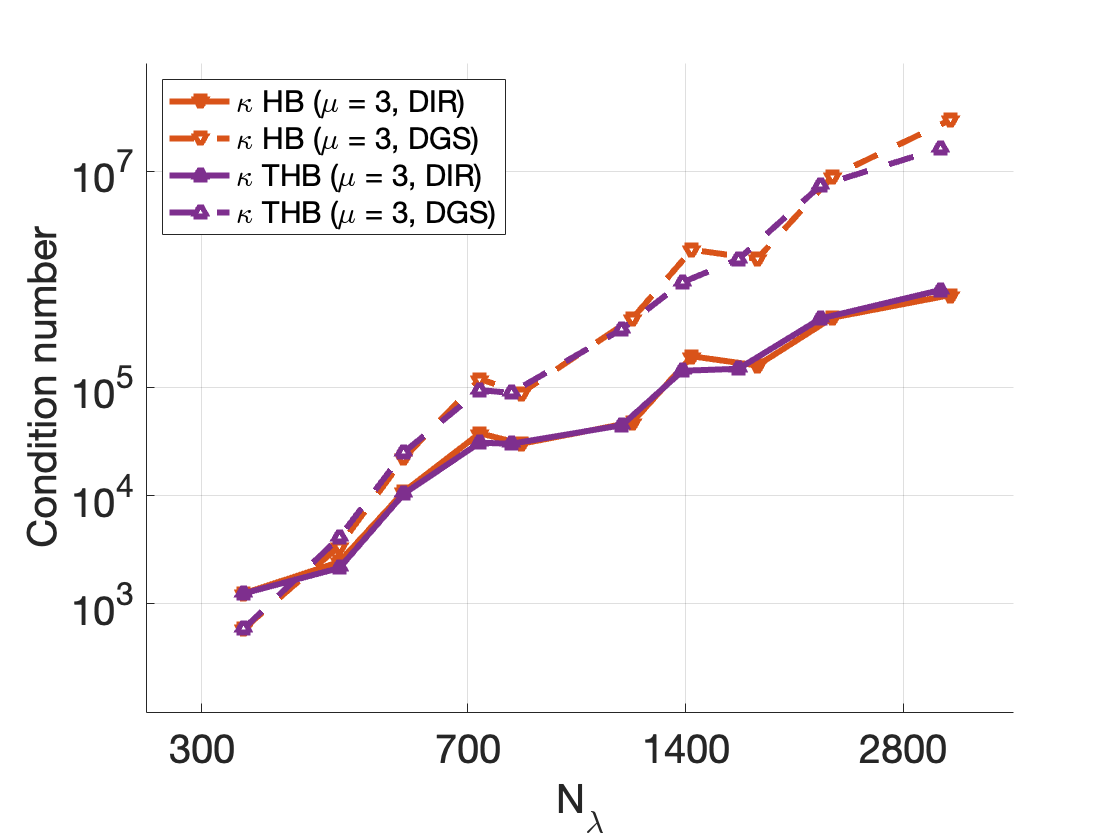}   & \hskip-1em
  \includegraphics[scale=0.19]{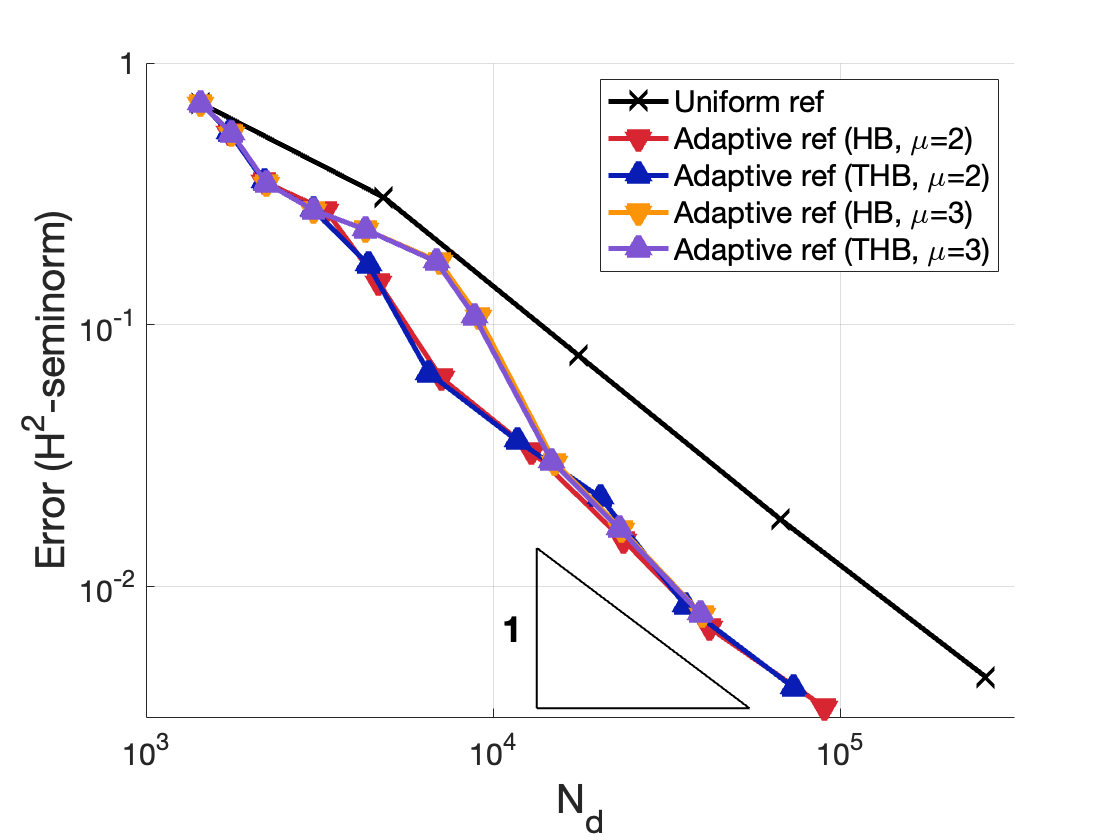} \\[0.2cm]
   Condition numbers $\kappa$ for $p=4$  & $p=3$ \\[0.8cm]
  \includegraphics[scale=0.19]{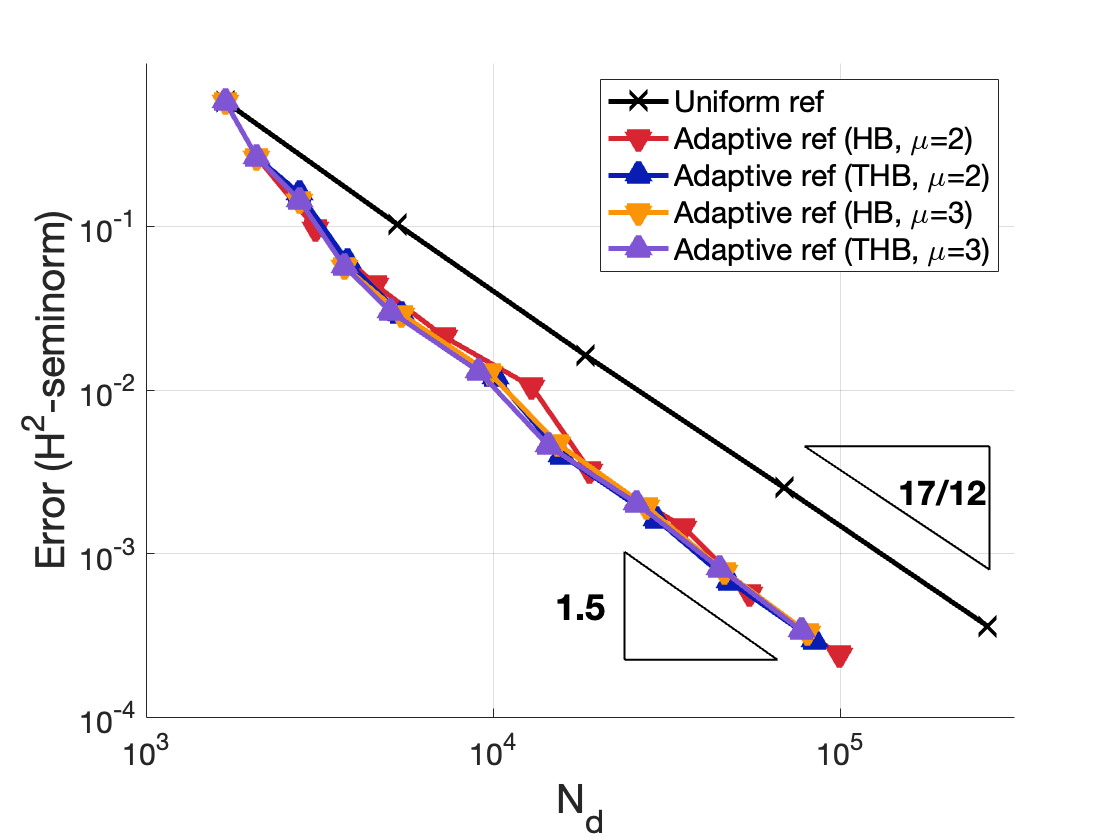}   & \hskip-1.em
  \includegraphics[scale=0.19]{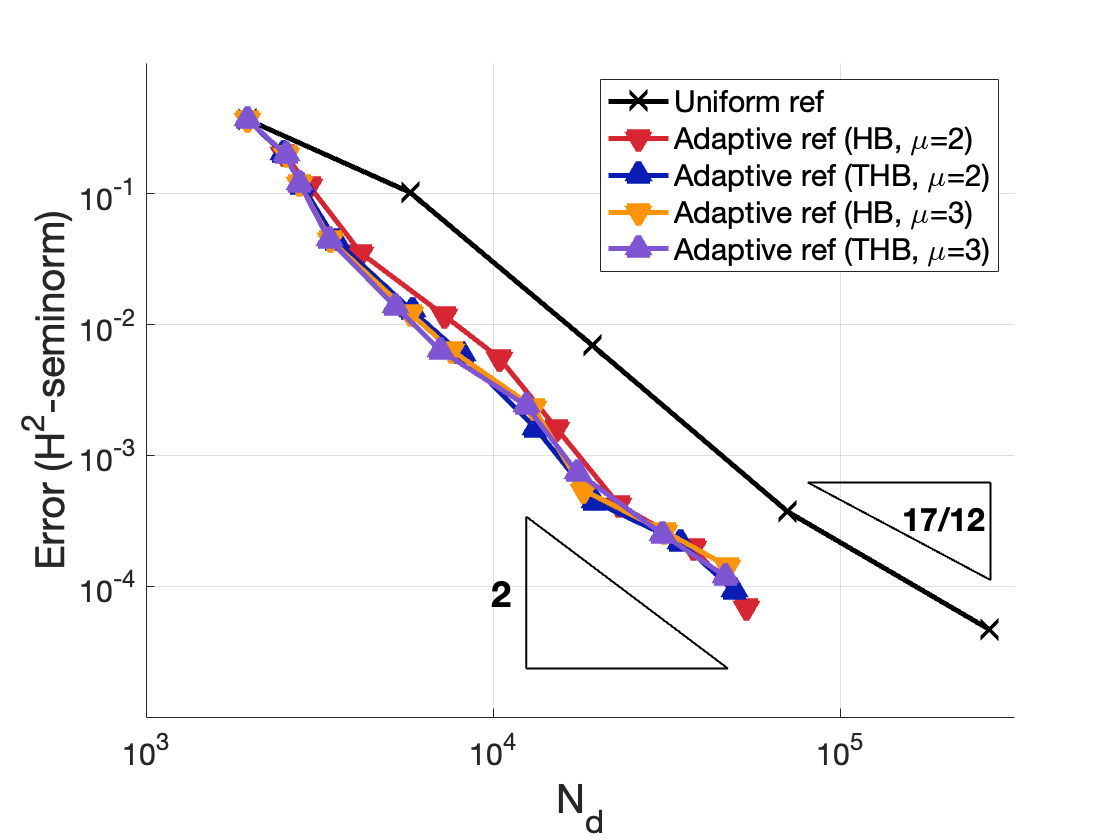}\\[0.4cm]
   $p=4$  &  $p=5$ \\[0.4cm]
  \end{tabular}
\caption{Example~\ref{ex:Example2}. {\new Condition numbers $\kappa$ (top left) for the linear system~\eqref{eq:spdFr} for $p=4$ versus 
the number of Lagrange multipliers $N_{\ab{\lambda}}$, 
{where the} dashed lines 
{correspond} to diagonal scaling (DGS) and the solid lines to the Dirichlet preconditioner (DIR). {The remaining panels show the} convergence plots for $p=3$ (top right), $p=4$ (bottom left), and $p=5$ (bottom right).}}
    \label{fig:example2}
\end{figure}

\begin{figure}[htb!]
    \centering
    \begin{tabular}{cc}
  \includegraphics[scale=0.28]{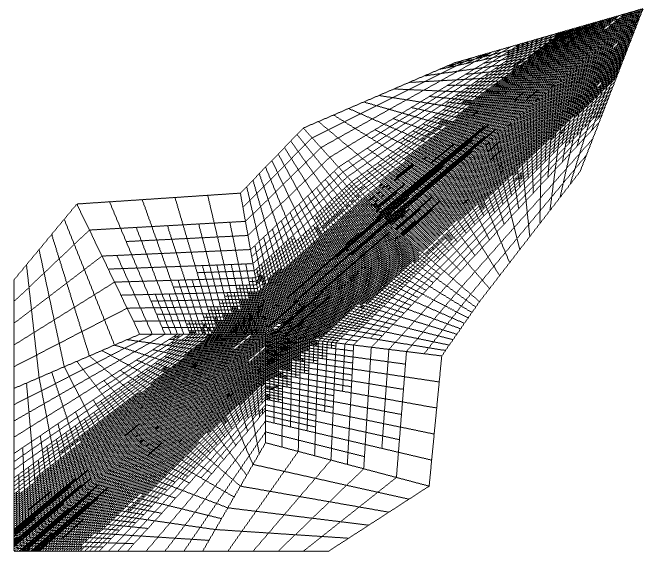}   & \hskip3em
  \includegraphics[scale=0.28]{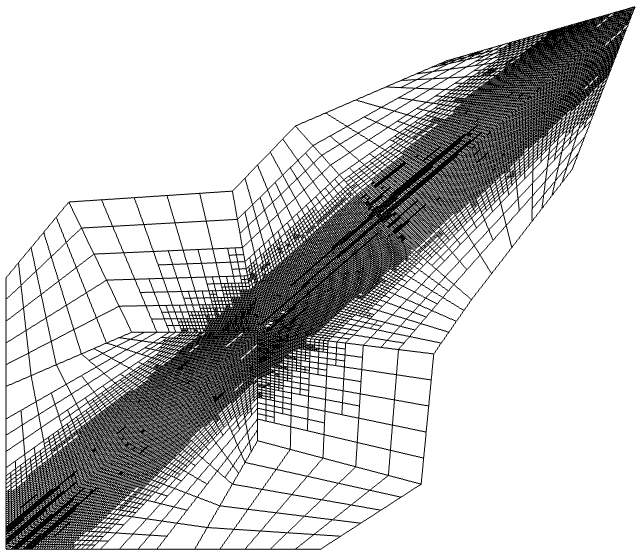} \\[0.2cm]
  \end{tabular}
\caption{Example~\ref{ex:Example2}. The admissible meshes for $\mu=3$ after nine levels of adaptive refinement when using hierarchical B-splines for $p=4$ (left), and truncated hierarchical B-splines for $p=4$ (right).}
    \label{fig:example2_meshes}
\end{figure}
\end{ex}

\begin{ex} \label{ex:Example3}
We consider the non-bilinear thirteen-patch controller domain presented in Fig.~\ref{fig:domains} (bottom left), parametrized by the analysis-suitable $G^1$ multi-patch geometry~$\ab{G}$ constructed in \cite[Example~2]{FaKaKoVi24}, given by individual geometry mappings~$\ab{G}^{(i)} \in \mathcal{S}_{1/3}^{\ab{3}, \ab{1}}([0,1]^2) \times \mathcal{S}_{1/3}^{\ab{3}, \ab{1}}([0,1]^2)$, $i \in \mathcal{I}_{\Omega}$. 
{\new We study the Kirchhoff--Love plate problem with a load given by 
\begin{equation} \label{eq:load_example3}
f(\ab{x}) = f_1(\ab{x}) + f_2(\ab{x}), \quad f_i(\ab{x}) = 16 A^2 e^{-A\, \|\ab{x} - \ab{\Xi}_i \|^2} \left( (A \|\ab{x} - \ab{\Xi}_i \|^2 - 2)^2-2 \right),\;\; i=1,2,  
\end{equation}
where $A=100$, and 
$
\ab{\Xi}_1 
$ and 
$
\ab{\Xi}_2 
$
are the two inner vertices of the multi-patch domain, cf. Fig.~\ref{fig:exactSolutions} (bottom left). 
{\new Moreover, we set $\Gamma_F=\emptyset$, while $\Gamma_C$ (red) and $\Gamma_S$ (green) together constitute the entire boundary $\partial\Omega$, as indicated in Fig.~\ref{fig:exactSolutions} (bottom left). We choose now the 
material and thickness parameters { (polycarbonate)} as}
$$
E = 2.3 \cdot 10^{9} N/m^2 , \quad \nu_{\mathrm p}= 0.37 ,\quad t = 0.01\, m,
$$
and impose homogeneous boundary functions
$g_0 = g_1 = g_2 = 0$. We start from a coarse mesh with initial mesh size $h_0=1/6$ and fix the degree to $p=3$. Since no exact solution is prescribed, we assess convergence using the $h$--$h/2$ relative difference
$e(u_h,u_{h/2})$. We perform nine adaptive refinement steps with both hierarchical and truncated hierarchical B-splines and compare the resulting behavior with uniform refinement, {see Fig.~\ref{fig:example3}~(bottom left)}.
}
{\new We select now the admissible class~$\mu=3$.} 
The computed errors decay with an optimal order of $\mathcal{O}(N_{d}^{-1})$, independent of an adaptive or uniform refinement strategy. However, the obtained errors for the adaptive cases, i.e.~for using hierarchical and truncated hierarchical B-splines, are significantly smaller than for the uniform case. Additionally to the convergence plots, Fig.~\ref{fig:example3} presents 
the computed condition numbers for the linear system~\eqref{eq:spdFr} using standard diagonally scaling and the proposed Dirichlet preconditioner~\eqref{eq:Dirichlet_prec} (bottom right), and the resulting admissible meshes after eight levels of adaptive refinement (top row). 

\begin{figure}[htb!]
    \centering
    \begin{tabular}{cc}
  \includegraphics[scale=0.115]{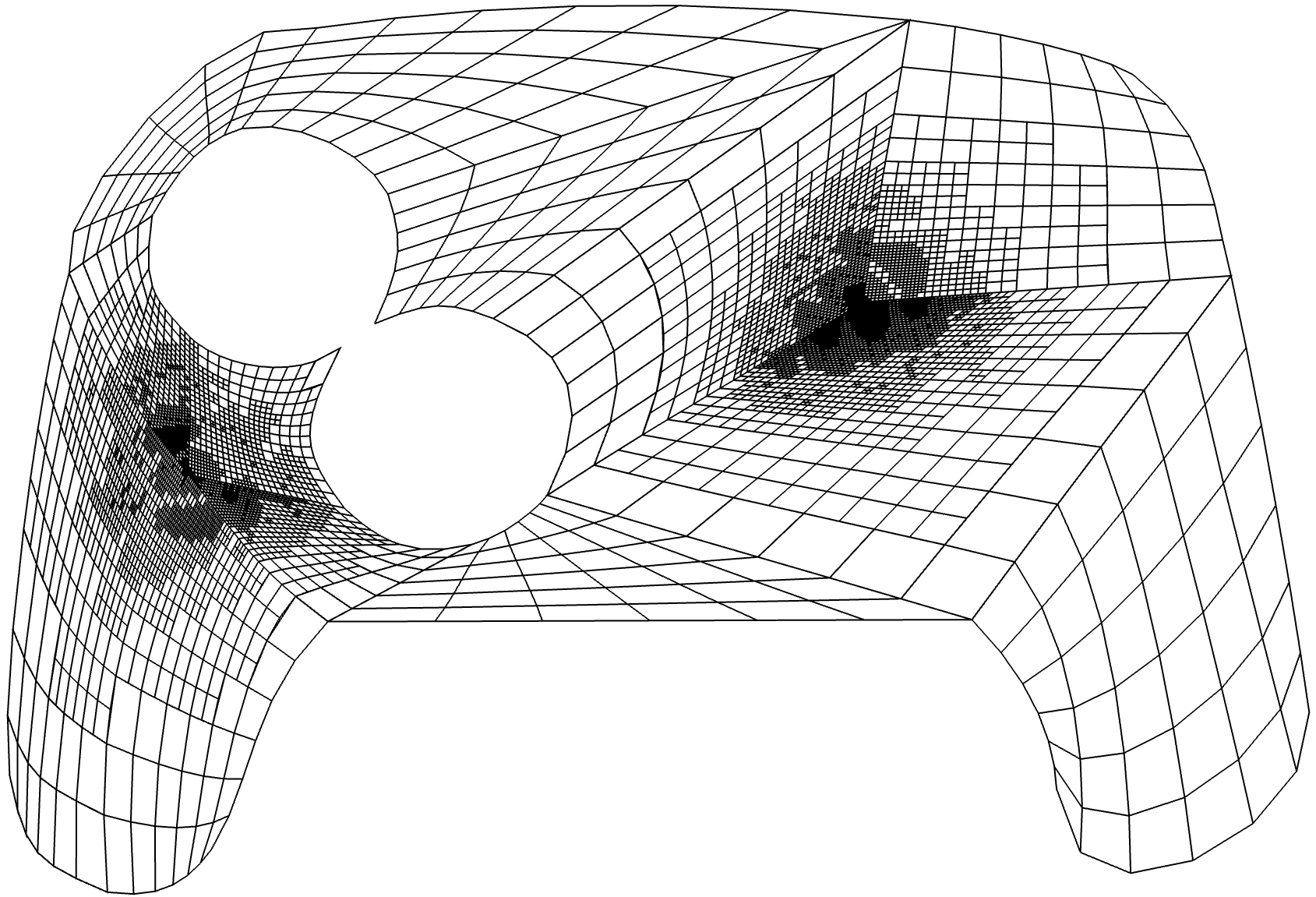}   &
  \hskip0.5em
  \includegraphics[scale=0.115]{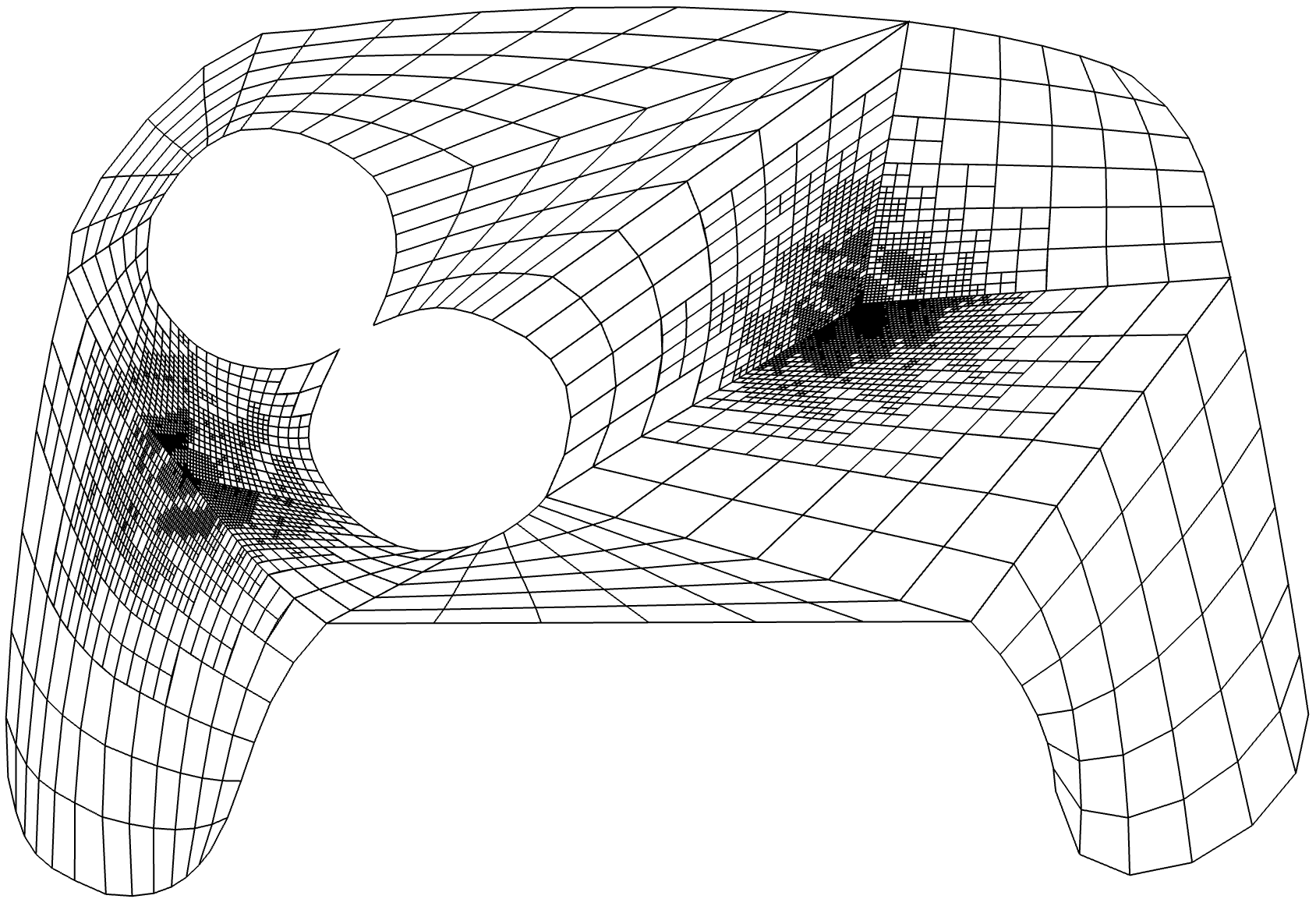} \\[0.2cm]
  \includegraphics[scale=0.36]{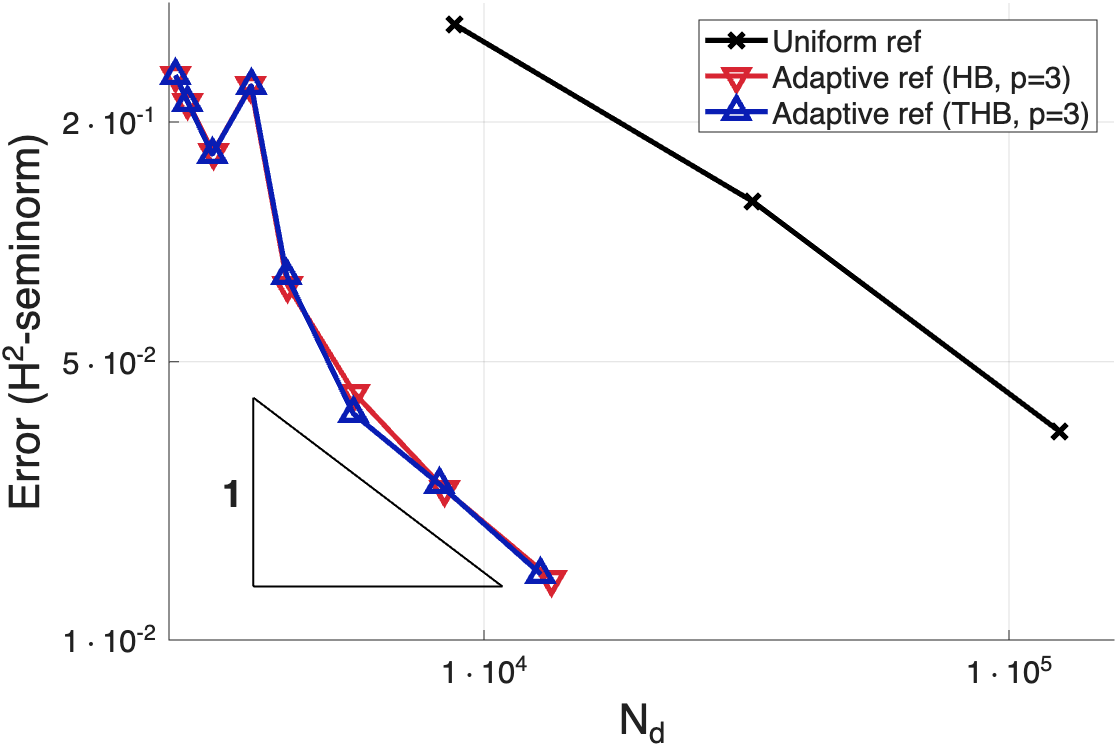}   & \hspace{0.2cm}
  \includegraphics[scale=0.36]{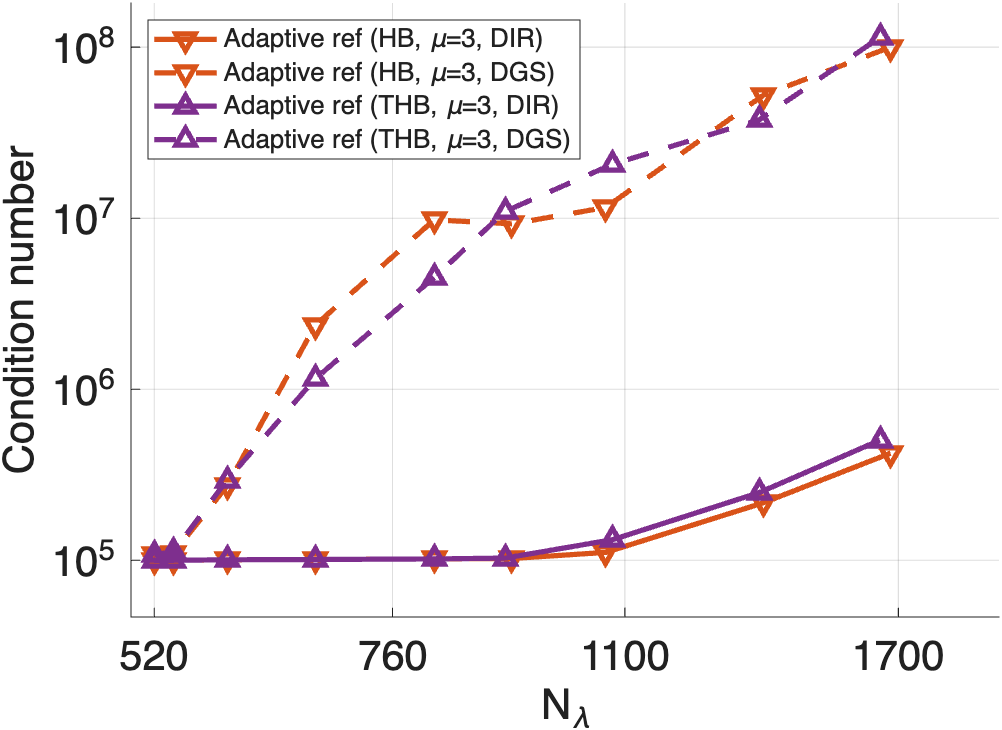} \\[0.2cm]
  $p=3$ & Condition numbers $\kappa$ for $p=3$ \\[0.8cm]
  \end{tabular}
\caption{Example~\ref{ex:Example3}. {\new The admissible meshes for $\mu=3$ after eight levels of adaptive refinement using hierarchical B-splines (top left) and truncated hierarchical B-splines (top right), the $h$--$h/2$ convergence plots for the different refinement strategies (bottom left), and the condition numbers $\kappa$ (bottom right) for the linear system~\eqref{eq:spdFr} versus the number of Lagrange multipliers $N_{\ab{\lambda}}$. The dashed lines correspond to diagonal scaling (DGS), and the solid lines to the Dirichlet preconditioner (DIR).}}
    \label{fig:example3}
\end{figure}
\end{ex}

\begin{ex} \label{ex:Example4}
{\new In the final example, we consider the bilinear eight-patch $L$-shaped domain shown in Fig.~\ref{fig:domains} (bottom right) and prescribe the exact solution~$u$ as}
\begin{equation} \label{eq:exact_sol_example4}
u(\rho, \theta) = \rho^{z+1} (C_1 F_1(\theta) -C_2 F_2 (\theta) )
\end{equation}
in polar coordinates~$(\rho, \theta)$ with 
\begin{equation*}
C_1  =  \frac{1}{z-1}\sin \left( \frac{3(z-1)\pi}{2} \right) - \frac{1}{z+1} \sin \left(\frac{3(z+1)\pi}{2} \right),
\end{equation*}
\begin{equation*}
C_2  =  \cos \left( \frac{3(z-1)\pi}{2} \right) - \cos \left(\frac{3(z+1)\pi}{2} \right),
\end{equation*}
\begin{equation*}
F_1(\theta)  =  \cos ((z-1) \theta) - \cos((z+1)\theta), 
\end{equation*}
\begin{equation*}
F_2(\theta)  =  \frac{1}{z-1}\sin ((z-1) \theta) - \frac{1}{z+1}\sin((z+1)\theta),
\end{equation*}
and $z=0.544483736782464$. {\new The exact solution~\eqref{eq:exact_sol_example4} has already been studied, e.g., in \cite[Section~3.4]{Gr92} and \cite[Example~$3$]{BrGiKaVa23}, and has a singularity at $(0,0)$, which coincides with the re-entrant corner of the $L$-shaped multi-patch domain in Fig.~\ref{fig:domains} (bottom right).}
{\new Thereby, the multi-patch domain 
is, as the exact solution~$u$, 
taken from \cite[Example~$3$]{BrGiKaVa23} with the aim of comparing the present results with those reported there.} For this purpose, we 
{\new
take 
$$
D_{\mathrm p}=1, \quad \nu_{\mathrm p}=0, \quad \Gamma_C=\partial \Omega,\quad \Gamma_S=\Gamma_F=\emptyset,
$$
start from a coarse mesh with initial mesh size $h_0=1/5$, and perform ten adaptive refinement steps with the proposed IETI-DP method. For comparison, we also perform uniform refinement.} Fig.~\ref{fig:example4} illustrates the errors in the $H^2$-seminorm with respect to the number of degrees of freedom $N_{d}$ obtained from our adaptive IETI-DP method using hierarchical and truncated hierarchical B-splines as well as from performing uniform refinement. 
{\new For all degrees $p \in \{3,4,5\}$, the adaptive computations attain the optimal convergence order $(p-1)/2$ for both hierarchical and truncated hierarchical B-splines, whereas uniform refinement yields clearly reduced rates. Furthermore, the errors obtained with adaptive refinement are of the same order of magnitude as the errors obtained with the method of~\cite{BrGiKaVa23} and reported in \cite[Fig.~9]{BrGiKaVa23}. This agreement is expected because the underlying hierarchical $C^1$-smooth spline spaces and meshes of the two approaches differ only slightly in the vicinity of an inner vertex of the multi-patch domain, cf.~Remark~\ref{rem:comparison} and Section~\ref{subsec:adaptive_ref}.} 

\begin{figure}[htb!]
    \centering
    \begin{tabular}{cc}
  \includegraphics[scale=0.19]{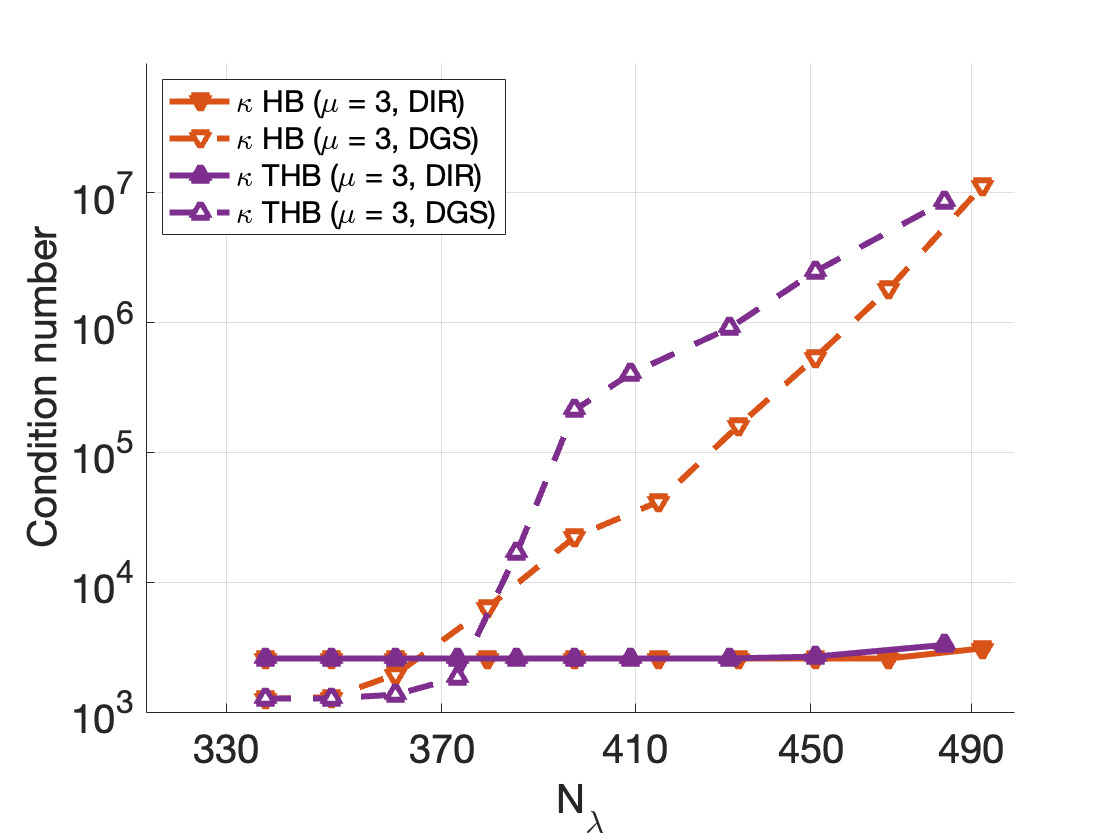}   & \hskip-1em
  \includegraphics[scale=0.19]{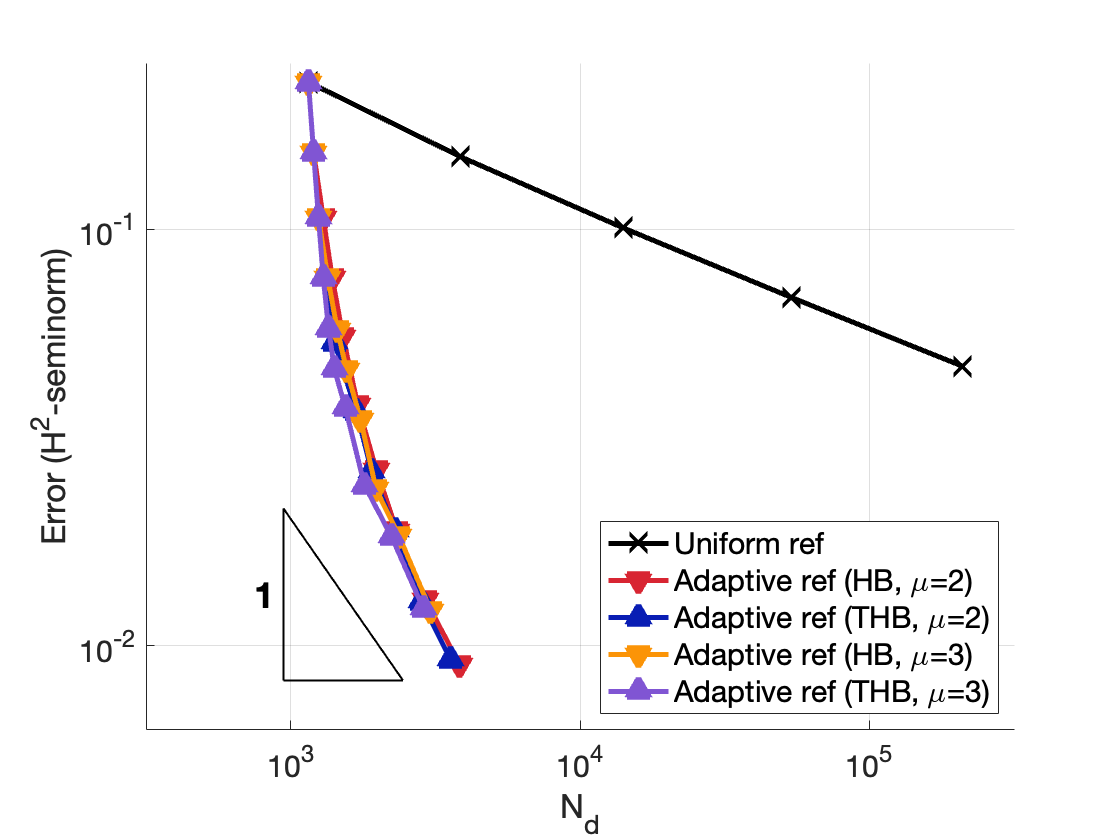} \\[0.2cm]
   Condition numbers $\kappa$ for $p=4$  & $p=3$ \\[0.8cm]
  \includegraphics[scale=0.19]{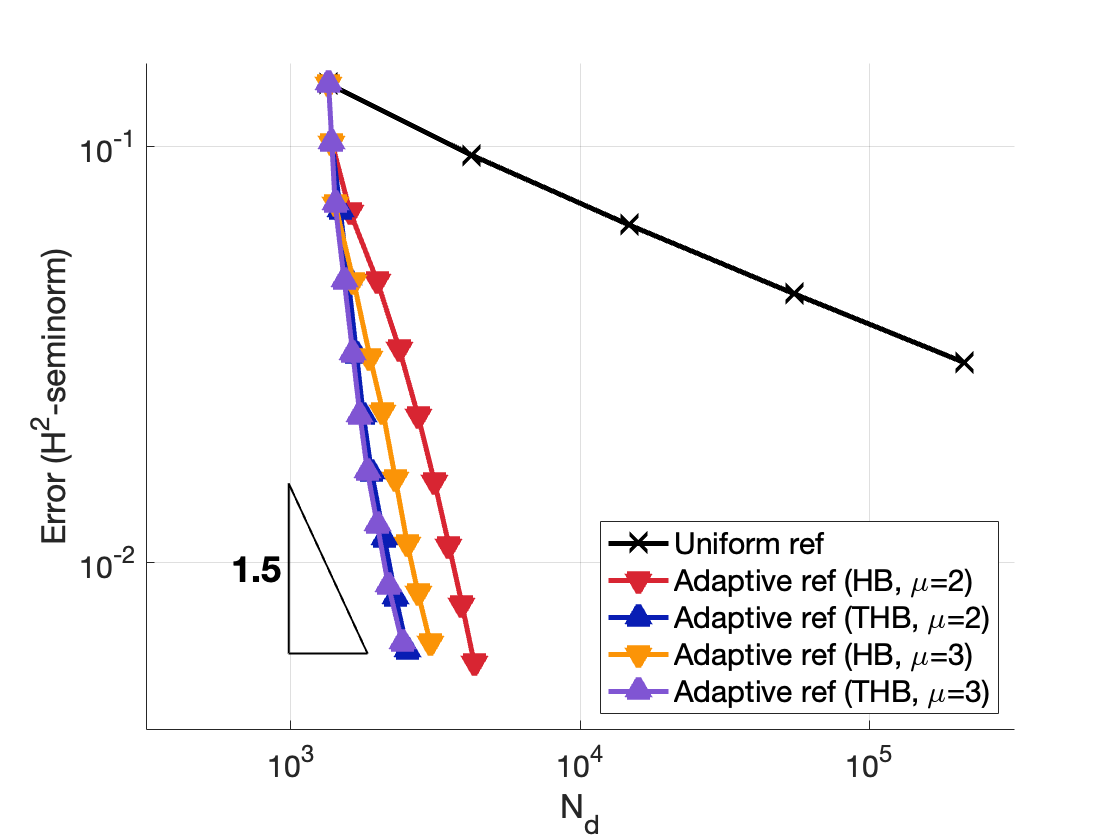}   & \hskip-1em
  \includegraphics[scale=0.19]{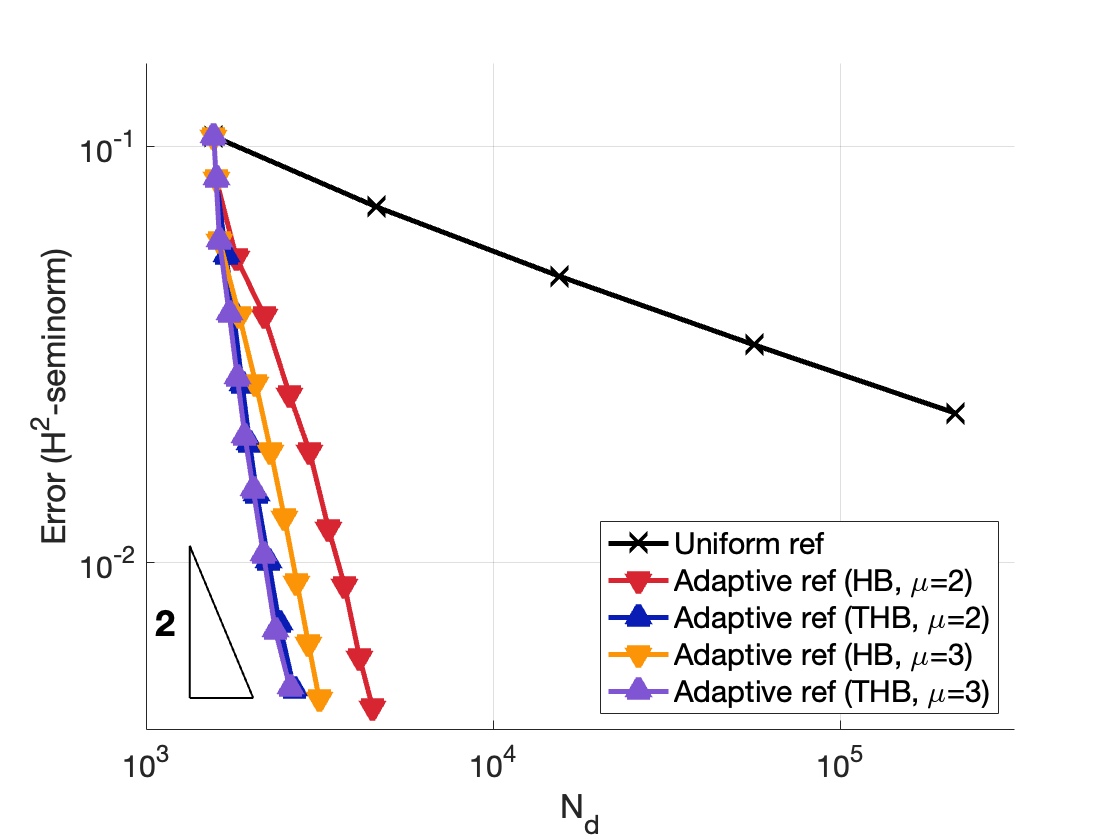}\\[0.4cm]
   $p=4$  &  $p=5$ \\[0.4cm]
  \end{tabular}
\caption{\new Example~\ref{ex:Example4}. Condition numbers $\kappa$ (top left) for $p=4$ 
{versus} the number of Lagrange multipliers $N_{\ab{\lambda}}$, 
{where the} dashed lines 
{correspond} to diagonal scaling (DGS) and the solid lines to the Dirichlet preconditioner (DIR). 
The remaining panels show the} convergence plots for $p=3$ (top right), $p=4$ (bottom left), and $p=5$ (bottom right).
    \label{fig:example4}
\end{figure}

\begin{figure}[htb!]
    \centering
    \begin{tabular}{cc}
  \includegraphics[scale=0.26]{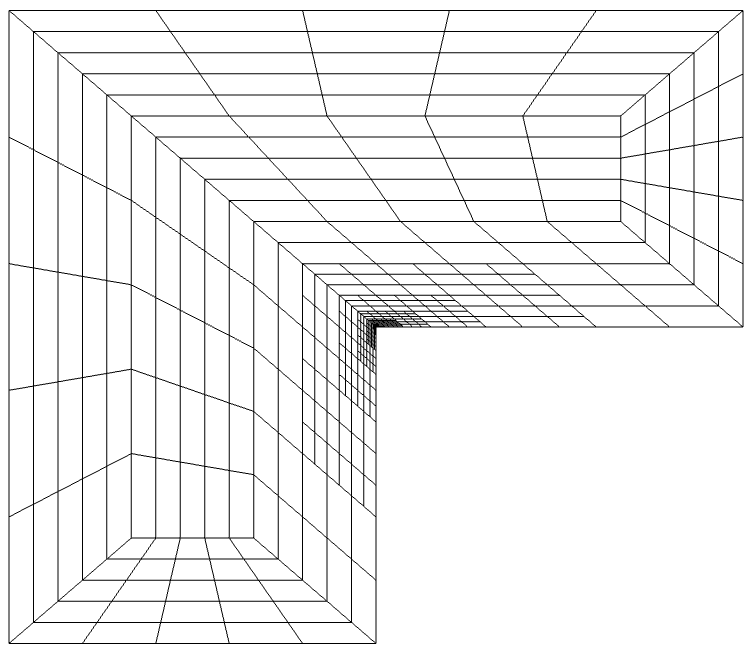}   &
   \hspace{0.9cm}
  \includegraphics[scale=0.26]{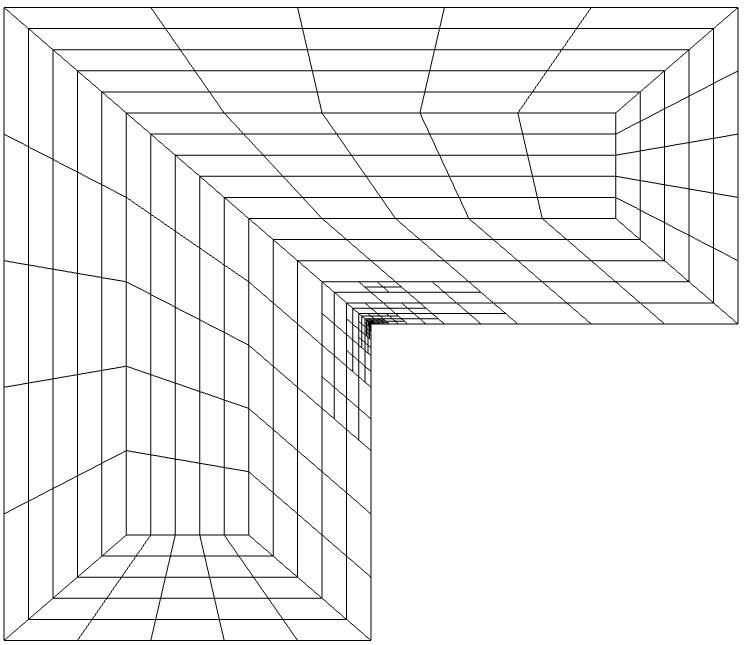} \\[0.2cm]
  \end{tabular}
\caption{Example~\ref{ex:Example4}. The admissible meshes for $\mu=3$ after six levels of adaptive refinement when using hierarchical B-splines for $p=4$ (left), and truncated hierarchical B-splines for $p=4$ (right).}
    \label{fig:example4_meshes}
\end{figure}

We consider in more detail the case $\mu=3$ and $p=4$. Figure~\ref{fig:example4} (top left) shows the condition numbers $\kappa$ of the linear system~\eqref{eq:spdFr} obtained with standard diagonal scaling and with the proposed Dirichlet preconditioner~\eqref{eq:Dirichlet_prec}, while Fig.~\ref{fig:example4_meshes} shows the admissible meshes after six adaptive refinement levels. Table~\ref{tab:ConditionAndStepsHierarchical} additionally reports the number of conjugate-gradient iterations required to solve the preconditioned system~\eqref{eq:preconditionedSystem} with the Dirichlet preconditioner $\ab{M}$. The iteration is terminated once the residual has been reduced by a factor of $10^{-6}$ relative to the initial residual. 

\begin{table}[htb!]
\centering
\begin{tabular}{|c||c|c|c|}
\hline
\multicolumn{4}{|c|}{Hierarchical refinement} \\
\hline
\hline
Level~$\ell$ & $p=3$ & $p=4$ & $p=5$ \\
 \hline        
0 & 34 & 45 & 59\\ 
1 & 36 & 44 & 59\\
2 & 35 & 44 & 62\\ 
3 & 36 & 46 & 63\\ 
4 & 35 & 46 & 66\\ 
5 & 35 & 45 & 69\\ 
6 & 38 & 47 & 72\\
7 & 36 & 47 & 76\\ 
8 & 42 & 49 & 81\\ 
9 & 49 & 53 & 89\\ 
\hline
\end{tabular}
\qquad
\begin{tabular}{|c||c|c|c|}
\hline
\multicolumn{4}{|c|}{Truncated hierarchical refinement} \\
\hline
\hline
Level~$\ell$ & $p=3$ & $p=4$ & $p=5$ \\
 \hline         
0 & 34 & 45 & 59\\ 
1 & 36 & 46 & 61\\
2 & 35 & 46 & 63\\ 
3 & 36 & 46 & 63\\ 
4 & 35 & 46 & 65\\ 
5 & 38 & 46 & 67\\ 
6 & 38 & 51 & 72\\
7 & 38 & 57 & 88\\ 
8 & 48 & 59 & 78\\ 
9 & 52 & 66 & 101\\ 
\hline
\end{tabular}
\caption{Example~\ref{ex:Example4}. Number of iterations required by the preconditioned conjugate-gradient method to solve~\eqref{eq:preconditionedSystem} for the admissibility class $\mu=3$, degrees $p=3,4,5$, and levels $\ell=0,1,\ldots,9$ of adaptive refinement, using hierarchical splines (left) and truncated hierarchical splines (right).}
\label{tab:ConditionAndStepsHierarchical}
\end{table}
\end{ex}

\section{Conclusion} \label{sec:Conclusion}

We presented a novel adaptive IETI-DP method for solving the 
{\new Kirchhoff--Love plate problem}, {\new where the mid-surface of the plate is described by a} 
planar analysis-suitable $G^1$ multi-patch geometry, which 
{\new is a particular} multi-patch parameterization with possibly extraordinary vertices. The proposed IETI-DP technique not only enables the challenging task of solving fourth-order PDEs {\new on} these multi-patch domains via adaptive refinement in a relatively technically simple way, but also makes it at the same time efficient. This is achieved by decomposing the initial large problem into a small linear problem for the Lagrange multipliers involved, which ensure that the numerical solution is $C^1$-smooth across the inner edges, and into small, local and parallelizable problems for the individual patches to compute the coefficients of the numerical solution, using a dual-primal formulation. The calculated numerical examples demonstrated the great potential of our adaptive IETI-DP method by showing, on the one hand, that the numerical results converge with optimal rates under adaptive refinement, and on the other hand, that the proposed Dirichlet preconditioner for the linear problem of determining the Lagrange multipliers works well. 

There are potential topics for further research that are worth investigating. The proposed adaptive IETI-DP technique could be extended to multi-patch surfaces with possibly extraordinary vertices or to multi-patch volumes with possibly extraordinary vertices and edges. Especially in the multi-patch volume case, which seems to be a rather challenging task and has so far only been limited to uniform refinement and, moreover, only to some initial, specific configurations of the multi-patch volume (see ~\cite{BiJuMa17,BiJuMa20,BiKa19,KaVi22}), the application of an adaptive IETI-DP framework could be beneficial. Furthermore, it could be worth to extend our work to other high-order PDEs such as fourth order problems like the Kirchhoff--Love shell problem~\cite{kiendl-bletzinger-linhard-09} and the Cahn-Hilliard equation~\cite{gomez2008isogeometric}, sixth order problems like the Phase-field crystal equation~\cite{Gomez2012}, the Kirchhoff--Love plate model based on the Mindlin’s gradient elasticity theory~\cite{Niiranen2016} and the gradient-enhanced continuum damage models~\cite{GradientDamageModels} or PDEs with even higher order like the polyharmonic equation~\cite{KaKoVi26} of an arbitrary order~${\new q} \geq 1$. 

\paragraph*{\bf Acknowledgment}

{\new The authors wish to thank the anonymous reviewers for their comments that helped to improve the paper.}
V.~Vitrih and A.~Kosma\v c have been partially supported by the Slovenian Research and Innovation Agency (research program P1-0404 and research project N1-0296). This support is gratefully acknowledged.

%

{\new

%

\appendix
\section{Computational realization}
\label{app:computational-realization}

The inverse matrices appearing in Sections~3 and~4 are only an algebraic
notation. In the implementation, matrices such as
$\ab{K}_{\Delta\Delta}^{(i)\,-1}$,
$\ab{K}_{RR}^{(i)\,-1}$, $\ab{M}_B^{(i)\,-1}$, and
$\ab{S}_{\Pi\Pi}^{-1}$ are never formed explicitly. Every occurrence of an
inverse is realized by solving the corresponding linear system. For the
problem sizes considered in Section~5, the auxiliary local and 
small global (for the primal degrees)
systems
are handled by direct factorizations, whereas the global interface problem
(30), or its preconditioned form (33), is solved iteratively. The
factorizations are performed once for a fixed hierarchical mesh and reused in
all subsequent solves. After an adaptive refinement step, only the matrices
corresponding to the new mesh have to be assembled and factorized again.

\paragraph{Local and global--primal solves}
For each patch $\Omega^{(i)}$, the matrix
$\ab{K}_{\Delta\Delta}^{(i)}$ is factorized once. Hence, expressions of the
form $\ab{K}_{\Delta\Delta}^{(i)\,-1}\ab{v}$ are evaluated by forward and
backward substitutions. The same factorization is reused for several right-hand
sides. For example, instead of forming
$\ab{K}_{\Delta\Delta}^{(i)\,-1}\ab{K}_{\Delta\Pi}^{(i)}$, we compute
$\ab{X}^{(i)}$ from
\[
    \ab{K}_{\Delta\Delta}^{(i)}\ab{X}^{(i)}
    = \ab{K}_{\Delta\Pi}^{(i)},
\]
and use it in the local contribution
\[
    \ab{S}_{\Pi\Pi}^{(i)}
    =
    \ab{N}_{\Xi\Pi}^{(i)T}
    \left(
        \ab{K}_{\Pi\Pi}^{(i)}
        -\ab{K}_{\Pi\Delta}^{(i)}\ab{X}^{(i)}
    \right)
    \ab{N}_{\Xi\Pi}^{(i)}.
\]
Analogous backsolves are used in the construction of $\ab{B}_{\Pi}$,
$\widetilde{\ab{f}}_{\Pi}$, $\ab{F}$, and $\ab{r}$. Since
$\ab{K}_{\Delta\Pi}^{(i)}$ has only a small number of columns compared with
the number of patch unknowns, this requires only a small number of right-hand
sides and does not amount to a large dense matrix--matrix operation.

The global primal matrix
$\ab{S}_{\Pi\Pi}=\sum_{i\in \mathcal{I}_\Omega} \ab{S}_{\Pi\Pi}^{(i)}$ has six primal degrees of
freedom per inner vertex. For the configurations considered here, it is
factorized once with a sparse direct solver and the factorization is reused
whenever $\ab{S}_{\Pi\Pi}^{-1}$ occurs. Similarly, the action of $\ab{F}$ in
\eqref{eq:spdMatrixF} can be evaluated without forming inverse matrices: one applies the
constraint matrices, performs independent patch-local backsolves with
$\ab{K}_{\Delta\Delta}^{(i)}$, performs one 
global--primal backsolve with
$\ab{S}_{\Pi\Pi}$, and multiplies by the corresponding transpose matrices.
For very large numbers of patches and/or for a very high number of degrees of freedom per patches, the direct solving (as used in the examples) could instead be replaced
by a suitably preconditioned iterative or multigrid techniques developed for isogeometric discretizations.
For tensor-product spline spaces, examples include symbol-based multigrid
methods~\cite{DoGaMaSeSp17}, robust multigrid based on stable splittings of
spline spaces~\cite{HoTa17MG}, and $p$-multigrid methods~\cite{TiMoGoVu20}.
For adaptively refined hierarchical spline spaces, geometric multigrid methods
for HB- and THB-splines were developed in~\cite{HoJuKiZu16}, while local
multigrid solvers for adaptive hierarchical spline spaces were proposed
in~\cite{HoMiSp22}. 

\paragraph{Dirichlet preconditioner}
The local Schur complement entering the preconditioner is
\[
    \ab{S}_{EE}^{(i)}
    =
    \ab{K}_{EE}^{(i)}
    -\ab{K}_{ER}^{(i)}
     \ab{K}_{RR}^{(i)\,-1}
     \ab{K}_{RE}^{(i)}.
\]
Again, no explicit inverse is needed. After factorizing
$\ab{K}_{RR}^{(i)}$ once, its action can be computed by
\[
    \ab{S}_{EE}^{(i)}\ab{v}
    =
    \ab{K}_{EE}^{(i)}\ab{v}
    -\ab{K}_{ER}^{(i)}
     \left[\ab{K}_{RR}^{(i)\,-1}
     \left(\ab{K}_{RE}^{(i)}\ab{v}\right)\right].
\]
Thus, one application of the Dirichlet preconditioner consists of local
matrix--vector products and local backsolves, which can be carried out patch
by patch and in parallel. The boundary system in Section~4.2 is treated in the
same way: the blocks $\MM^{(i)}$ are factorized independently, while the
small reduced system involving $\ab{B}\MM^{-1}\ab{B}^{T}$ is solved
directly.

\paragraph{Null spaces and normalization}
The null-space computations in Sections~4.2 and~4.3 are small local dense
problems, not global null-space computations. A basis $\ab{N}^{(i)}$ in (37)
can be obtained by a QR factorization or an SVD. Likewise, for
an inner vertex of valency $\nu_i$, the reduced matrix
$\widehat{\ab{C}}_{\Xi,\Pi}^{(i)}\in
\mathbb{R}^{6(\nu_i-1)\times 6\nu_i}$ has a six-dimensional kernel, and the
basis $\widehat{\ab{N}}_{\Xi,\Pi}^{(i)}\in
\mathbb{R}^{6\nu_i\times 6}$ is computed locally by QR or SVD, and then
assembled into $\ab{N}_{\Xi,\Pi}$. Since the vertex valencies are small, this
cost is negligible compared with the factorization of the patch matrices.

The normalization in (23) is also a one-time setup operation. It may be
computed exactly as written in the paper from the eigendecomposition of
$\ab{C}_{\Gamma,E}\ab{C}_{\Gamma,E}^{T}$. Alternatively, a QR
factorization can be used to obtain an equivalent normalization without
forming an explicit inverse. The small full-rank systems involving
$\ab{D}_{2}^{(i)}$ in Section~4.4 are treated analogously by a local direct
factorization.

\paragraph{Overall cost}
It is useful to distinguish a setup phase from the preconditioned conjugate gradient (PCG) iteration.
The setup contains matrix assembly, patch-local and 
global--primal factorizations,
and the small edge- and vertex-local QR/SVD or normalization problems. During
the PCG solution of (33), the dominant operations are sparse
matrix--vector products, independent local backsolves with already factorized
matrices, and one 
global--primal backsolve. Hence, the inverse symbols and products in
the derivation do not imply explicit matrix inversions or large dense
matrix--matrix multiplications. For the moderate-size decompositions studied
here, direct methods are used for these auxiliary problems and PCG is reserved
for the global multiplier system. Iterative or multilevel auxiliary solves
would become relevant for substantially larger decompositions and/or for a very high number of degrees of freedom per patches.
}




\begin{thebibliography}{10}

\bibitem{AnBuCo20}
P.~Antolin, A.~Buffa, and L.~Coradello.
\newblock A hierarchical approach to the a posteriori error estimation of
  isogeometric {K}irchhoff plates and {K}irchhoff–{L}ove shells.
\newblock {\em Comput. Methods Appl. Mech. Engrg.}, 363:112919, 2020.

\bibitem{ArReKlSi23}
J.~Arf, M.~Reichle, S.~Klinkel, and B.~Simeon.
\newblock Scaled boundary isogeometric analysis with ${C}^1$ coupling for
  {K}irchhoff plate theory.
\newblock {\em Comput. Methods Appl. Mech. Engrg.}, 415:116198, 2023.

\bibitem{BaSm93}
R.~E. Bank and R.~K. Smith.
\newblock A posteriori error estimates based on hierarchical bases.
\newblock {\em SIAM J. Numer. Anal.}, 30(4):921--935, 1993.

\bibitem{VePaScWiZa14}
L.~Beir{\~a}o~da Veiga, L.~F. Pavarino, S.~Scacchi, O.~B. Widlund, and
  S.~Zampini.
\newblock Isogeometric {BDDC} preconditioners with deluxe scaling.
\newblock {\em SIAM J. Sci. Comput.}, 36(3):1118--1139, 2014.

\bibitem{VePaScWiZa17}
L.~Beir{\~a}o~da Veiga, L.~F. Pavarino, S.~Scacchi, O.~B. Widlund, and
  S.~Zampini.
\newblock Adaptive selection of primal constraints for isogeometric {BDDC}
  deluxe preconditioners.
\newblock {\em SIAM J. Sci. Comput.}, 39(1):281--302, 2017.

\bibitem{BeLoSaTa23}
A.~Benvenuti, G.~Loli, G.~Sangalli, and T.~Takacs.
\newblock Isogeometric multi-patch ${C}^1$-mortar coupling for the biharmonic
  equation.
\newblock {\em arXiv}, 2303.07255, 2023.

\bibitem{BiJuMa17}
K.~Birner, B.~J\"uttler, and A.~Mantzaflaris.
\newblock Bases and dimensions of ${C}^1$-smooth isogeometric splines on
  volumetric two-patch domains.
\newblock {\em Graphical Models}, 99:46 -- 56, 2018.

\bibitem{BiJuMa20}
K.~Birner, B.~J\"uttler, and A.~Mantzaflaris.
\newblock Approximation power of ${C}^1$-smooth isogeometric splines on
  volumetric two-patch domains.
\newblock In {\em Isogeometric Analysis and Applications 2018}, pages 27--38.
  Springer, LNCSE, 2021.

\bibitem{BiKa19}
K.~Birner and M.~Kapl.
\newblock The space of ${C}^1$-smooth isogeometric spline functions on
  trilinearly parameterized volumetric two-patch domains.
\newblock {\em Comput. Aided Geom. Des.}, 70:16--30, 2019.

\bibitem{Bouclier2017}
R.~Bouclier, J.~C. Passieux, and M.~Sala{\"{u}}n.
\newblock {Development of a new, more regular, mortar method for the coupling
  of NURBS subdomains within a NURBS patch: Application to a non-intrusive
  local enrichment of NURBS patches}.
\newblock {\em Comput. Methods Appl. Mech. Engrg.}, 316:123--150, 2017.

\bibitem{BrFaGiKaVa24}
C.~Bracco, A.~Farahat, C.~Giannelli, M.~Kapl, and R.~V\'azquez.
\newblock Adaptive methods with ${C}^1$ splines for multi-patch surfaces and
  shells.
\newblock {\em Comput. Methods Appl. Mech. Engrg.}, 431:117287, 2024.

\bibitem{BrGiKaVa20}
C.~Bracco, C.~Giannelli, M.~Kapl, and R.~V\'azquez.
\newblock Isogeometric analysis with ${C}^1$ hierarchical functions on planar
  two-patch geometries.
\newblock {\em Comput. Math. Appl.}, 80(11):2538--2562, 2020.

\bibitem{BrGiKaVa23}
C.~Bracco, C.~Giannelli, M.~Kapl, and R.~V{\'a}zquez.
\newblock Adaptive isogeometric methods with ${C}^1$ (truncated) hierarchical
  splines on planar multi-patch domains.
\newblock {\em Math. Models Methods Appl. Sci.}, 33(9):1829--1874, 2023.

\bibitem{BrGiReToVa23}
C.~Bracco, C.~Giannelli, A.~Reali, M.~Torre, and R.~V{\'a}zquez.
\newblock Adaptive isogeometric phase-field modeling of the {C}ahn–{H}illiard
  equation: {S}uitably graded hierarchical refinement and coarsening on
  multi-patch geometries.
\newblock {\em Comput. Methods Appl. Mech. Engrg.}, 417:116355, 2023.

\bibitem{BuGaGiPrVa22}
A.~Buffa, G.~Gantner, C.~Giannelli, D.~Praetorius, and R.~V\'azquez.
\newblock Mathematical foundations of adaptive isogeometric analysis.
\newblock {\em Arch. Comput. Methods Eng.}, 29:4479--4555, 2022.

\bibitem{BuGi16}
A.~Buffa and C.~Giannelli.
\newblock Adaptive isogeometric methods with hierarchical splines: error
  estimator and convergence.
\newblock {\em Math. Models Methods Appl. Sci.}, 26(1):1--25, 2016.

\bibitem{CaWeToLiHuKiZh20}
H.~Casquero, X.~Wei, D.~Toshniwal, A.~Li, T.J.R Hughes, J.~Kiendl, and Y.~J.
  Zhang.
\newblock Seamless integration of design and {K}irchhoff–{L}ove shell
  analysis using analysis-suitable unstructured {T}-splines.
\newblock {\em Comput. Methods Appl. Mech. Engrg.}, 360:112765, 2020.

\bibitem{ChAnRa18}
C.~L. Chan, C.~Anitescu, and T.~Rabczuk.
\newblock Isogeometric analysis with strong multipatch {C}$^1$-coupling.
\newblock {\em Comput. Aided Geom. Design}, 62:294--310, 2018.

\bibitem{ChAnRa19}
C.~L. Chan, C.~Anitescu, and T.~Rabczuk.
\newblock Strong multipatch {C}$^1$-coupling for isogeometric analysis on {2D}
  and {3D} domains.
\newblock {\em Comput. Methods Appl. Mech. Engrg.}, 357:112599, 2019.

\bibitem{CoSaTa16}
A.~Collin, G.~Sangalli, and T.~Takacs.
\newblock Analysis-suitable {G}$^1$ multi-patch parametrizations for {C}$^1$
  isogeometric spaces.
\newblock {\em Comput. Aided Geom. Des.}, 47:93 -- 113, 2016.

\bibitem{CoAnVaBu20}
L.~Coradello, P.~Antolin, R.~V\'azquez, and A.~Buffa.
\newblock Adaptive isogeometric analysis on two-dimensional trimmed domains
  based on a hierarchical approach.
\newblock {\em Comput. Methods Appl. Mech. Engrg.}, 364:112925, 2020.

\bibitem{DiSchWoHe19}
M.~Dittmann, S.~Schu{\ss}, B.~Wohlmuth, and C.~Hesch.
\newblock Weak ${C}^n$ coupling for multipatch isogeometric analysis in solid
  mechanics.
\newblock {\em Int. J. Numer. Methods Eng.}, 118(11):678--699, 2019.

\bibitem{Do03}
C.~R. Dohrmann.
\newblock A preconditioner for substructuring based on constrained energy
  minimization.
\newblock {\em SIAM J. Sci. Comput.}, 25(1):246--258, 2003.

\bibitem{DoGaMaSeSp17}
M.~Donatelli, C.~Garoni, C.~Manni, S.~Serra-Capizzano, and H.~Speleers.
\newblock Symbol-based multigrid methods for galerkin {B}-spline isogeometric
  analysis.
\newblock {\em SIAM Journal on Numerical Analysis}, 55(1):31--62, 2017.

\bibitem{Do96}
W.~D{\"o}rfler.
\newblock A convergent adaptive algorithm for {Poisson}'s equation.
\newblock {\em SIAM J. Numer. Anal.}, 33(3):1106--1124, 1996.

\bibitem{ErGaPr20}
C.~Erath, G.~Gantner, and D.~Praetorius.
\newblock Optimal convergence behavior of adaptive fem driven by simple
  (h-h/2)-type error estimators.
\newblock {\em Comput. Math. Appl.}, 79(3):623--642, 2020.

\bibitem{FaJuKaTa22}
A.~Farahat, B.~J\"uttler, M.~Kapl, and T.~Takacs.
\newblock Isogeometric analysis with ${C}^1$-smooth functions over multi-patch
  surfaces.
\newblock {\em Comput. Methods Appl. Mech. Engrg.}, 403:115706, 2023.

\bibitem{FaKaKoVi24}
A.~Farahat, M.~Kapl, A.~Kosma\v{c}, and V.~Vitrih.
\newblock A locally based construction of analysis-suitable ${G}^1$ multi-patch
  spline surfaces.
\newblock {\em Comput. Math. Appl.}, 168:46--57, 2024.

\bibitem{FarhatFETI-DP}
C.~Farhat, M.~Lesoinne, P.~LeTallec, K.~Pierson, and D.~Rixen.
\newblock {FETI-DP}: a dual–primal unified {FETI} method—part {I}: {A}
  faster alternative to the two-level {FETI} method.
\newblock {\em Int. J. Numer. Methods Eng.}, 50(7):1523--1544, 2001.

\bibitem{FaMa98}
C.~Farhat and J.~Mandel.
\newblock The two-level {FETI} method for static and dynamic plate problems
  {Part I}: {An} optimal iterative solver for biharmonic systems.
\newblock {\em Comput. Methods Appl. Mech. Engrg.}, 155(1):129--151, 1998.

\bibitem{Farhat1991}
C.~Farhat and F.-X. Roux.
\newblock A method of finite element tearing and interconnecting and its
  parallel solution algorithm.
\newblock {\em Int. J. Numer. Methods Eng.}, 32(6):1205--1227, 1991.

\bibitem{GiJuSp2012}
C.~Giannelli, B.~J{\"u}ttler, and H.~Speleers.
\newblock {THB}-splines: The truncated basis for hierarchical splines.
\newblock {\em Computer Aided Geometric Design}, 29(7):485--498, 2012.

\bibitem{gomez2008isogeometric}
H.~Gomez, V.~M Calo, Y.~Bazilevs, and T.~J.R. Hughes.
\newblock Isogeometric analysis of the {Cahn--Hilliard} phase-field model.
\newblock {\em Comput. Methods Appl. Mech. Engrg.}, 197(49):4333--4352, 2008.

\bibitem{Gomez2012}
H.~Gomez and X.~Nogueira.
\newblock An unconditionally energy-stable method for the phase field crystal
  equation.
\newblock {\em Comput. Methods Appl. Mech. Engrg.}, 249 -- 252:52 -- 61, 2012.

\bibitem{Gr92}
P.~Grisvard.
\newblock {\em Singularities in boundary value problems}, volume~22 of {\em
  Rech. Math. Appl.}
\newblock Paris: Masson; Berlin: Springer-Verlag, 1992.

\bibitem{HaTaTa26}
F.~Hasanova, S.~Takacs, and T.~Takacs.
\newblock Robust approximation error estimates for analysis-suitable {$G^1$}
  isogeometric multi-patch discretizations.
\newblock {\em arXiv}, 2605.13270, 2026.

\bibitem{HoLa17}
C.~Hofer and U.~Langer.
\newblock Dual-primal isogeometric tearing and interconnecting solvers for
  multipatch {dG}-{Iga} equations.
\newblock {\em Comput. Methods Appl. Mech. Eng.}, 316:2--21, 2017.

\bibitem{HoJuKiZu16}
C.~Hofreither, B.~J{\"u}ttler, G.~Kiss, and W.~Zulehner.
\newblock Multigrid methods for isogeometric analysis with {THB}-splines.
\newblock {\em Computer Methods in Applied Mechanics and Engineering},
  308:96--112, 2016.

\bibitem{HoMiSp22}
C.~Hofreither, L.~Mitter, and H.~Speleers.
\newblock Local multigrid solvers for adaptive isogeometric analysis in
  hierarchical spline spaces.
\newblock {\em IMA Journal of Numerical Analysis}, 42(3):2429--2458, 2022.

\bibitem{HoTa17MG}
C.~Hofreither and S.~Takacs.
\newblock Robust multigrid for isogeometric analysis based on stable splittings
  of spline spaces.
\newblock {\em SIAM Journal on Numerical Analysis}, 55(4):2004--2024, 2017.

\bibitem{HuCoBa04}
T.~J.~R. Hughes, J.~A. Cottrell, and Y.~Bazilevs.
\newblock Isogeometric analysis: {CAD}, finite elements, {NURBS}, exact
  geometry and mesh refinement.
\newblock {\em Comput. Methods Appl. Mech. Engrg.}, 194(39-41):4135--4195,
  2005.

\bibitem{KaKoVi24}
M.~Kapl, A.~Kosma\v{c}, and V.~Vitrih.
\newblock Isogeometric collocation for solving the biharmonic equation over
  planar multi-patch domains.
\newblock {\em Comput. Methods Appl. Mech. Engrg.}, 424:116882, 2024.

\bibitem{KaKoVi24b}
M.~Kapl, A.~Kosma\v{c}, and V.~Vitrih.
\newblock A ${C}^s$-smooth mixed degree and regularity isogeometric spline
  space over planar multi-patch domains.
\newblock {\em J. Comput. Appl. Math.}, 473:116836, 2026.

\bibitem{KaKoVi24c}
M.~Kapl, A.~Kosma\v{c}, and V.~Vitrih.
\newblock Isogeometric collocation with smooth mixed degree splines over planar
  multi-patch domains.
\newblock {\em J. Comput. Appl. Math.}, 210:89--112, 2026.

\bibitem{KaKoVi26}
M.~Kapl, A.~Kosma\v{c}, and V.~Vitrih.
\newblock An isogeometric tearing and interconnecting ({IETI}) method for
  solving high order partial differential equations over planar multi-patch
  geometries.
\newblock {\em Comput. Methods Appl. Mech. Engrg.}, 452:118769, 2026.

\bibitem{KaSaTa17b}
M.~Kapl, G.~Sangalli, and T.~Takacs.
\newblock Construction of analysis-suitable {G}$^1$ planar multi-patch
  parameterizations.
\newblock {\em Comput.-Aided Des.}, 97:41--55, 2018.

\bibitem{KaSaTa19b}
M.~Kapl, G.~Sangalli, and T.~Takacs.
\newblock Isogeometric analysis with {C}$^{1}$ functions on unstructured
  quadrilateral meshes.
\newblock {\em The SMAI journal of computational mathematics}, 5:67--86, 2019.

\bibitem{KaSaTa19a}
M.~Kapl, G.~Sangalli, and T.~Takacs.
\newblock An isogeometric {C}$^{1}$ subspace on unstructured multi-patch planar
  domains.
\newblock {\em Comput. Aided Geom. Des.}, 69:55--75, 2019.

\bibitem{KaVi17b}
M.~Kapl and V.~Vitrih.
\newblock Space of {C}$^2$-smooth geometrically continuous isogeometric
  functions on planar multi-patch geometries: {D}imension and numerical
  experiments.
\newblock {\em Comput. Math. Appl.}, 73(10):2319--2338, 2017.

\bibitem{KaVi17c}
M.~Kapl and V.~Vitrih.
\newblock Dimension and basis construction for ${C}^{2}$-smooth isogeometric
  spline spaces over bilinear-like ${G}^{2}$ two-patch parameterizations.
\newblock {\em J. Comput. Appl. Math.}, 335:289--311, 2018.

\bibitem{KaVi19a}
M.~Kapl and V.~Vitrih.
\newblock Solving the triharmonic equation over multi-patch planar domains
  using isogeometric analysis.
\newblock {\em J. Comput. Appl. Math.}, 358:385--404, 2019.

\bibitem{KaVi20}
M.~Kapl and V.~Vitrih.
\newblock Isogeometric collocation on planar multi-patch domains.
\newblock {\em Comput. Methods Appl. Mech. Engrg.}, 360:112684, 2020.

\bibitem{KaVi20b}
M.~Kapl and V.~Vitrih.
\newblock ${C}^s$-smooth isogeometric spline spaces over planar multi-patch
  parameterizations.
\newblock {\em Adv. Comput. Math.}, 47:47, 2021.

\bibitem{KaVi22}
M.~Kapl and V.~Vitrih.
\newblock ${C}^1$ isogeometric spline space for trilinearly parameterized
  multi-patch volumes.
\newblock {\em Comput. Math. Appl.}, 117:53--68, 2022.

\bibitem{KaPe17}
K.~Kar{\v c}iauskas and J.~Peters.
\newblock Refinable ${G}^1$ functions on ${G}^1$ free-form surfaces.
\newblock {\em Comput. Aided Geom. Des.}, 54:61--73, 2017.

\bibitem{KaPe18}
K.~Kar{\v c}iauskas and J.~Peters.
\newblock Refinable bi-quartics for design and analysis.
\newblock {\em Comput.-Aided Des.}, 102:204--214, 2018.

\bibitem{kiendl-bletzinger-linhard-09}
J.~Kiendl, K.-U. Bletzinger, J.~Linhard, and R.~W{\"u}chner.
\newblock Isogeometric shell analysis with {K}irchhoff-{L}ove elements.
\newblock {\em Comput. Methods Appl. Mech. Engrg.}, 198(49):3902--3914, 2009.

\bibitem{Klawonn2000}
A.~Klawonn and O.~B. Widlund.
\newblock A domain decomposition method with {L}agrange multipliers and inexact
  solvers for linear elasticity.
\newblock {\em SIAM Journal on Scientific Computing}, 22(4):1199--1219, 2000.

\bibitem{Klawonn2006}
A.~Klawonn and O.~B. Widlund.
\newblock Dual-primal {FETI} methods for linear elasticity.
\newblock {\em Communications on Pure and Applied Mathematics},
  59(11):1523--1572, 2006.

\bibitem{KlPeSaJu12}
S.~K. Kleiss, C.~Pechstein, B.~J\"uttler, and S.~Tomar.
\newblock {IETI} -- isogeometric tearing and interconnecting.
\newblock {\em Comput. Methods Appl. Mech. Engrg.}, 247-248:201--215, 2012.

\bibitem{Kr97}
R.~Kraft.
\newblock Adaptive and linearly independent multilevel {B}-splines.
\newblock In {\em Surface fitting and multiresolution methods. Vol. 2 of the
  proceedings of the 3rd international conference on Curves and surfaces, held
  in Chamonix-Mont-Blanc, France, June 27--July 3, 1996}, pages 209--218.
  Nashville, TN: Vanderbilt University Press, 1997.

\bibitem{MaDo03}
J.~Mandel and C.~R. Dohrmann.
\newblock Convergence of a balancing domain decomposition by constraints and
  energy minimization.
\newblock {\em Numerical Linear Algebra with Applications}, 10(7):639--659,
  2003.

\bibitem{MaDoTe05}
J.~Mandel, C.~R. Dohrmann, and R.~Tezaur.
\newblock An algebraic theory for primal and dual substructuring methods by
  constraints.
\newblock {\em Appl. Numer. Math.}, 54(2):167--193, 2005.

\bibitem{MaMaMo2024}
M.~Marsala, A.~Mantzaflaris, and B.~Mourrain.
\newblock ${G}^1$ spline functions for point cloud fitting.
\newblock {\em Applied Mathematics and Computation}, 460:128279, 2024.

\bibitem{MiZoScBoTh21}
D.~Miao, Z.~Zou, M.~A. Scott, M.~J. Borden, and D.~C. Thomas.
\newblock Isogeometric {B{\'e}zier} dual mortaring: the {Kirchhoff}-{Love}
  shell problem.
\newblock {\em Comput. Methods Appl. Mech. Eng.}, 382:113873, 2021.

\bibitem{MoSaSchTaTa23}
M.~Montardini, G.~Sangalli, R.~Schneckenleitner, S.~Takacs, and M.~Tani.
\newblock A {IETI}-{DP} method for discontinuous {Galerkin} discretizations in
  isogeometric analysis with inexact local solvers.
\newblock {\em Math. Models Methods Appl. Sci.}, 33(10):2085--2111, 2023.

\bibitem{mourrain2015geometrically}
B.~Mourrain, R.~Vidunas, and N.~Villamizar.
\newblock Dimension and bases for geometrically continuous splines on surfaces
  of arbitrary topology.
\newblock {\em Comput. Aided Geom. Des.}, 45:108 -- 133, 2016.

\bibitem{NgPe16}
T.~Nguyen and J.~Peters.
\newblock Refinable ${C}^{1}$ spline elements for irregular quad layout.
\newblock {\em Comput. Aided Geom. Des.}, 43:123--130, 2016.

\bibitem{Niiranen2016}
J.~Niiranen, J.~Kiendl, A.~H. Niemi, and A.~Reali.
\newblock Isogeometric analysis for sixth-order boundary value problems of
  gradient-elastic {K}irchhoff plates.
\newblock {\em Comput. Methods Appl. Mech. Engrg.}, 316:328--348, 2017.

\bibitem{ReArSiKl23}
M.~F.~M. Reichle, J.~Arf, B.~Simeon, and S.~Klinkel.
\newblock {S}mooth multi-patch scaled boundary isogeometric analysis for
  {K}irchhoff–{L}ove shells.
\newblock {\em Meccanica}, 58(8):1693--1716, 2023.

\bibitem{Re97}
U.~Reif.
\newblock A refinable space of smooth spline surfaces of arbitrary topological
  genus.
\newblock {\em Journal of Approximation Theory}, 90(2):174--199, 1997.

\bibitem{IETI_LowRank2024}
T.-C. Riemer, A.~B{\"u}nger, and M.~Stoll.
\newblock {IETI}-based low-rank method for {PDE}-constrained optimization.
\newblock {\em SIAM Journal on Scientific Computing}, 47(3):1758--1782, 2025.

\bibitem{WeightedCollocation2013}
D.~Schillinger, J.~A. Evans, A.~Reali, M.~A. Scott, and T.~J.~R. Hughes.
\newblock Isogeometric collocation: {C}ost comparison with {G}alerkin methods
  and extension to adaptive hierarchical {NURBS} discretizations.
\newblock {\em Computer Methods in Applied Mechanics and Engineering},
  267:170--232, 2013.

\bibitem{SchDiWoKlHe19}
S.~Schu{\ss}, M.~Dittmann, B.~Wohlmuth, S.~Klinkel, and C.~Hesch.
\newblock Multi-patch isogeometric analysis for {K}irchhoff–{L}ove shell
  elements.
\newblock {\em Comput. Methods Appl. Mech. Engrg.}, 349:91--116, 2019.

\bibitem{ScSiEv13}
M.A. Scott, R.N. Simpson, J.A. Evans, S.~Lipton, S.P.A. Bordas, T.J.R. Hughes,
  and T.W. Sederberg.
\newblock Isogeometric boundary element analysis using unstructured
  {T}-splines.
\newblock {\em Comp. Methods Appl. Mech. Engrg.}, 254:197--221, 2013.

\bibitem{SoTa23}
J.~Sogn and S.~Takacs.
\newblock Stable discretizations and {IETI}-{DP} solvers for the {Stokes}
  system in multi-patch {IgA}.
\newblock {\em ESAIM, Math. Model. Numer. Anal.}, 57(2):921--952, 2023.

\bibitem{SoglTakacs_IETI_Elasticity}
J.~Sogn and S.~Takacs.
\newblock Isogeometric tearing and interconnecting solvers for linear
  elasticity in multi-patch isogeometric analysis with theory for two
  dimensional domains.
\newblock {\em Comput. Methods Appl. Mech. Engrg.}, 418:116482, 2024.

\bibitem{Ta25}
S.~Takacs.
\newblock An isogeometric tearing and interconnecting method for conforming
  discretizations of the biharmonic problem.
\newblock {\em arXiv}, 2511.05247, 2025.

\bibitem{TiMoGoVu20}
R.~Tielen, M.~M{\"o}ller, D.~G{\"o}ddeke, and C.~Vuik.
\newblock $p$-multigrid methods and their comparison to $h$-multigrid methods
  within isogeometric analysis.
\newblock {\em Computer Methods in Applied Mechanics and Engineering},
  372:113347, 2020.

\bibitem{ToSpHu17}
D.~Toshniwal, H.~Speleers, and T.~J.~R. Hughes.
\newblock Smooth cubic spline spaces on unstructured quadrilateral meshes with
  particular emphasis on extraordinary points: Geometric design and
  isogeometric analysis considerations.
\newblock {\em Comput. Methods Appl. Mech. Engrg.}, 327:411--458, 2017.

\bibitem{GradientDamageModels}
C.~V. Verhoosel, M.~A. Scott, T.~J.~R. Hughes, and R.~de~Borst.
\newblock An isogeometric analysis approach to gradient damage models.
\newblock {\em Internat. J. Numer. Methods Engrg.}, 86(1):115--134, 2011.

\bibitem{We22}
X.~Wei.
\newblock {THU}-splines: {H}ighly localized refinement on smooth unstructured
  splines.
\newblock In Carla Manni and Hendrik Speleers, editors, {\em Geometric
  Challenges in Isogeometric Analysis}, pages 305--332, Cham, 2022. Springer
  International Publishing.

\bibitem{WeLiQiHuZhCa22}
X.~Wei, X~Li, K.~Qian, T.~J.~R Hughes, Y.~J Zhang, and H.~Casquero.
\newblock Analysis-suitable unstructured {T}-splines: {M}ultiple extraordinary
  points per face.
\newblock {\em Comput. Methods Appl. Mech. Engrg.}, 391:114494, 2022.

\bibitem{WeZhLiHu17}
X.~Wei, Y.~Zhang, L.~Liu, and T.~J.~R. Hughes.
\newblock Truncated {T}-splines: {F}undamentals and methods.
\newblock {\em Comput. Methods Appl. Mech. Engrg.}, 316:349--372, 2017.

\bibitem{WiZaScPa21}
O.~B. Widlund, S.~Zampini, S.~Scacchi, and L.~F. Pavarino.
\newblock Block {FETI}-{DP}/{BDDC} preconditioners for mixed isogeometric
  discretizations of three-dimensional almost incompressible elasticity.
\newblock {\em Math. Comput.}, 90(330):1773--1797, 2021.

\bibitem{WeFaLiWeCa23}
Wen. Z., Md.~S. Faruque, X.~Li, X.~Wei, and H.~Casquero.
\newblock Isogeometric analysis using {G}-spline surfaces with arbitrary
  unstructured quadrilateral layout.
\newblock {\em Comput. Methods Appl. Mech. Engrg.}, 408:115965, 2023.

\end{thebibliography}
%

\end{document}